%% file: orga.tex
\documentclass[10pt]{report}

\usepackage{german}

\begin{document}
\input{tita}

\tableofcontents
\newtheorem{satz}{Satz}[section]
\newtheorem{lem}[satz]{Lemma}
\newtheorem{kor}[satz]{Korollar}
\newtheorem{pro}[satz]{Proposition}
\input{intra}

\input{tast}

\input{last}

\input{tafi}

\input{rang}

\input{ausb}

\input{lita}
\end{document}

%% file: tita.tex
\begin{titlepage}
\begin{center} \Large RUHR-UNIVERSIT"AT BOCHUM    \\
FAKULT"AT F"UR MATHEMATIK \end{center}
\vspace{5.5cm}
\begin{center} \huge TOPOLOGISCHE UND ALGEBRAISCHE FILTER \end{center}
\vspace{5cm}
\begin{center} \normalsize Diplomarbeit von Holger Brenner \\
Betreuer: Prof. Dr. U. Storch \\
25.10.1994 
\end{center}
\end{titlepage}

%% file: intra.tex
\newpage
\section*{Einleitung}
\addcontentsline{toc}{section}{Einleitung}
In dieser Arbeit wird ein recht allgemein gefa\3ter Begriff von einem
Filter
vorgestellt und diskutiert, der es erlaubt, in jedem kommutativen Monoid
von Filtern zu sprechen. Er beinhaltet unter anderem die aus der Topologie
bekannten Filter, aber auch die in der kommutativen Algebra auftretenden
saturierten multiplikativen Systeme, die beispielsweise bei der
Konstruktion von Bruch\-ringen eine Rolle spielen. Der Inhalt der Arbeit
sei im folgenden kurz skizziert.
\par
Im ersten Kapitel werden nach der allgemeinen Definition und einigen
Beispielen die Filter in kommutativen Ringen untersucht, insbesondere ihre
Beziehung zu Primelementen und Primidealen und der Zusammenhang mit
Irreduzibilit"ats- und Teilerfremdheits\-eigenschaften. Diese "Uberlegungen
wer\-den im vierten Kapitel wieder aufgegriffen und durch geometrische 
Deutungen im Spektrum des Ringes erg"anzt.
\par
Im zweiten Kapitel wird in weitgehender Analogie zur Definition des
Spektrums eines kommutativen Ringes die Menge der Filter eines Monoids
mit einer Topologie versehen und dieser Konstruktionsproze\3 als Funktor
nachgewiesen. Die Tatsache, da\3 es sich hierbei, anders als beim Spektrum,
nicht nur um einen kontravarianten, sondern auch um einen kovarianten
Funktor handelt, erlaubt es, die Fil\-ter zwei\-er durch einen Mor\-phis\-mus
ver\-bun\-dener Mono\-ide wech\-sel\-sei\-tig in Be\-zieh\-ung zu set\-zen und zu dem Begriff
des Fixfilters zu gelangen. Er wird im weiteren eine wichtige Rolle spielen.
Die Analogie des Filtrums eines kommutativen Ringes zu seinem Spektrum wird
im vierten Kapitel weitergef"uhrt, indem es mit einer Strukturgarbe versehen
und nachgewiesen wird, da\3 dieser beringte Raum in der Kate\-gorie aller
beringten R"aume die gleiche uni\-ver\-sel\-le Ei\-gen\-schaft hat wie das Spektrum
in der Kate\-gorie der lokal bering\-ten R"aume. Da\3 die nat"ur\-lich
vorkom\-menden beringten R"aume alle\-samt lokal beringt sind, h"angt damit
zusammen, da\3 die Filtrums\-r"aume als solche vom geo\-metri\-schen
Standpunkt aus be\-trach\-tet abseits stehen und ziem\-lich unver\-traute
Eigen\-schaften haben. Die R"aume, die sich 
als Filtrum eines kommu\-ta\-tiven Mono\-ids er\-ge\-ben k"on\-nen, wer\-den topo\-lo\-gisch
cha\-rak\-ter\-i\-siert.
\par
Das dritte Kapitel soll zeigen, da\3 diese R"aume dennoch eini\-ge
be\-mer\-kens\-wer\-te Ei\-gen\-schaf\-ten
haben. Dort wird das Filtrum eines
topologi\-schen Raumes untersucht, das in nat"urlicher Weise eine Ausdehnung
des Ausgangs\-raumes ist und das viele bekannte topologische Erweiterungen
als Unterr"aume bein\-haltet. Topologische Filter erlauben es somit,
ver\-schie\-dene Kon\-struk\-tio\-nen in einem ein\-heit\-lichen Rahmen abzuhandeln. Eine
weitere interes\-sante Eigen\-schaft des Filtrums eines topologischen Raumes
ergibt sich im Zusammenhang mit Garben. Es wird gezeigt, da\3 die Bildgarben
(auch die h"oheren) auf dem Filtrum im wesentlichen bereits die Information
"uber s"amtliche Bildgarben in verschiedenen R"aumen enth"alt.
\par
Das vierte Kapitel besch"aftigt sich mit beringten R"aumen und versucht, das
Verh"altnis zwischen algebraischen und topologischen Filtern zu kl"aren.
Damit f"ugt es sich ein in das allgemeine Konzept, zu untersuchen, welche
algebraischen Objekte des globalen Schnittringes mit welchen topologischen
Objekten korrespondieren. Dabei werden die hier vorgestellten
"Uberlegungen anhand konkreterer Situationen verdeutlicht, und es ergeben sich
Beziehung\-en zu vollst"andig regul"aren und zu kompakten R"aumen, zu
Nullstellens"atzen, zur Strukturgarbe auf dem Spektrum, zu Affinit"atsfragen,
zu Trennungs\-ei\-gen\-schaften und anderem.
\par
Aus dieser kurzen "Ubersicht d"urfte hervorgehen, da\3 sich diese Arbeit nicht
eindeutig einer mathematischen Disziplin zuordnen l"a\3t. Sie ber"uhrt neben
anderem die Gebiete kommutative  Algebra, Topologie, algebraische Geometrie,
Theorie der beringten R"aume, wobei jeweils im wesentlichen nur
Grundkenntnisse aus diesen Gebieten vorausgesetzt werden. Bei spezielleren
Dingen wird auf Literatur verwiesen. 
\par
Ich danke meinen Eltern f"ur die Unterst"utzung w"ahrend meines Studi\-ums 
und Herrn Prof. Dr. U. Storch f"ur die gute Betreuung der Arbeit und
be\-sonders f"ur die M"oglichkeit, eine Arbeit zu schreiben, bei der ich
eigenen "Uberlegungen nachgehen konnte, ohne an feste Vorgaben gebunden zu
sein. Ferner den Teilnehmern an dem Oberseminar `Algebra und Geometrie'
an der Ruhr-Universit"at Bochum f"ur verschiedene Hinweise und
Anregungen.

%% file: tast.tex
\chapter {Filter in einem kommutativen Monoid}
\section {Definition und Beispiele}
In diesem Abschnitt wird der Begriff des Filters f"ur ein beliebiges
kommuta\-ti\-ves
Monoid definiert und erl"autert. Die Verkn"upfung im Monoid wird multipli\-ka\-tiv
geschrieben mit dem neutralen Element 1. Im mengentheoreti\-sch\-en Kontext wird
stets der Durchschnitt (nicht die Vereinigung) als die
Ver\-kn"upf\-ung behandelt
und in einem kommutativen Ring die Multiplikation. Ist in einem Monoid von
einem Nullelement 0 die Rede, so ist ein (dadurch ein\-deutig bestimmtes)
Element gemeint, das jedes Element annulliert. Dabei wird $0 = 1$ zugelassen.
\par\bigskip\noindent
{\bf Definition} Sei $M$ ein kom\-muta\-tives Monoid. Eine Teil\-menge $ F
\subseteq M $ hei\3t Filter, wenn fol\-gende Beding\-ungen erf"ullt sind:
\par\smallskip\noindent
(1) $ 1 \in F $
\par\smallskip\noindent
(2) Mit $ f \in F $ und $ g \in F $ ist stets auch $ f \cdot g \in F $
\par\smallskip\noindent
(3) Ist $ f \in F $ und $g$ ein Teiler von $f$, so ist $ g \in F $ .
\par\bigskip\noindent
Ein Filter ist also eine multi\-pli\-kativ ab\-geschlos\-sene, teiler\-stabile
Teil\-menge  
von $M$, die die 1 enth"alt. Da die 1 jedes Ele\-ment aus $M$ teilt, kann
 (1) durch
die For\-derung ersetzt werden, da\3 $F$ nicht leer ist. Ring\-theo\-retisch ist 
ein
Filter das\-selbe wie ein satur\-iertes multi\-pli\-katives System. Be\-sitzt
 das Monoid
ein Null\-ele\-ment 0, so ist das ganze Monoid $M$ der ein\-zige Filter, der
 die Null
ent\-h"alt, da jedes Ele\-ment $f \in M $ wegen $ f \cdot 0 = 0 $ ein Teiler von
Null ist. In diesem Fall hei\3en die Filter, die die Null nicht enthalten,
konsistente Filter, und der gesamte Filter inkonsistent. H"aufig werden bei der
Definition von Filtern in der Topologie oder in der mathematischen Logik nur
die konsistenten Filter zugelassen, alles in allem ist es aber vorteilhaft, auch
den inkonsistenten Filter zuzulassen und in bestimmten Situationen explizit
auszuschlie\3en.
\par \bigskip \noindent
{\bf Beispiel 1}   Ist $N$ eine Menge, so wird gew"ohnlich ein Filter in $N$ als
nichtleere Teilmenge 
der Potenz\-menge von $N$ de\-finiert, die mit zwei Teil\-mengen von $N$ auch 
ihren Durch\-schnitt und mit einer Teil\-menge auch jede gr"o\3ere enth"alt.
Fa\3t man $ M = \wp (N) $ als kommuta\-tives Monoid mit dem Durch\-schnitt als
Ver\-kn"upfung und der ganzen Menge $N$ als neu\-tralem Element auf,
so ist diese
Defini\-tion mit der obigen ver\-tr"aglich. In diesem Fall bedeutet
 n"am\-lich die Teilbarkeit
von $A$ durch $B$, also $ B \cap C =A $ f"ur ein gewisses $C$, die 
Teil\-mengen\-bezieh\-ung
$ A \subseteq B $, aus der umge\-kehrt die Teilbar\-keits\-bezieh\-ung
$ B \cap A = A $ resultiert.
\par \medskip \noindent
{\bf Beispiel 2}   Die Einheitengruppe $ M ^\times $ von $ M $ ist ein Filter.
Sind $ f $ und $ g $
Ein\-heiten
in $ M $ mit den Inversen $ f ^{-1} $ und $ g ^{-1} $ , so ist
$ g ^{-1} \cdot f ^{-1} $ das Inverse von $ f \cdot g $ , also $ f \cdot
g \in M ^\times $ . Ist $ g \cdot a = f $ und $ f \cdot f^{-1} = 1 $ , so ist
auch $ g $ inver\-tierbar mit dem Inversen $ a \cdot f^{-1} $. Da jeder Filter die
1 enth"alt und jede Ein\-heit die $1$ teilt, ist die Ein\-heiten\-gruppe $ M 
^\times $ in jedem Filter ent\-halten. In Beispiel 1 ist die Einheiten\-gruppe
einfach gleich $ \{N\}$.
\par \medskip \noindent
{\bf Beispiel 3}   Sei $M$ ein kommu\-tatives Monoid mit Nullelement.
Dann bilden die
Nicht\-null\-teiler
$ M^\ast $ von $M$ einen Filter. Sind n"amlich $ f $ und $ g $ 
Nicht\-null\-teiler, so folgt aus $ f \cdot g \cdot a = 0 $ zun"achst 
$ g \cdot a = 0 $ 
und daraus $ a = 0 $, also ist $ f \cdot g $ ein Nichtnullteiler. Ist $ g $ 
ein Teiler des Nicht\-null\-teilers $ f $, etwa $ a \cdot g = f $, so
impli\-ziert
$ g \cdot b = 0 $ sofort $ a \cdot g \cdot b = f \cdot b = 0 $, also ist 
$ b = 0 $ und $ g $ ein Nichtnull\-teiler.
\par \medskip \noindent
{\bf Beispiel 4}   Sei $ S \subseteq M $ ein Untermonoid. So bildet die Menge
aller Teiler 
von $S$ einen Filter, und zwar den kleinsten, der $S$ umfa\3t. Hierzu hat man sich 
nur klar zu machen, da\3 die Teil\-bar\-keit eine tran\-sitive Bezieh\-ung ist 
und da\3, wenn $ f $ von $ a $ und $ g $ von $ b $ geteilt wird, dann
$ f \cdot g $
von $ a \cdot b $ geteilt wird.
\par \medskip \noindent
{\bf Beispiel 5}   Sind $ F_i $ ,$ i \in I, $ Filter in $M$, so ist auch der 
Durchschnitt $ \cap 
_{i \in I} F_i $ ein Filter. Man kann daher zu einer beliebigen Teil\-menge 
$ S \subseteq M $ von dem von $S$ erzeugten Filter $ F(S) $ 
sprechen. $F(S)$ besteht aus allen Teilern von Produkten mit Faktoren aus
$S$, wobei man wie "ublich das leere Produkt gleich 1 setzt. Ist $ S = 
\{f\} $ einelementig, so wird hierf"ur $ F(f) $ statt $ F(\{f\}) $
geschrie\-ben. Diese Filter hei\3en Hauptfilter.
$ F(f) $ besteht aus allen Teilern von Potenzen von $ f $.
\par\bigskip
Bevor wir auf die enge Beziehung zwischen Primidealen und Filtern in
kommu\-tativen Ringen eingehen, sei an eine verallgemeinerte Version des
Lemmas von Krull erinnert: Ist $ {\bf a} \subseteq A $ ein Ideal und
$ S \subseteq A $ ein multiplikatives System mit $ {\bf a} \cap S = \emptyset $,
so gibt es ein Primideal ${\bf p} $ mit $ {\bf a} \subseteq {\bf p} $ und
$ {\bf p} \cap S = \emptyset $.
\begin{satz}
Sei $ A $ ein kommutativer Ring. Dann ist das Komple\-ment ein\-er be\-lie\-bi\-gen
Ver\-ein\-i\-gung von Primidealen ein Fil\-ter. Um\-ge\-kehrt ist das
Kom\-ple\-ment eines Filters
darstellbar als Ver\-einigung ge\-eig\-ne\-ter Primideale.
\end{satz}
Beweis. Sei zun"achst $ {\bf p} \subseteq A $ 
ein Primideal, 
so ist nach Defini\-tion 
das Komplement davon ein
multi\-pli\-katives System ist. Die Teiler\-stabi\-lit"at ergibt sich,
da ein Ideal
unter Vielfach\-bildung abgeschlossen ist. F"ur eine beliebige Verein\-igung
ergibt sich die Behauptung aus Beispiel 5. Der leeren Verein\-igung ent\-spricht
hierbei der ganze Filter $ F=A $. Sei umgekehrt ein Filter $ F \subseteq A $
gegeben, und sei $ a \in A, a \not\in F $. Wegen der Teil\-erstabilit"at ist 
damit auch das von $ a $ erzeugte Hauptideal zu $ F $ disjunkt. Dann gibt
es nach dem Lemma von Krull ein Primideal $ {\bf p} $ mit $ a \in {\bf p}$
und $ {\bf p} \cap F = \emptyset $. Die Verein\-igung der so gefundenen
Primideale leistet das Gew"unschte.
\par \medskip \noindent
{\bf Bemerkung} Als Folgerung erh"alt man sofort: In einem kommutativen Ring
gibt es genau dann nur endlich viele Filter, wenn es nur endlich viele 
Prim\-ide\-a\-le gibt. Dagegen kann es in einem solchen Ring unendlich viele
multipli\-ka\-ti\-ve
Systeme geben.
\par\bigskip
Die im Satz angesprochene Darstellung mit Primidealen ist
nat"urlich nicht eindeutig. Beispielsweise ist das Komplement
der Einheitengruppe gleich der Vereinigung aller Primideale, aber auch
gleich der Vereinigung aller maximalen Ideale.
In Spezialf"allen lassen sich die Filter eines Ringes genauer durch
Primelemente bzw. bestimmte Primideale beschreiben. In noetherschen Ringen kann
man sich aufgrund des Krull'schen Hauptideal\-satz\-es auf Primideale der H"ohe
$ \le 1 $ beschr"anken:
\begin{satz}
In einem noetherschen Ring $ A $ ist das Komplement eines
Filters $ F $ Ver\-einigung 
von Primidealen der H"ohe $ \le 1 $.
\end{satz}
Beweis. Zu $ f \not\in F $ w"ahlt man ein minimales Prim\-ober\-ideal $ {\bf p} $
von $ Af $, das zu $ F $ disjunkt ist. Der Haupt\-ideal\-satz besagt nun, da\3
$ {\bf p} $ die ange\-sprochene H"ohen\-eigen\-schaft hat.
\par\bigskip
Die Menge der Filter eines Monoids ist bez"uglich der Inklusion geordnet, mit
der Einheitengruppe als kleinstem und dem gesamten Monoid als gr"o\3tem
Element. Im mengentheoretischen Kontext spielen die unter allen konsisten\-ten,
also von der
gesamten Potenzmenge verschiedenen, Filtern maximalen Elemente eine besondere
Rolle und werden Ultrafilter genannt. Begriff, Existenzbeweis und einige
charakteristische Eigenschaften lassen sich auf beliebige Monoide mit einem
Nullelement "ubertragen. Bei Ringen ergibt sich eine Beziehung zu minimalen
Primidealen.
\par\medskip\noindent
{\bf Definition} Sei $ M $ ein kommutatives Monoid mit einem Nullelement
$ 0 \neq 1 $. Dann hei\3en die
maximalen, konsistenten Filter Ultrafilter.
\par\medskip\noindent
Die Existenz von Ultrafiltern beweisen wir sofort in einer versch"arften
Schach\-telungsform:
\begin{lem}
Sei $ M $ ein kommutatives Monoid, $ S \subset M $ ein multiplikatives System,
$ {\bf a} \subset M $ ein `Pseudoideal', $ {\bf a } $ sei also abge\-schlos\-sen
unter Mul\-ti\-pli\-ka\-tion mit Elementen aus $ M $, 
$ S \cap {\bf a} = \emptyset $, und es sei
$ {\cal Z} = \{ F \subseteq M \mbox{ \rm Filter} :
 S \subseteq F, F \cap {\bf a} = \emptyset \} $. Dann ist $ \cal Z $ induktiv
geordnet und enth"alt insbesondere maximale Elemente. Dabei ist ein Filter
$ F \in \cal Z $ genau dann maximal, wenn gilt: Ist $ g \not\in F $, so gibt es
ein $ n \in {\bf N} $ und ein $ f \in F $ mit $ g^n \cdot f \in {\bf a} $.
\end{lem}
Beweis. Da $ S \cap {\bf a} = \emptyset $ ist $ F(S) \in \cal Z $, also
$ \cal Z \neq \emptyset $. Sei $ F_i $, $i \in I, $ eine Kette in $ \cal Z $. Dann
ist offenbar $ \bigcup _{i \in I} F_i $ ein Filter aus $ \cal Z $ und eine
obere Schranke. Nach dem Lemma von Zorn besitzt daher $ \cal Z $ maximale
Elemente. Den Zusatz beweist man so: Sei $ F \in  \cal Z $ maximal und
angenommen, es gebe ein $ g \not \in F $ und f"ur alle  $ n \in {\bf N} $ und
alle $ f \in F $ sei $ g^n \cdot f \not\in {\bf a} $. Dann w"are
$ F(F \cup \{g\}) $ ein gr"o\3erer Filter in $ \cal Z $. Gebe es n"amlich
einen Teiler $ x $ eines Produktes $ g^n \cdot f $ mit $ x \in {\bf a} $,
so w"urde schon die\-ses Produkt selbst zu $ {\bf a} $ ge\-h"or\-en. Hat um\-gekehrt ein
Filter $ F \in \cal Z $ diese Eigenschaft, so mu\3 er maximal sein. Aus
$ F \subset G $ ergibt sich die Existenz eines $ g \in G $,
$g \not\in F $. Dann ist nach Voraussetzung "uber $ F $ f"ur gewisse
$ n \in {\bf N} $ und $ f \in F $ $ g^n \cdot f \in {\bf a} \cap G $, also
$ G \not\in \cal Z $.
\begin{kor}
Sei $ M $ ein kommutatives Monoid mit Nullelement $ 0 \neq 1 $. Dann ist jeder 
konsistente Filter in einem Ultrafilter enthalten. Insbesondere gibt es in
$ M $ Ultrafilter, und ein Filter $ F $ ist genau dann Ulrafilter,
wenn es f"ur alle
$ g \not\in F $ ein $ n \in {\bf N} $ und ein $ f \in F $ gibt mit
$ g^n \cdot f = 0 $.
\end{kor}
Beweis. Man nehme im obigem Lemma einfach $ {\bf a} = \{0\} $ und f"ur den
Existenz\-beweis $ S = A^\times $.
\begin{kor}
Sei $ A $ ein kommutativer Ring und {\bf a} ein Ideal in $ A $. Dann gibt es
(relativ) maximale Filter $F$ in $ A $ mit $ F \cap {\bf a} = \emptyset $. Die
Komplemente dieser Filter sind minimale Primoberideale zu {\bf a}.
\end{kor}
Beweis. Nur der Zusatz ist hier neu, wobei lediglich die
additive Abge\-schlos\-sen\-heit
zu zeigen ist. Seien $ x,y \not\in F $. Dann gibt es $ n,m \in {\bf N} $ und
$ f,g \in F $ mit $ x^n \cdot f = a \in {\bf a} $ und $ y^m \cdot g = b \in
{\bf a} $. F"ur geeignetes $ k \in {\bf N} $ gilt dann
$ (x+y)^k = x^n \cdot r + y^m \cdot s $. Dies multipliziert mit $ f \cdot g $
ergibt ein Element aus {\bf a} und geh"ort daher nicht zu $ F $, also kann auch
$ x + y $ nicht zu $ F $ geh"oren. Die Minimalit"atseigenschaft ist klar.
\begin{kor}
In einem kommutativen Ring $ A $ entsprechen sich Ultrafilter und minimale 
Primideale bijektiv mittels Komplementbildung.
\end{kor}
Beweis. Sei $ {\bf p} $ ein minimales Primideal, sei $ F := A -{\bf p} $ und sei
$ F \subset G $.
Dann ist $ {\bf p} \supset A-G = \bigcup_{i \in I} {\bf p}_i $
nach Satz 1.1.1 und wegen der Minimalit"at von ${\bf p}$ mu\3 $ I $ leer sein.
 Also ist $ G = A $.
\par\bigskip
In Integrit"atsbereichen oder allgemeiner in integren Monoiden mit Null ist
nat"urlich $ A^ \ast $ der einzige Ultrafilter. Generell gilt
in einem Monoid $M$
mit Null stets $ M^ \ast \subseteq F $ f"ur jeden Ultrafilter $ F $, da
ein Nichtnullteiler zu einem konsistenten Filter stets konsistent
hinzugef"ugt werden kann. Insbesondere ist $M^\ast$ im Durchschnitt aller
Ultrafilter enthalten. Besitzt M keine nilpoten\-ten Elemente, so gilt hier
sogar die Gleichheit.
\begin{kor}
Sei $ M $ ein kommutatives, reduziertes Monoid mit $ 0 \neq 1 $.
Dann ist $ M ^\ast $ gleich
dem Durchschnitt aller Ultrafilter.
\end{kor}
Beweis. Ist $ f \in M $ ein Nullteiler, so gibt es ein $ g \in M , g \neq 0 $
mit $ f \cdot g =0 $. $ g $ ist nicht nilpotent, daher ist der Filter $ F(g) $
konsistent und in einem Ultrafilter $ F $ enthalten. Daher ist $ f \not\in F $.
\par\bigskip
In booleschen Ringen ergeben sich einige Besonderheiten, die im folgen\-den Satz
zum Ausdruck kommen.
\begin{satz}
In einem booleschen Ring $ A $ gelten folgende Aussagen:
\par\smallskip\noindent
(1) Ideale und Filter entsprechen sich bijektiv mittels
`innerer Komple\-ment\-bildung'
$ {\bf a} \longmapsto \{1-e : e \in {\bf a} \} $ und
$ F \longmapsto \{1-e : e \in F \} $.
\par \smallskip \noindent
(2) Ein Filter $ F $ ist genau dann ein Ultrafilter, wenn f"ur alle $ e \in A $
gilt: entweder $ e \in F $ oder $ 1-e \in F $.
\par\smallskip\noindent
(3) Jeder Filter ist Durchschnitt von Ultrafiltern.
\end{satz}
Beweis. (1) Die Bijektivit"at der Zuordnung ist klar. Sei {\bf a} ein Ideal
und $ F $ die ihm zugeordnete Menge. F"ur jeden Ring ist $ F $ ein
multiplikatives System. Zum Nachweis der Teilerstabilit"at
sei $ f \cdot g = 1-e $,
$ e \in {\bf a} $ angenommen. Dann ist $ e \cdot(1-f) = (1-f \cdot g)
\cdot (1-f) = 1-f- f\cdot g + f^2 \cdot g = 1-f \in {\bf a} $, also $ f \in F $.
Sei umgekehrt $ F $ ein Filter und $ {\bf a} $ die ihm zugeordnete Menge. Sie
enth"alt die Null. Zum Nachweis der Additivit"at seien $ e,f \in F $, dann ist
$ 1-e+1-f = e+f $ und $ 1-e-f \in F $ ist zu zeigen, was sich sofort aus
$ (1-e-f)\cdot e \cdot f = -e\cdot f $ ergibt. Schlie\3lich sei
$ g = f \cdot (1-e) $ mit $ e \in F $. Dann ist
$ (1-g) \cdot e = (1-f-e \cdot f) \cdot e = e $, also $ g \in{\bf a} $.
\par\noindent
(2) Sei $ F $ ein Ultrafilter und $ e \in A $. Wegen $ e \cdot (1-e) = 0 $
k"onnen nicht beide zu $ F $ geh"oren, und da das Komplement von $ F $ ein
Primideal ist, k"onnen nicht beide nicht zu $ F $ geh"oren. Umgekehrt, sei
$ F $ ein Filter mit der angegebenen Eigenschaft. Dann ist $F$ konsistent.
Ist $ e \not\in F $, so ist
$ 1-e \in F $ und $ e \cdot (1-e) = 0 $, d.h. jedes Element, das nicht zu dem
Filter geh"ort, wird von einem Element aus dem Filter annulliert, $ F $ ist
also ein Ultrafilter.
\par\noindent
(3) Sei $ F $ ein Filter, die Behauptung ist $ F = \bigcap G_i $ mit
$ F \subseteq G_i , G_i $ Ultrafilter. Die eine Inklusion ist klar, und der
inkonsistente Filter ist gleich dem leeren Durchschnitt, sei also $ F $ 
konsistent und $ e \not\in F $. Dann ist auch
$ G := F (F \cup \{1-e \}) $
konsistent, denn $ f \cdot (1-e) = 0 $ f"ur ein $ f \in F $ bedeutet
$ f = f \cdot e $ , also $ e \in F $, was nach Voraussetzung nicht der Fall 
ist. $ G $ ist daher in einem
Ultrafilter $ G_i $ enthalten, der $ e $ nicht enth"alt.
\par\bigskip\noindent
{\bf Bemerkung} Die Aussage (3) im vorstehenden Satz gilt allgemeiner in
null\-dimensionalen Ringen, wozu die booleschen Ringe geh"oren. Die Aussage beruht
dann auf der Tatsachen, da\3 Filter stets Durchschnitte von Primideal\-komplementen
sind, die im Nulldimensionalen alle Ultrafilter sind.
\par\smallskip\noindent
{\bf Bemerkung} Die Zuordnungen im vorangegangenen Satz sind Spezialf"alle
folgender Zuordnungen zwischen Idealen und Filtern in beliebigen kommuta\-ti\-ven
Ringen. Einem Ideal kann man mittels der Restklassenbildung
$\varphi : A \longrightarrow A_{/{\bf a}} $ den Filter
$ \varphi^{-1}((A_{/{\bf a}})^\times) = F(1+{\bf a}) $ zuordnen, und einem 
Filter $F \subseteq A $ den Kern der Abbildung $ A \longrightarrow A_F $.
Im allgemeinen handelt es sich hierbei nicht um eine Korrespondenz wie oben.
\par\bigskip
In Integrit"atsbereichen spielen natur\-gem"a\3 Ultra\-filter kei\-ne gro\3e Rol\-le.
Geht man eine Stufe tiefer, so erh"alt man eine Beziehung zwischen den Filtern,
die genau eine echte, konsistente Erweiterung besitzen (n"amlich den Filter
$ A^\ast $), und den Primidealen der H"ohe 1, die analog zu der zwischen
Ultrafiltern und minimalen Primidealen in beliebigen Ringen ist.
\begin{satz}
Sei $ A $ ein Integrit"atsbereich. Dann entsprechen sich mittels
Kom\-plementbildung
die Primideale der H"ohe 1 und die Filter, die nur $ A^\ast $ als echte,
konsistente Erweiterung besitzen.
\end{satz}
Beweis. Sei {\bf p} ein Primideal der H"ohe 1, und sei $ A- {\bf p} \subset G $,
wobei $ G $ ein konsistenter Filter sei. Dann ist
$ {\bf p} \supset A-G = \bigcup_{i \in I} {\bf p}_i $, $ I \neq \emptyset $, also
$ {\bf p}_i = \langle 0 \rangle $, also $ G = A^\ast $.
Ist umgekehrt $ F $ ein Filter mit der
angegebenen Eigenschaft, so ist $ A-F = \bigcup_{i \in I} {\bf p}_i $ , und ein
$ {\bf p}_i =: {\bf q} $ ist vom Nullideal verschieden. Damit ist
$ F \subseteq A-{\bf q} \subset A^\ast $, also gilt vorne die Gleichheit und
$ {\bf q} $ kann au\3er dem Nullideal kein weiteres Primideal echt umfassen.
\par\bigskip\noindent
\section{Filter und Teilbarkeitseigenschaften}
In diesem Abschnitt wird untersucht, durch welche Elemente ein Filter
cha\-rak\-terisiert werden kann, und wie ein von Elementen eines bestimmten Typs erzeugter
Filter aussieht. Die Beantwortung dieser Frage h"angt dabei wesent\-lich von den
Teilbarkeitseigenschaften des betrachteten Monoids bzw. Ringes ab. Umgekehrt
erlauben Filter die Neubildung von Begriffen im Umkreis der Teilbarkeitslehre,
deren Verh"altnis zu den "ublichen gekl"art werden soll. Im Vordergrund werden
dabei die Integrit"atsbereiche stehen, in denen jede von Null verschiedene
Nichteinheit Produkt unzerlegbarer Elemente ist (hierf"ur ist die Bezeichnung
$ {\bf F_0} $ im Gebrauch). Wir beginnen mit einer Liste
einfacher Sachverhalte.
\begin{pro}
Sei $ M $ ein kommutatives Monoid, $ g,f \in M $. Dann gilt:
\par\smallskip\noindent
(1) Es sind "aquivalent: (a) $ g $ teilt eine Potenz von $ f $,
(b) $ g \in F(f) $,
(c) $ F(g) \subseteq F(f) $ , (d) $ {\rm Rad}(f) \subseteq {\rm Rad}(g) $,
(e) $ D(f) \subseteq D(g) $ in {\rm Spek}$ M $, wobei in (d) und (e) $ M $ als Ring
vorausgesetzt wird.
\par\smallskip\noindent
(2) $ F(1) = M^\times $, $ F(0) = M $.
\par\smallskip\noindent
(3) $ F(F(f) \cup F(g)) = F(f,g) = F(f \cdot g) $.
\end{pro}
\begin{lem}
Sei $ A $ ein $ {\bf F_0} $-Bereich, $ P $ ein Repr"asentantensystem assoziier\-ter
irreduzibler Elemente und $ F $ ein konsistenter Filter.
Dann ist $ F = F(F \cap P) $,
das hei\3t F l"a\3t sich aus seinen irreduziblen Elementen rekonstruieren.
\end{lem}
Beweis. Die eine Inklusion ist klar, zum Beweis der anderen sei $ f \in F $
vorgegeben. Wegen $ f \neq 0 $ und der Voraussetzung "uber den Ring gibt es
irreduzible Elemente $ p_1 , ... , p_n $ aus $ P $ und eine
Einheit $ e $ mit $ f = e \cdot p_1 \cdot ... \cdot p_n $. Dann geh"oren die
$ p_i $ zu $ F $ und damit ist $ f \in F(F \cap P) $.
\begin{lem}
Sei $ A $ ein Integrit"atsbereich und $ T $ eine Menge von Primele\-menten aus
$ A $. Dann besteht $ F(T) $ aus s"amtlichen Produkten von Elementen aus $ T $
und Einheiten, also $ F(T) = \{ e \cdot p_1 \cdot ... \cdot p_n :
e \in A^\times, n \in {\bf N}, p_i \in T \} $.
\end{lem}
Beweis. Es ist klar, da\3 die angegebenen Produkte ein multiplikatives System
bilden und zum Filter $ F(T) $ geh"oren. Sei $ f \in F(T) $. Dann teilt $ f $
eines der Produkte, also $ f \cdot g = e \cdot p_1 \cdot ... \cdot p_n $.
Jedes $ p_i $ teilt also das Produkt $ f \cdot g $, und damit, da es prim ist,
einen der Faktoren, woraus sich durch k"urzen die Behauptung ergibt. Formal
beweist man die Aussage durch Induktion "uber $ n $ : Ist $ n = 0 $, so ist
$ f $ eine Einheit und daher besitzt $ f $ die gew"unschte Darstellung, sei
also $ n \ge 1 $ und $ q = p_n $. Teilt $ q $ den Faktor $ g $, so ergibt sich
nach k"urzen $ f \cdot \tilde{g} = e \cdot p_1 \cdot ... p_{n-1} $ und die
Behauptung aus der Induktions\-vor\-aussetzung. Teilt hingegen $ q $ nicht den
Faktor $ g $, so wird $ f $ von $ q $ geteilt, also $ f = \tilde{f} \cdot q $.
Nach k"urzen ergibt sich $ \tilde{f} \cdot g = e \cdot p_1 \cdot ... p_{n-1} $, also
nach Induktionsvoraussetzung eine Darstellung der gew"unschten Form f"ur
$ \tilde{f} $ und nach Multiplikation mit $ q $ eine solche f"ur $ f $.
\par\bigskip\noindent
Aus diesen zwei Lemmata gewinnt man sofort eine "Ubersicht "uber die Filter
in faktoriellen Ringen. Diese haben die Eigenschaft $ {\bf F_0} $ und in ihnen
ist jedes irreduzible Element prim.
\begin{kor}
In einem faktoriellen Ring $ A $ entsprechen die konsistenten
Filter ein\-deutig den Teil\-mengen eines Repr"a\-sen\-tanten\-systems $ P $ von 
"Aquiva\-lenz\-klassen assozi\-ierter Prim\-elemente, indem man 
einem konsistenten Filter $ F $ die Menge $ F \cap P $ 
zuordnet, und umge\-kehrt einer Teilmenge $ T \subseteq P $ den von $ T $  
erzeugten Filter. Das Komplement eines Filters $ F \neq A^\ast $ ist dar\-stellbar
als Verein\-igung ein\-deutig bestimmter Prim\-haupt\-ideale.
\end{kor}
\par\bigskip\noindent
{\bf Bemerkung} Ist $M$ das freie Monoid zu den `Basiselementen'
$X_i , i \in I $, so gilt die gleiche Charakterisierung der Filter wie in
einem faktoriellen Ring, d.h. die Filter entsprechen den Teilmengen von $I$.
\par\smallskip\noindent
{\bf Bemerkung}
Ist $ A $ kein faktorieller Bereich, so l"a\3t sich der von einer Menge von
irreduziblen Elementen erzeugte Filter im allgemeinen nicht ein\-fach
be\-schrei\-ben. Schon wenn man mit einem einzelnen irreduziblen Element startet,
kann der erzeugte Filter zus"atzliche (nicht assoziierte) irreduzible 
Elemente enthalten. In $ D = {\bf Z}[i \sqrt {5} ] $ hat man beispielsweise 
wegen $ (1 + i \sqrt {5})^2 = -4+2 \sqrt {5} $ , $ 2 \in F(1 + i \sqrt {5}) $.
Die\-ses Ph"anomen tritt in einem Integrit"atsbereich $ A $ mit der Eigenschaft
$ {\bf F_0} $ sogar immer 
dann auf, wenn $ A $ nicht faktoriell ist. Sind n"amlich $ p_1,...,p_n $ die 
unzerlegbaren Elemente einer Zerlegung von $ f $, und besitzt der von diesen 
unzerlegbaren Elementen erzeugte Filter keine weiteren, so ist die Zerlegung 
bis auf Reihenfolge und Assoziiertheit eindeutig, woraus die Faktorialit"at von
$ A $ folgt (\cite{rei}, Satz 160).
\par\bigskip\noindent
Nach obigem Lemma besteht der von einem Primelement $ p $ 
erzeugte Fil\-ter ein\-fach aus Produkten von Einheiten und Potenzen von $ p $.
Ist nun $ p $ 
Nicht\-null\-teiler und Nichteinheit mit dieser Eigenschaft, so ist $ p $ 
irreduzibel: aus $ x \cdot y = p $ ergibt sich unter dieser Voraussetzung 
$ x = e_1 \cdot p^n $ und $ y = e_2 \cdot p^m $ mit $ e_1 , e_2 $ Einheiten. 
K"urzen mit p liefert $ n $ oder $ m = 0 $, also ist $ x $ oder $ y $ eine 
Einheit. Somit l"a\3t sich mit Filtern eine Unzerleg\-barkeits\-eigenschaft
beschreiben, die sich (echt, wie die beiden folgenden Beispiele zeigen)
zwischen prim und irreduzibel befindet.
\par\bigskip\noindent
{\bf Beispiel 1} Sei $ A $ ein kommutativer Ring mit einem Element $ a \in A $
mit $ a \neq 0 $, $ a^2 = 0 $. In $ A[X] $ ist X irreduzibel, wegen
$ (X+a) \cdot (X-a) = X^2 $ aber ist $ X+a \in F(X) $, und $ X+a $ ist nicht
assoziiert zu einer Potenz von $ X $.
\par\smallskip\noindent
{\bf Beispiel 2} In $ D $ ist 2 nicht prim, aber $ F(2) $ besteht nur aus 
Assoziierten von Zweierpotenzen. Der Beweis ist nicht ganz einfach und
erfordert einige Grundtatsachen "uber Dedekindringe. Zun"achst
ist $ {\bf p} = \langle 2,1+i \sqrt{5} \rangle $ ein maximales Ideal,
kein Hauptideal, und es
ist $ 2D = {\bf p}^2 $, (\cite{rei}, p. 148). Aus $ x \in F(2) $ ,
$ x \neq D^\times $
ergibt sich Rad$(2) \subseteq $ Rad$(x)$, und
wegen der Maximalit"at von
Rad$(2) = {\bf p} $ ergibt sich Rad$(2) =$
Rad$(x) $, und damit
$ xD = {\bf p}^n $. $ n =1 $ kann nicht sein, da {\bf p} kein Hauptideal ist.
Ebenso kann $ n $ und $ 2 $ nicht teilerfremd sein. W"are
n"amlich $ n $ ungerade, so gebe es eine nat"urliche Zahl $ k $
mit $ n-k2 = 1 $,
und dann h"atte der Divisor von $ x/2^k $ an der Stelle $ {\bf p} $ den
Wert $ 1 $ und w"are sonst "uberall null. Dann w"are $ x/2^k \in D $
und w"urde {\bf p} erzeugen. Also ist $ xD = {\bf p}^{2k} = ({\bf p}^2)^k = 
(2D)^k = (2^k)D $, und damit ist $ x $ zu einer Potenz von $ 2 $ assoziiert.
Als Nebenresultat erhalten wir, da\3 die Gleichung $ a^2+5b^2 = 2^n $
f"ur $ b \neq 0 $ in $\bf Z $ keine L"osung besitzt.
\par\bigskip 
In einem noetherschen Ring $ A $ wird jede absteigende Kette von Hauptfil\-tern
station"ar, da die entsprechende Aussage f"ur die aufsteigende Kette der
zuge\-h"o\-ri\-gen Radikalideale gilt. Daraus ergibt sich aber nicht, da\3 es in einem
Haupt\-filter nur endlich viele wesentlich
verschiedene irre\-duzible Elemente gibt.
Unter st"arkeren Voraussetzungen l"a\3t sich dies beweisen.
\begin{satz}
In einem normalen, noetherschen Integrit"atsbereich $ A $ besitzt ein
konsistenter Hauptfilter
$ F(f) $ nur endlich viele nicht assoziierte irreduzible Elemente.
\end{satz}
Beweis. Sei $ K = Q(A) $ der Quotientenk"orper von $ A $. Aufgrund der
Voraus\-setz\-ung
"uber den Ring hat man einen Gruppenmorphismus von $ K^\ast $
in die Di\-vi\-sor\-engruppe von $ A $. Diese Abbildung besitzt folgende Eigenschaften,
die wir zum Beweis ben"otigen:
\par\noindent
(1) F"ur $ q \in K^\ast $ gilt: $ (q) \ge 0 $ genau dann, wenn $ q \in A $.
\par\noindent
(2) F"ur $ a,b \in A $ gilt: $ (a) = (b) $ genau dann, wenn $ a $ und $ b $ in
$ A $ assoziiert sind.
\par\noindent
(3) F"ur $ a \in A $ gilt: $ a $ ist irreduzibel genau dann,
wenn $ a \neq 0 $,
$ a \not\in A^\times $ und $ (a) $ minimal unter allen effektiven
Hauptdivisoren ist.
\par\noindent
Zu (1) verweisen wir auf die Literatur (\cite{mat}, Satz 38),
(2) und (3) sind direkte Folgerungen
davon. Sei nun $ f \in A , f \neq 0 $ vorgegeben. Der Divisor von $ f $ hat nur an
endlich vielen Stellen (= Primidealen der H"ohe 1) einen Wert ungleich null, und
die Divisoren von $ f^n , n \in {\bf N}, $ haben an denselben Stellen einen Wert
ungleich null. Der Divisor eines Teilers $ g $ einer Potenz von $ f $ ist daher
allenfalls an diesen endlich vielen Stellen von null verschieden, wobei die
Divisoren von irreduziblen Teilern noch zus"atzlich minimal sind unter allen
$ (g) , g\in F(f), g \not\in A^\times $. Die Behauptung ergibt sich dann aus
folgendem Lemma.
\begin{lem}
Sei $ M \subseteq {\bf N}^n $ eine Teilmenge, die mit der nat"urlichen Ord\-nung
versehen sei. Dann besitzt $ M $ nur endlich viele minimale Elemente.
\end{lem}
Beweis. Induktion "uber $ n $, wobei $ n=0,1 $ klar sind, also $ n \ge 2 $.
Ohne Beschr"ankung der Allgemeing"ultigkeit seien die Elemente aus $ M $
paarweise unvergleichbar. F"ur $ i = 1,...,n $ und $ a \in M $ sei $ a_i $ die
$ i-te $ Komponente von $ a $ und $ a^i $ das $ (n-1)- $ Tupel, das aus $ a $
entsteht,
wenn man die $ i-te $ Komponente wegl"a\3t. $ M^i $ sei die Menge aller $ a^i $
und
$ \tilde{M}^i $ die Menge der minimalen Elemente von $ M^i $. F"ur minimale
Elemente
$ z \in M $ sind nun zwei F"alle zu unterscheiden.
\par\noindent
1. Fall. $ z^i $ ist f"ur ein $ i $ minimal in $ M^i $. Nach Induktionsvoraussetzung
ist f"ur jedes $ i $ die Menge $ \tilde{M}^i $ endlich, und wegen
$ 1 \le i \le n $ kann es nur endlich viele $ z \in M $ von diesem Typ geben.
\par\noindent
2. Fall. $ z^i $ ist f"ur alle $ i $ nicht minimal in $ M^i $. Zu jedem $ i $
gibt es also ein $ a \in M $ mit $ z^i > a^i $, wobei man $ a^i $ als minimal
annehmen darf. Da $ z $ aber minimal in $ M $ ist, gilt $ z_i < a_i $, und
damit $ z_i < max \{a_i : a\in M, a^i \mbox{ist minimal in} M^i \} $.
Diese Menge ist
aber endlich, da $ \tilde{M}^i $ nach Induktions\-voraus\-setzung endlich ist,
und $ a=b $ aus $ a^i = b^i $ folgt, da unter dieser Voraus\-setzung $ a $ und
$ b $ vergleichbar sind. Damit finden wir f"ur die $ z \in M $ von diesem
Typ obere Schranken f"ur jede Komponente, also kann es auch von diesem Typ
nur endlich viele geben.
\begin{kor}
Sei $ A $ ein normaler noetherscher Integrit"atsbereich. Dann ist der Durchschnitt
eines konsistenten Hauptfilters mit einem beliebigen Filter ein
Hauptfilter.
Insbesondere ist der Durchschnitt zweier Hauptfilter wieder ein Hauptfilter.
\end{kor}
Beweis. Seien $ F(f) $ ein konsistenter Hauptfilter und $ G $ ein Filter.
Ein noetherscher
Ring hat die Eigenschaft $ {\bf F_0} $, $ P $ sei ein Repr"asentantensystem
nicht assoziierter irreduzibler Elemente. Nach Lemma 1.2.1 gilt dann
$ F(f) \cap G = F(F(f) \cap G \cap P) $, und nach obigem Satz ist das
erzeugende System aus irreduziblen Elementen endlich, weshalb man es durch ihr
Produkt ersetzen kann. Zum Zusatz hat man nur zu beachten, da\3 die Aussage, wenn
einer der Hauptfilter inkonsistent ist, trivial ist. 
\par\bigskip\noindent
{\bf Beispiel 3} Sei $A$ der Ring der ganz-analytischen Funktionen auf {\bf C},
und $ f=\sin $. Dann geh"oren zu F(f) die unendlich vielen, nicht assoziierten
irreduziblen Polynome $ X-n \cdot \pi , n \in {\bf Z} $. $A$ ist nicht
noethersch.
\par\medskip\noindent
{\bf Beispiel 4} Sei
$ {\bf K}[[X]] , {\bf K} = {\bf R} \mbox{ oder } = {\bf C} $
der Ring der formalen Potenz\-reihen in einer Variablen und $A$ der Unterring von
$ {\bf K}[[X]] $, der aus allen Potenzreihen besteht, deren linearer
Koeffizient verschwindet. Ein Element $ f \in A $ ist genau dann eine
Einheit in $A$, wenn der konstante Term von null verschieden ist, da sich diese
Eigenschaft aus der entsprechenden in $ {\bf K}[[X]] $ "ubertr"agt. $A$ ist ein
lokaler, integrer, noetherscher ($A$ ist ein Restklassenring von
$ {\bf K}[[X,Y]] $), eindimensionaler Ring. Wir zeigen, da\3 es in $A$ ein
Element gibt, das unendlich viele, nicht assoziierte irreduzible Teiler hat,
was dann erst recht auf den davon erzeugten Hauptfilter zutrifft. Dabei lassen
wir uns vom vorangegangenen Beispiel leiten und versuchen, die dort
beschriebenen Tatsachen im Ring $A$ zu repr"asentieren. Die
Potenzreihe des Sinus ist $ X - 1/3! X^3 + 1/5! X^5 - 1/7! X^7 + ... $, und
geh"ort nat"urlich nicht zu $A$, wohl aber die Potenzreihe
$ P := X^4 \cdot \sin X $. Da f"ur jedes $ n \in {\bf Z} $
die Polynomfunktion
$ X-n \cdot \pi $ die Funktion
$ \sin x $ teilt, gibt es auch eine Potenzreihe $ Q \in {\bf K}[[X]] $ mit
$ (X- n \cdot \pi ) \cdot Q = \sin X $, und damit gilt
$ X^2 \cdot (X - n \cdot \pi ) \cdot X^2 \cdot Q = X^4 \cdot \sin X $, wobei
$ X^2 \cdot (X - n \cdot \pi ) = X^3 - n \cdot \pi \cdot X^2 $ und
$ X^2 \cdot Q $ zu $A$ geh"oren. $ X^3 - n \cdot \pi \cdot X^2 , n \in {\bf Z} $
sind in $A$ nicht assoziierte, irredu\-zible Elemente: Aus
$( X^3 - n \cdot \pi \cdot X^2 )\cdot f = X^3 - m \cdot \pi \cdot X^2 $
mit $ f = a_0 + a_2 X^2 + ... $ ergibt sich n"amlich
$ ( X^3 - n \cdot \pi \cdot X^2 ) \cdot f = -n \cdot \pi \cdot a_0 X^2 + a_0 X^3 + ... $,
also ist $a_0 = 1 $ und damit $ n = m $. Zum Beweis der Irreduzibilit"at sei
$ X^3 - n \cdot \pi \cdot X^2 = f \cdot g $ angenommen, und ohne
Einschr"ankung sei $n \neq 0 $. K"urzen mit $X^2$ liefert
in $ {\bf K}[[X]] $ die Beziehung
$ X - n \cdot \pi  = \tilde{f} \cdot \tilde{g} $,
wobei dann $\tilde{f}$ und $\tilde{g}$ Einheiten sind, also einen nicht
verschwindenden konstanten Term haben. $f = X \cdot \tilde{f} $ kann nicht
sein, da ja $f$ zu $A$ geh"ort. Es ist also ohne Einschr"ankung
$ f = \tilde{f} $ , und $f$ ist eine
Einheit in $A$.
\par\bigskip
In einem kommutativen Ring $ A $ hei\3t ein Element $ g $ gr"o\3ter
gemeinsamer Teiler
von Elementen $ f_i , i \in I $, wenn $ g $ jedes $ f_i $ teilt und jeder
gemeinsame Teiler $ h $ der $ f_i $ bereits $ g $ teilt. Da die Eigenschaft,
gemeinsamer Teiler von Elementen $ f_i , i \in I $ zu sein, "aquivalent ist zu
$ \langle f_i , i \in I \rangle \subseteq \langle g \rangle $, ist die
Eigenschaft, da\3 $ g $ ein gr"o\3ter gemeinsamer Teiler von $ f_i $ ist,
gleichbedeutend zu
\begin{eqnarray*}
\langle f_i , i \in I \rangle \subseteq \langle g \rangle , \mbox{ und }
\langle f_i , i \in I \rangle \subseteq \langle h \rangle \mbox{ impliziert }
\\ \langle g \rangle \subseteq \langle h \rangle \hbox{ f"ur beliebiges }
h \in A.
\end{eqnarray*}
Ein ggT erzeugt also das kleinste Hauptideal, das alle $ f_i $ umfa\3t, und ist
in einem Integrit"atsbereich bis auf Assoziiertheit eindeutig bestimmt. In der
Regel gibt es zu einem Ideal kein kleinstes Hauptoberideal. Die Frage, wann es
zu einem Ideal ein kleinstes Hauptoberradikal gibt, f"uhrt uns zur"uck zu
Durchschnitten von Hauptfiltern.
\begin{lem}
Sei $ A $ ein kommutativer Ring, $ g $ und $ f_i , i \in I $ Elemente aus $ A $.
Dann sind "aquivalent:
\par\smallskip\noindent
(1) $ g \in \bigcap_{i \in I} F(f_i) $
\par\smallskip\noindent
(2) $ \langle f_i , i \in I \rangle \subseteq \mbox{\rm Rad}(g) $.
\par\smallskip\noindent
(3) In {\rm Spek}$ A $ ist $ \bigcup_{ i \in I } D(f_i) \subseteq D(g) $.
\end{lem}
Beweis.Alle drei Formulierungen besagen, da\3 $ F(g) \subseteq F(f_i) $ ist f"ur
alle $ i \in I $.
\begin{satz}
Sei $ A $ ein kommutativer Ring, $ g $ und $ f_i $, $i \in I $, Elemente aus $ A $,
so sind "aquivalent:
\par\smallskip\noindent
(1) $ F(g) = \bigcap_{i \in I} F(f_i) $
\par\smallskip\noindent
(2) $ \langle f_i , i \in I \rangle \subseteq \mbox{\rm Rad}(g) $, und f"ur jedes $ h \in A $
mit dieser Eigenschaft ist $ \mbox{\rm Rad}(g) \subseteq \mbox{\rm Rad}(h) $, $ g $ definiert also
ein kleinstes Hauptoberradikal zu den $ f_i , i \in I $.
\end{satz}
Beweis. Aus (1) folgt (2). Da\3 Rad$(g) $ das Ideal der $ f_i $, $i \in I $,
umfa\3t, besagt das Lemma, und aus $ \langle f_i , i \in I \rangle
\subseteq \mbox{Rad}(h) $ folgt aus ihm
$ h \in \bigcap_{i \in I} F(f_i) = F(g) $, also
wiederum Rad$(g) \subseteq \mbox{Rad}(h) $. Aus (2) folgt (1).
Aus dem Lemma folgt
zun"achst $ F(g) \subseteq \bigcap_{i \in I} F(f_i) $, und sei
$ h \in \bigcap_{i \in I} F(f_i) $. Dann ist erneut
$ \langle f_i , i \in I \rangle \subseteq \mbox{Rad}(h) $, woraus aufgrund von (2)
Rad$(g) \subseteq \mbox{Rad}(h) $ folgt. Dies bedeutet aber $ h \in F(g) $.
\par\bigskip\noindent
{\bf Bemerkung} Ist $A$ ein Hauptradikalring - also ein Ring, in dem jedes Radikal
gleich dem Radikal eines Hauptideals ist, oder, anders formuliert, ein Ring,
f"ur den alle offenen Mengen im Spektrum Basismengen sind -, so folgt aus dem Satz
sofort, da\3 der Durchschnitt von Hauptfiltern wieder ein Hauptfilter ist. Dies
gilt auch, wenn $A$ ein eindi\-men\-sio\-na\-ler, noe\-ther\-scher In\-te\-gri\-t"ats\-be\-reich ist.
In diesem Fall besitzt ein vom Null\-ideal ver\-schie\-denes Ideal nur endlich viele
Primoberideale, damit nur endlich viele Oberradikale und insbesondere nur
endlich viele Oberhaupt\-radikale, deren Schnitt das gesuchte kleinste
Oberhauptradikal ist.
\par\medskip\noindent
{\bf Beispiel 5} Sei wieder $ D = {\bf Z}[i\sqrt{5}] $. Dann haben $ 3 $ und
$ 1+i\sqrt{5} $ au\3er Einheiten keine gemeinsamen Teiler. Dagegen gilt
$ F(3) \cap F(1+i\sqrt{5}) = F(2-i\sqrt{5}) $. Wegen $ 9 = (2-i\sqrt{5}) \cdot
(2+i\sqrt{5}) $ und $ (1+i\sqrt{5})^2 = -2 \cdot (2-i\sqrt{5}) $ haben wir
$ 2-i\sqrt{5} \in F(3) \cap F(1+i\sqrt{5}) $, und wegen $ 2-i\sqrt{5} \in
\langle 3, 1+i\sqrt{5} \rangle $ sind die zugeh"origen Radikale gleich, also
auch die Filter. Im allgemeinen kann aber der Durchschnitt zweier Hauptfilter
$ F(f) $ und $ F(g) $ ein Hauptfilter $ F(h) $ sein, ohne da\3 das Radikal von
$ h $ gleich dem von $ f $ und $ g $ ist.
\par\bigskip
Eine besondere Situation liegt vor, wenn der Durchschnitt zweier Hauptfil\-ter
nur aus den Einheiten besteht. So gelangt man zu einem neuen Begriff von
teilerfremd, der in Beziehung zu den einschl"agigen gesetzt werden soll. Zuerst
sei daran erinnert, da\3 zwei von null verschiedene Elemente $f$ und $g$
eines Integrit"atsbereiches
$ A $ teilerfremd hei\3en, wenn $ Af \cap Ag = Afg $. Die Teilerfremdheit l"a\3t
sich auch so ausdr"ucken:
\begin{lem}
Sei $ A $ ein Integrit"atsbereich, $ Q $ sein Quotientenk"orper, und $f,g$
Elemente aus $A^\ast$. Dann sind "aquivalent:
\par\smallskip\noindent
(1) $f$ und $g$ sind teilerfremd.
\par\smallskip\noindent
(2) $f^n$ und $g^m$ sind teilerfremd f"ur alle $m,n \in {\bf N}$.
\par\smallskip\noindent
(3) $A_f \cap A_g = A$, wobei die Ringe als Unterringe von $Q$ aufzufassen sind.
\end{lem}
Beweis. Aus (1) folgt (2) durch Induktion "uber $n+m$, wobei man $n\ge 2$ und
$m\ge 1$ annehmen darf. Sei $x=a\cdot f^n=b\cdot g^m$. Nach Induktionsvoraussetzung
ist $x=c\cdot f^{n-1}\cdot g^m$, und k"urzen mit $f$ liefert $a\cdot f^{n-1}=
c\cdot f^{n-2}\cdot g^m$, woraus aus der Induktionsvoraussetzung wiederum
$a\cdot f^{n-1}=d\cdot f^{n-1}\cdot g^m$ folgt. Multiplikation mit $f$ liefert
dann $ x=d\cdot f^n \cdot g^m$.
\par\noindent
Aus (2) folgt (3). Sei $q \in A_f \cap A_g$, also
$q= a/f^n = b/g^m$ und damit $a\cdot g^m = b\cdot f^n$. Da nach
Voraussetzung $f^n$ und $g^m$ teilerfremd sind, wird $a$ von $f^n$ geteilt, also
$a=c \cdot f^n$ mit $c \in A $ und damit $q= a/f^n=c \in A$.
\par\noindent
Aus (3) folgt
(1). Sei $x=b\cdot f = a\cdot g$. Dann gilt nach Voraussetzung in $Q$:
$ a/f= b/g=c \in A$, und mit $b=c\cdot g $ erh"alt man $x=c\cdot g
\cdot f $.
\begin{satz}
Sei $A$ ein Integrit"atsbereich, und $f,g \in A^\ast$ teilerfremde Ele\-men\-te.
Dann ist $F(f) \cap F(g) = A^\times $.
\end{satz}
Beweis. Sei $h \in F(f) \cap F(g)$, also $h \cdot a= f^n$ und $h \cdot b=g^m$.
Dann ist in $Q$ $ 1/h= a/f^n = b/g^m $, woraus aufgrund
der dritten Formulierung von teilerfremd im obigen Lemma $1/h \in A$
folgt, $h$ ist also eine Einheit.
\par\bigskip
Ist $A$ ein Hauptradikalbereich, so ist $F(f) \cap F(g) = A^\times $ sogar
"aquivalent dazu, da\3 $f$ und $g$ das Einheitsideal erzeugen, woraus sich
die Teilerfremdheit ergibt. In einem faktoriellen Bereich gilt
$F(f) \cap F(g) = F(ggT(f,g))$, und damit auch die Umkehrung des
vorstehenden Satzes. Ebenso unter folgender Bedingung, die schw"acher als die
Faktorialit"at ist, und zum Beispiel auch f"ur Bezoutbereiche erf"ullt ist.
\begin{satz}
Sei $A$ ein Integrit"atsbereich, indem jedes Element des Quotien\-ten\-k"orpers $Q$
eine Darstellung mit teilerfremden Z"ahler und Nenner besitzt.
Sind $f,g \in A^\ast $
mit $F(f) \cap F(g) = A^\times $, so sind $f$ und $g$ teilerfremd.
\end{satz}
Beweis. Sei $q \in A_f \cap A_g \subseteq Q$,
also $q= a/f^n = b/g^m$, und es sei $q = c/h$ eine
teilerfremde Darstellung von $q$. Dann ist $ c \cdot f^n = a \cdot h $, und da
$h$ und $c$ teilerfremd sind, wird $f^n$ von $h$ geteilt, also $ h \in F(f) $.
Ebenso erh"alt man $ h \in F(g) $, und damit ist $h$ nach Voraussetzung eine
Einheit, also ist
$q \in A $, und $f$ und $g$ sind teilerfremd.
\par\bigskip\noindent
{\bf Bemerkung} Die im vorausgegangenen Satz formulierte Bedingung impliziert
bereits die Normalit"at von $A$, da eine Ganzheitsgleichung f"ur $q =f/g$
mit teilerfremden $f$ und $g$ durch Multiplikation mit $g^n$ ergibt, da\3
$f^n$ von $g$ geteilt wird, also $ g \in F(f)$, woraus aufgrund der
Teilerfremdheit und Satz 1.2.11 folgt, da\3 $g$ eine Einheit ist.
\par\medskip\noindent
{\bf Beispiel 6} Sei $ A={\bf Z}[X,6/X] $. In $A$ sind $2$ und $X$
nicht teilerfremd, da $6=X \cdot 6/X $ von $2$ geteilt wird, aber weder
$X$ noch $6/X $ wird von $2$ geteilt, wie man sieht, wenn man einen
entsprechenden Ansatz mit einer geeigneten Potenz von $X$ multipliziert. Ebenso
sieht man, da\3 $F(2)$ beziehungsweise $F(X)$ nur aus den positiven und
negativen Potenzen von $2$ beziehungsweise $X$ besteht; daher ist $F(2) \cap
F(X) = A^\times$.
\par\medskip\noindent

%% file: last.tex
\chapter{Das Filtrum eines Monoids}
\section{Definition und allgemeine Eigenschaften}
Die Menge der Primideale eines kommutativen Ringes $A$ kann man mit einer
Topologie versehen, indem man die Mengen $ D(f) := \{ {\bf p} : {\bf p}
\mbox{ Primideal }, f \not\in {\bf p} \} $ f"ur $f \in A $ als Basismengen
nimmt, was m"oglich
ist wegen $D(f) \cap D(g) = D(f \cdot g) $. Analog hierzu l"a\3t sich die
Menge der Filter eines beliebigen kommutativen Monoids zu einem topologischen
Raum machen.
\par\bigskip\noindent
{\bf Definition und Erl"auterung} Sei $M$ ein kommutatives Monoid,
und sei Filt$M$ die Menge
der Filter von $M$. F"ur $f \in M $ sei $ D(f) :=
\{ F \in \mbox{Filt} M : f \in F \} $.
Wegen der multiplikativen Abgeschlossenheit und der Teilerstabilit"at 
eines Filters gilt
$D(f \cdot g) =D(f) \cap D(g) $. Somit bilden die
Mengen $D(f), f \in M $, die Basis einer Topologie. Filt$M$, versehen mit
dieser Topologie, hei\3t das Filtrum von $M$.
\par\bigskip\noindent
{\bf Bemerkung} Man kann auch einer beliebigen Teilmenge $ N \subseteq M $ eine
offene Teilmenge zuordnen, indem man
\begin{displaymath}
D(N) := \bigcup _ { f \in N} D(f) =
\{ F : F \mbox{ Filter } , F \cap N \neq \emptyset \}
\end{displaymath}
setzt, doch wird hiervon nur selten Gebrauch gemacht werden.
\begin{pro}
Sei $M$ ein kommutatives Monoid, $ X = \mbox{\rm Filt} M $, $ U \subseteq X $ eine
offene Teilmenge, $ f,g \in M $, und $F,G$ seien Filter. Dann gilt:
\par\smallskip\noindent
(1) $ D(f) = X $ genau dann, wenn $f \in M^\times $.
\par\smallskip\noindent
(2) Ist $ F \in U $, und $ F \subseteq G $, dann ist auch $ G \in U $.
\par\smallskip\noindent
(3) $ F(f) \in U $ genau dann, wenn $ D(f) \subseteq U $.
\par\smallskip\noindent
(4) $ D(f) \subseteq D(g) $ genau dann, wenn $ F(g) \subseteq F(f) $.
\par\smallskip\noindent
(5) Ist $ D(f) = \bigcup _{i \in I} U_i $ eine offene "Uberdeckung, so ist
$ D(f) = U_i $ f"ur ein $ i \in I $. Jede offene "Uberdeckung einer Basismenge
besitzt also eine Teil"uberdeckung aus einer Menge.
\par\smallskip\noindent
(6) Eine offene Menge $U$ ist genau dann eine Basismenge, wenn es einen Filter
$F \in U $ gibt, so da\3 $F$ in jedem anderen Filter aus $U$ enthalten ist.
\par\smallskip\noindent
(7) $X$ ist ein $ T_0 $-Raum. 
\par\smallskip\noindent
(8) $ M^\times $ ist der einzige abgeschlossene Punkt in $X$.
\end{pro}
Beweis. (1) Ist $f$ eine Einheit, so ist $f$ in jedem Filter enthalten, und ist
$f$ keine Einheit, so ist $ M^\times \not\in D(f) $. (2) Da $U$ die Vereinigung
von gewissen Basismengen ist, gibt es ein $f$ mit $ F \in D(f) \subseteq U $,
und damit ist auch $f \in G$, also $ G \in U $. (3) Ist $G \in D(f)$, so ist
$ F(f) \subseteq G $, also nach (2) ist $ G \in U $; die Umkehrung ist klar. (4)
Aus $ D(f) \subseteq D(g) $ folgt $ F(f) \in D(g) $, also $ g \in F(f) $; ist
umgekehrt $ g \in F(f) $ und $ G \in D(f) $, so ist auch $ g \in G $. (5) Es ist
$ F(f) \in U_i $ f"ur ein $ i \in I $, und damit nach (3) $ D(f) \subseteq U_i $,
wobei nach Voraussetzung sogar die Gleichheit gilt. (6) Ist $U= D(f)$, so ist
$F(f)$ der kleinste Filter in $U$. Umgekehrt, $F$ sei der kleinste Filter in $U$,
und es sei $ U = \bigcup _{i \in I} D(f_i) $. Dann ist $ f_i \in F $ f"ur ein
$i \in I $, und f"ur einen beliebigen Filter $G \in U $ ist $ f_i \in G $,
also $ U = D(f_i) $. (7) Ist $ F \neq G $, so gibt es ein $ f \in M $ mit
$ f \in F,
f \not\in G $ oder umgekehrt, und $ D(f) $ ist eine trennende Umgebung. (8)
Ist $ F \neq M^\times $, so gibt es eine Nichteinheit $ f \in F $, und $ D(f) $
ist eine offene Umgebung von $F$, in der $ M^\times $ nicht liegt. Ein anderer
Filter kann nicht abgeschlossen sein, da wegen (2) der gesamte Raum die einzige,
offene Umgebung des Einheitsfilters ist.
\begin{kor}
Sei $X$ das Filtrum eines kommutativen Monoids. Dann ist $X$ quasikompakt,
zusammenh"angend, und f"ur jede abelsche Garbe $ \cal G $ ist $ H^i({ \cal G}) =
0 $ f"ur $ i > 0 $.
\end{kor}
Beweis. Die ersten beiden Aussagen ergeben sich sofort aus (5) der obigen
Proposition. Ist $ \cal G $ eine abelsche Garbe auf $X$, so ist der Halm von
$ \cal G $ in einem Punkt gleich dem induktiven Limes "uber allen offenen
Umgebungen des Punktes. Die einzige offene Umgebung des Filters $ M^\times $
ist aber die Gesamtmenge $X$, und daher ist der Halm von $ \cal G $ in
$ M^\times $ gleich der globalen Schnittgruppe $ {\cal G}(X) $. Eine injektive
Aufl"osung ist nun in allen Halmen ex\-akt, insbesondere auch in $ M^\times $, und
daher ist sie auch global exakt, also ver\-schwinden die Kohomologiegruppen.
\par\bigskip\noindent
{\bf Bemerkung} Es ist im allgemeinen nicht einfach, sich eine Vorstellung
von der Topologie auf einem Filtrum zu machen. F"ur ein freies Monoid $M$
mit den Basiselementen $X_i , i \in I $, haben wir die Korrespondenz
Filt$M \cong {\cal P}(I) $, da hier ein Filter eindeutig durch die in ihm
enthaltenen Basiselemente bestimmt ist. Auf der Potenzmenge ${\cal P}(I)$
ist die "ubertragene Topologie gegeben durch die Basismengen
$D(K) = \{ J \in {\cal P}(I) : K  \subseteq J \} $ zu endlichem
$ K \subseteq I $. Da es zu jedem kommutativen Monoid $N$ einen surjektiven
Monoidmorphismus $ M \longrightarrow N $ eines freien Monoids $M$ gibt, ist
wegen Satz 2.2.6 weiter unten Filt$N$ ein\-bett\-bar in eine der\-art to\-po\-lo\-gi\-sier\-te
Po\-tenz\-men\-ge. Dies ist aller\-dings nicht sonderlich erhellend, da "uberhaupt
jeder topologische $T_0-$Raum $X$ diese Eigenschaft hat, wie die
Abbildung $ X \longrightarrow {\cal P}(\mbox{Top}(X) )$ mit
$ x \longmapsto U(x)=\{ U \subseteq X \mbox{ offen } : x \in U \} $ zeigt. 
\par\bigskip
Schon jetzt d"urfte klar sein, da\3 das Filtrum, verglichen mit dem Spek\-trum
eines kommutativen Ringes, ein relativ unkompliziertes Objekt ist. So lassen
sich beispielsweise die Basismengen rein topologisch charakterisieren, wie
der Punkt (7) in obiger Proposition zeigt. Im dritten Abschnitt dieses
Kapitels werden dar"uber\-hinaus
die topologischen R"aume charakterisiert, die als Filtrum auftreten. Dagegen
enth"alt das Filtrum oft interes\-san\-te Teil\-r"aume, wie wir im folgenden noch
sehen werden. Ist $A$ ein kommutativer Ring, so ist Spek$A$ ein Unterraum von
Filt$A$, indem man zu einem Primideal sein Komplement betrachtet. Unterr"aume
von Filt$M$ werden stets mit der induzierten Topologie versehen, und die
induzierten Basismengen werden wieder mit $ D(f) $ bezeichnet, wenn das nicht zu
Mi\3verst"andnissen f"uhrt.
\par\bigskip\noindent
{\bf Definition} Sei $M$ ein kommutatives Monoid mit einem Nullelement. Dann wird mit
Filt$^oM$ der Unterraum der konsistenten Filter bezeichnet.
\begin{pro}
Sei $M$ ein kommutatives Monoid mit Null. Dann gilt
\par\smallskip\noindent
(1) {\rm Filt}$^oM$ ist abgeschlossen in {\rm Filt}$M$.
\par\smallskip\noindent
(2) In {\rm Filt}$^oX$ ist $ D(f) = \emptyset $ genau dann, wenn $f$ nilpotent ist.
\end{pro}
Beweis. (1) Filt$M$ - Filt$^oM$ =\{ $M$ \} = $ D(0) $ ist offen. (2) Ist $f$
nilpotent, so erzeugt $f$ den inkonsistenten Filter, ist also in keinem
konsistenten Filter enthalten. Ist $f$ nicht nilpotent, so ist $F(f)$ konsistent.
\begin{satz}
Sei $M$ ein kommutatives Monoid mit einem Nullelement 0, und sei $Y$ der
Teilraum von {\rm Filt}$M$, der aus allen Ultrafiltern besteht. Dann gilt:
\par\smallskip\noindent
(1) $Y$ ist ein Hausdorffraum.
\par\smallskip\noindent
(2) Die Topologie auf $Y$ hat eine randlose Basis; insbesondere ist $Y$ total
zusammenhangslos.
\par\smallskip\noindent
(3) $Y$ ist dicht in {\rm Filt}$^oM$.
\end{satz}
Beweis. (1) Seien $G \neq F $ Ultrafilter. Dann gibt es ein
$ g \in G, g \not \in F $,
und daher auch aufgrund der Charakterisierung von Ultrafiltern ein
$ n \in {\bf N} $ und ein $ f \in F $ mit $ g^n \cdot f = 0 $. Also ist
$ F \in D(f) $ und $G \in D(g) $, und der Durchschnitt der beiden offenen Mengen
ist leer in $Y$. (2) Es ist zu zeigen, da\3 $ D(f)$ in $Y$ auch abgeschlossen
ist. Hierzu zeigen wir 
\begin{displaymath}
Y-D(f) = \{F \mbox{ Ultrafilter } : f \not\in F \} =
\bigcup_{ g \in M, \hspace{0.5em} 0 \in F(f \cdot g)} D(g)
\end{displaymath}
Ist $F$ ein Ultrafilter mit $ f \not\in F $, so gibt es ein $ g \in M $ und ein
$ n \in {\bf N} $ mit $ f \cdot g^n = 0 $. Es ist also $ F \in D(g) $, und $g$
geh"ort zu der angegebenen Indexmenge. Ist umgekehrt $F$ ein Ultrafilter in
$ D(g) $, und ist $ F(f \cdot g) $ inkonsistent, so kann $f$ nicht zu $F$
geh"oren. Daraus und aus der Hausdorffeigenschaft folgt nun, da\3 in $Y$ nur die
einelementigen Unterr"aume zusammenh"angend sind. (3) Sei $ D(f) $ eine
nichtleere offene Menge in Filt$^oM$. Dann ist $f$ nicht nilpotent, also ist
$ F(f) $ konsistent und in einem Ultrafilter enthalten.
\begin{kor}
Die minimalen Primideale eines kommutativen Ringes $A$ bilden einen dichten,
hausdorffschen, total zusammenhangslosen Unterraum von {\rm Spek}$A$.
\end{kor}
\section{Das Filtrum als Funktor}
Mit dem Spektrum wird einem kommutativen Ring nicht nur ein topologi\-scher Raum
zugeordnet, sondern auch jedem Ringmorphismus eine stetige Abbildung. Dies
l"a\3t sich auch f"ur Monoide und ihre Filtren durchf"uhren. Der wesentliche
Unterschied dabei ist, da\3 dies bei den Filtren nicht nur wie bei den
Spektren in kontravarianter, sondern auch in kovarianter Weise m"oglich ist.
Dadurch lassen sich die Filter zweier Monoide, die durch einen Monoidmorphismus
verbunden sind, wechselseitig in Beziehung setzen.
\begin{satz}
Seien $M$ und $N$ kommutative Monoide und $ \varphi :M \longrightarrow N $ ein
Monoidmorphismus. Dann ist die Abbildung 
\begin{displaymath}
\varphi^\ast : \mbox{\rm Filt} N
\longrightarrow \mbox{\rm Filt} M \mbox{ mit } F \longmapsto \varphi ^{-1} (F) 
\end{displaymath}
stetig. Das Filtrum ist ein kontravarianter Funktor von der Kategorie der
kommutativen Monoide in die Kategorie der topologischen R"aume.
\end{satz}
Beweis. Zun"achst ist zu zeigen, da\3 $ G := \varphi^{-1}(F) $ ein Filter ist. Wegen
$ \varphi(1) = 1 $ ist $ 1 \in G $, wegen $ \varphi(f \cdot g ) = 
\varphi (f) \cdot \varphi (g) $ ist $G$ multiplikativ abgeschlossen und auch
teilerstabil. F"ur eine Basismenge $ D(f) \subseteq \mbox{ Filt}M $ gilt
\begin{eqnarray*}
( \varphi^\ast )^{-1} (D(f) =
\{ F \in \mbox{ Filt}N : \varphi ^{-1}(F) \in D(f) \} \\
= \{ F \in \mbox{ Filt}N : \varphi(f) \in F \} = D( \varphi (f)).
\end{eqnarray*}
Daher ist die Abbildung stetig. Da\3 die Zuordnung mit der
Hinter\-ein\-ander\-schaltung von Monoidmorphismen vertr"aglich ist, gilt
ganz allgemein f"ur das Urbildnehmen.
\par\bigskip
Sind $A$ und $B$ kommutative Ringe und ist $ \varphi : A \longrightarrow B $
ein Ringmorphis\-mus, so f"allt die Einschr"ankung der soeben
definierten Filtrumsabbildung auf die Spektren nat"urlich mit der bekannten
Spektrumsabbildung zusammen. Sind $M$ und $N$ Monoide mit Nullelementen und ist
$ \varphi : M \longrightarrow N $ ein Monoidmor\-phismus, der zu\-s"atz\-lich die
Null auf die Null ab\-bil\-det, so l"a\3t sich die Filtrums\-abbildung zu einer
stetigen Abbildung
$ \varphi^\ast : \mbox{ Filt}^oN \longrightarrow \mbox{ Filt}^oM $
einschr"anken, da $ 0 \in \varphi^{-1}(F) $  sofort $ F = N $ ergibt.
\par\bigskip\noindent
{\bf Bemerkung} Ist $A$ ein kommutativer Ring, so l"a\3t sich mittels der
kanoni\-schen
Einbettung Spek$A \longrightarrow $ Filt$A$ die Strukturgarbe $ \tilde{A} $
auf Filt$A$ "uber\-tra\-gen, so da\3 das Filtrum ebenfalls zu einem beringten
Raum wird. Schon jetzt sei erw"ahnt, da\3 dieser beringte Raum in der Kategorie
der beringten R"aume die gleiche universelle Eigenschaft hat wie das Spektrum
in der Kategorie der lokal beringten R"aume: Ist $ X, {\cal O}_X $ ein beringter
Raum mit $ \Gamma ( {\cal O}_X ) = B $ und ist $ \psi : A \longrightarrow B $
ein Ringmorphismus, so gibt es genau einen zugeh"origen Morphismus beringter
R"aume von $ X, {\cal O}_X $ nach Filt$A, \tilde{A} $. Dies wird in Kapitel 4
genauer ausgef"uhrt.
\begin{satz}
Seien $ M,N $ kommutative Monoide und $ \varphi : M \longrightarrow N $ ein
Mono\-id\-morphismus, so ist die Abbildung
\begin{displaymath}
\varphi_ \ast : \mbox{\rm Filt}M \longrightarrow \mbox{\rm Filt}N \mbox{ mit }
F \longmapsto F( \varphi (F))
\end{displaymath}
stetig, und liefert einen kovarianten Funktor von der Kategorie der 
kommuta\-tiven Monoide in die Kategorie der topologischen R"aume.
\end{satz}
Beweis. Zun"achst sei bemerkt, da\3 das Bild $ \varphi (F) $ eines Filters
$ F \in \mbox{ Filt}M $ ein multiplikatives
System ist; daher besteht $ \varphi_ \ast (F) $ aus allen Teilern von Elementen
$ \varphi (f) $ mit $ f \in F $, was im folgenden st"andig benutzt wird.
Zum Nachweis der Stetigkeit sei $ D(g) $ mit $ g \in N $ vorgegeben. Dann ist
\begin{eqnarray*}
\varphi_ \ast ^{-1} (D(g)) = \{ F \in \mbox{ Filt}M :
g \in \varphi_ \ast (F) \} \\
= \{ F \in \mbox{ Filt}M : \mbox{ es gibt ein } f \in F 
\mbox{ mit } g \mbox{ teilt } \varphi(f) \} \\
= \bigcup_{f \in M : \hspace{0.5em} g \mbox{ teilt } \varphi (f) } D(f)
\end{eqnarray*}
offen. Somit ist noch die Vertr"aglichkeit der Zuordnung mit Verkn"upfungen von
Morphismen zu pr"ufen, seien hierzu
\begin{displaymath}
M_1 \stackrel{ \varphi _1} { \longrightarrow} M_2 \stackrel{ \varphi_2}
{\longrightarrow} M_3 
\end{displaymath}
Monoidmorphismen und $ \psi $ die Hintereinanderschaltung von $ \varphi _1 $
und $ \varphi _2 $, und sei $F$ ein Filter in $ M_1 $. Wir haben zu zeigen,
da\3 $ H_1 := \varphi_{2 \ast} ( \varphi _{1 \ast} (F)) =
F( \varphi_2 (F( \varphi_1(F)))) $ gleich $ H_2 := \psi_\ast(F) =
F( \varphi_2 ( \varphi_1(F))) $ ist, wobei $ H_2 \subseteq H_1 $ klar ist.
Sei also $ z \in H_1 $. Dann gibt es ein $ \acute{z} \in M_3 $ und ein
$ y \in F( \varphi_1 (F)) $ mit $ z \cdot \acute{z} = \varphi_2 (y) $. F"ur $y$
gibt es dann ein $ \acute{y} \in M_2 $ und ein $ x \in F $ mit 
$ y \cdot \acute{y} = \varphi _1 (x) $, und damit ist
\begin{displaymath}
\varphi_2 ( \varphi_1(x)) = \varphi_2 (y \cdot \acute{y}) =
\varphi_2(y) \cdot \varphi_2(\acute{y}) =
z \cdot \acute{z}  \cdot \varphi_2 (\acute{y}) ,
\end{displaymath}
das hei\3t $z$ teilt ein Element aus $ \varphi_2 ( \varphi_1(F)) $ und geh"ort
zu $ H_2 $.
\begin{satz}
Seien $ M_1 $ und $M_2 $ kommutative Monoide und $ M := M_1 \times M_2 $ das
Produktmonoid mit komponentenweiser Verkn"upfung. Dann ist {\rm Filt}$M$ hom"oomorph
zu {\rm Filt}$M_1 \times $ {\rm Filt}$M_2$, versehen mit der Produkttopologie.
\end{satz}
Beweis. Seien $ p_i : M \longrightarrow M_i $ die kanonischen Projektionen,
$ i = 1,2 $, und seien
$ (p_i)_\ast : $ Filt$M \longrightarrow \mbox{\rm Filt}M_i $ die
zugeh"origen stetigen Abbildungen. Dann ist die stetige Abbildung
\begin{displaymath}
(p_1)_\ast \times (p_2)_\ast : \mbox{ Filt}M \longrightarrow
\mbox{ Filt} M_1 \times \mbox{Filt} M_2
\end{displaymath}
als Hom"oomorphie nachzuweisen. Es ist $ (p_i)_\ast (F) = p_i(F) $ f"ur einen
Filter aus $M$, da aus $ g \in (p_1)_\ast (F) $, also $ g \cdot h = f_1 $, mit
$ (f_1,f_2) \in F $, sofort $ (g,1) \cdot (h, f_2) = (f_1,f_2) $, also
$ (g,1) \in F $ und damit $ g \in p_1(F) $ folgt. Sodann ist
$ F = p_1(F) \times p_2(F) $, wobei `$ \subseteq $' klar ist. Zum Beweis der
anderen Inklusion seien
$ f_1 = p_1 (f_1,g_2) , \hspace{0.5em} (f_1,g_2) \in F $ und
$ f_2 = p_2(g_1,f_2) , \hspace{0.5em} (g_1,f_2) \in F $. Dann ist auch
$ (f_1 \cdot g_1, g_2 \cdot f_2) \in F $ und damit $ (f_1,f_2) \in F $.
Andererseits ist $ F = F_1 \times F_2 $ f"ur $ F_i \in \mbox{ Filt} M_i $ ein
Filter in $M$, daher ist die angegebene Abbildung bijektiv, und so bleibt die
Offenheit zu zeigen. Es ist aber
\begin{eqnarray*}
(p_1)_\ast \times (p_2)_\ast (D((f_1,f_2)))  = \{ p_1(F) \times p_2(F) :
(f_1,f_2) \in F \} \\ = \{ F_1 \times F_2 : F_1 \in D(f_1) , F_2 \in D(f_2) \}
= D(f_1) \times D(f_2).
\end{eqnarray*}
\par\bigskip
Sind $ M$ und $N$ kommutative Monoide und ist
$ \varphi : M \longrightarrow N $
ein Monoid\-mor\-phismus, so f"uhrt die Hinter\-ein\-an\-der\-schal\-tung von
$ \varphi _\ast $ und $ \varphi^\ast $ zu stetigen Abbil\-dungen
\begin{displaymath}
\varphi^\ast \varphi_\ast : \mbox{ Filt}M \longrightarrow \mbox{ Filt}M
\mbox { und }
\varphi_\ast \varphi^\ast : \mbox{ Filt}N \longrightarrow \mbox{ Filt}N ,
\end{displaymath}
denen wir uns nun unter Beibehaltung der Bezeichnungen zuwenden wollen.
\par\bigskip\noindent
{\bf Definition} Ein Filter $ F \in \mbox{ Filt}M $ hei\3t
ein Fixfilter von $M$ bez"uglich $ \varphi $, wenn
$ \varphi^\ast \varphi_\ast (F) = F $, und $ G \in \mbox{ Filt}N $ hei\3t ein
Fixfilter von $N$ bez"uglich $ \varphi $, wenn
$ \varphi_\ast \varphi^\ast (G) = G $.
\begin{lem}
Die Abbildung $ \varphi^\ast \varphi_\ast $ ist aufsteigend, d.h. es ist
$ F \subseteq \varphi^\ast \varphi_\ast (F) $ f"ur $ F \in \mbox{ \rm Filt}M $,
und die Abbildung $ \varphi_\ast \varphi^\ast $ ist absteigend, d.h. es ist
$ G \supseteq \varphi_\ast \varphi^\ast (G) $ f"ur $ G \in \mbox{ \rm Filt}N $.
\end{lem}
Beweis. Ist $ f \in F $, so ist
$ \varphi (f) \in \varphi (F) \subseteq \varphi _\ast (F) $, also
$ f \in \varphi^\ast \varphi_\ast (F) $. Zum Beweis der zweiten Behauptung sei
$ g \in \varphi_\ast \varphi^\ast (G) $. D.h. es gibt ein $ h \in N $ mit
$ g \cdot h = \varphi(f) $ und mit $ f \in \varphi^\ast (G) $. Damit ist
$ \varphi (f) \in G $, also auch $ g \in G $.
\begin{satz}
Sei wieder $ \varphi : M \longrightarrow N $ ein Monoidmorphismus. Ist $ F $ ein
Filter in $M$, so ist $ \varphi_\ast (F) $ ein Fixfilter von $N$, und ist $G$ ein
Filter in $N$, so ist $ \varphi^\ast (G) $ ein Fixfilter von $M$.
$ \varphi_\ast $ definiert eine Hom"oomorphie der Fixfilter von $M$ in die
Fixfilter von $N$ mit der Umkehrabbildung $ \varphi^\ast $.
\end{satz}
Beweis. Wegen der Eigenschaft von $ \varphi_\ast \varphi^\ast $, absteigend zu
sein, ist
$ \varphi_\ast (F) \supseteq \varphi_\ast \varphi^\ast \varphi_\ast (F) $. Da
$ \varphi^\ast \varphi_\ast $ aufsteigend ist, ist
$ F \subseteq \varphi^\ast \varphi_\ast (F)$ , und die Monotonie von
$ \varphi_\ast $ liefert die andere Inklusion. V"ollig analog argumentiert man
in der zweiten Situation; der Zusatz ist klar.
\par\bigskip
Wir bestimmen jetzt in einigen spezielleren Situationen die Fixfilter. Die
Inklusionen, die sich dabei direkt aus der auf- oder absteigenden Eigenschaft
der Abbildungen ergeben, werden in den Beweisen gew"ohnlich "ubergangen.
\begin{satz}
Ist $ \varphi : M \longrightarrow N $ ein surjektiver Monoidmorphismus, so ist
jeder Filter aus $N$ ein Fixfilter.
\end{satz}
Beweis. Sei $ g \in G $, $G$ ein Filter in $N$. Wegen der Surjektivit"at
gibt es ein $ f \in M $ mit $ g = \varphi(f) $ . Dann ist
$ f \in \varphi^{-1}(G) = \varphi^\ast (G) $, und damit
$ g = \varphi(f) \in \varphi_\ast \varphi^\ast (G) $.
\begin{satz}
Seien $A,B$ kommutative Ringe und $ \varphi : A \longrightarrow B $ ein 
surjektiver Ringmorphismus mit dem Kern {\bf a}. Dann ist ein Filter $F$ von
$A$ genau dann Fixfilter, wenn f"ur alle $ f \in F $ und $a \in {\bf a} $ gilt
$f+a \in F $.
\end{satz}
Beweis. Da alle Filter aus $B$ fix sind, geht es um die Bestimmung der Filter
$ \varphi^\ast (G) $, $ G \in \mbox{ Filt}B $. Da $f$ und $f+a$,
$ a \in {\bf a} $ dasselbe Bild in $B$ haben, m"ussen die Fixfilter von $A$
die angegebene Eigenschaft haben. Sei nun $F$ ein Filter mit dieser Eigenschaft,
und sei $ h \in \varphi^\ast \varphi_\ast (F) $. Dann ist
$ \bar{h} \in \varphi_\ast (F) $, d.h. $ \bar{h} $ teilt ein $ \bar{f} , f \in F $.
Dann gibt es ein $ x \in A $ und ein $a \in {\bf a} $ mit $ h\cdot x = f + a $.
Nach Voraussetzung "uber $F$ ist dann $ f+a \in F$ und als Teiler ist
$ h \in F $.
\par\bigskip\noindent
{\bf Bemerkung} Die Filter eines kommutativen Ringes $A$, die fix bez"uglich
der Restklassenbildung nach einem Ideal {\bf a} sind, werden auch einfach fix
modulo {\bf a} genannt. Das Urbild der Einheitengruppe von $ A_{/ {\bf a}} $ ist
der kleinste Filter, der fix modulo $ {\bf a} $ ist. Er ist gleich
$ F(1+ {\bf a}) $, wie man sofort nachpr"uft. 
Bez"uglich des Nilradikals ist jeder Filter fix. Ist
n"amlich $ f \in F $ und $a$ nilpotent mit $ a^n = 0 $, so ist wegen
\begin{eqnarray*}
(f - a) \cdot (f^{n-1} + f^{n-2} \cdot a + ... + f \cdot a^{n-2} + a^{n-1}) \\
= f^n - a^n = f^n
\end{eqnarray*}
$ f - a $ ein Teiler einer Potenz von $f$ und geh"ort daher zu $F$. Es ist also
Filt $A$ $ \cong $ Filt$ A_{red} $. Das Komplement eines Primideals
{\bf p} von $A$ ist genau
dann fix modulo {\bf a}, wenn $ {\bf a} \subseteq {\bf p} $ gilt.
Ist n"amlich {\bf p} fix, so folgt aus $ a \in {\bf a} $ und der Annahme
$ a \not\in {\bf p} $ sofort der Widerspruch $ 0 = a-a \not \in {\bf p}$. Ist
umgekehrt $ a \in {\bf a} \subseteq {\bf p} $ und $ f \not\in {\bf p} $, so ist
auch $ a+f \not\in {\bf p} $.
\begin{satz}
Sei $M$ ein kommutatives Monoid, $F$ ein Filter in $M$ und
$ \varphi : M \longrightarrow M_F $ die kanonische Abbildung in die
Brucherweiterung, die genau wie im Fall eines kommutativen Ringes zu
definieren ist. Dann sind alle Filter von $M_F$ Fixfilter, und ihnen entsprechen
in $M$ die Filter, die $F$ umfassen.
\end{satz}
Beweis. Jeder Filter $G$ in $M_F$ umfa\3t die Einheitengruppe, daher ist
$ F \subseteq \varphi^\ast (G) $. Ist nun $ g = a/f \in G $ mit $f \in F$, so ist
$ a/1 \in G $, also $ a \in \varphi^\ast (G) $ und wiederum
$ a/1 \in \varphi_\ast \varphi^\ast (G) $. Wegen $ a/f \cdot f/1 = a/1 $ ist
dann auch $ a/f \in \varphi_\ast \varphi^\ast (G) $. Somit ist nur noch zu
zeigen, da\3 ein Filter $ H \supseteq F $ aus $M$ fix ist. Sei hierzu
$ a \in \varphi^\ast \varphi_\ast (H) $ gegeben, also eine Gleichung in $M_F$
mit $ a/1 \cdot b/f = h/1 $ mit $ b \in M $, $f \in F $ und $ h \in H $. Dies
besagt in $M$, da\3 $ a\cdot b \cdot \acute{f} = h \cdot f \cdot \acute{f} $
f"ur ein $ \acute{f} \in F $. Wegen der Voraussetzung "uber $H$ steht rechts
ein Element aus $H$, das von $a$ geteilt wird, also $a \in H $.
\par\bigskip
Die S"atze 2.2.6 und 2.2.8 haben gemeinsam, da\3 ein Monoidmorphis\-mus
$ \varphi : M \longrightarrow N $ mittels $ \varphi_\ast $ zu einer
Einbettung von Filt$N$ in Filt$M$ f"uhrt. Sind $M$ und $N$ kommutative Ringe,
so erh"alt man aus diesen Abbildungen durch Einschr"ankung auf die
Spektren die kanonischen
Einbettungen. Aus dem letzten Satz ergibt sich speziell
$ D(f) \cong \mbox{ Filt}M_f $, die Basismengen eines Filtrums sind also selbst
wieder Filtren von Monoiden. Wir wenden uns nun Morphismen zu, bei denen
s"amtliche Filter des Definitionsmonoids fix sind.
\begin{satz}
Ist $\varphi : M \longrightarrow N $ ein Monoidmorphismus, dann sind "aquivalent
\par\smallskip\noindent
(1) Alle Filter von $M$ sind fix.
\par\smallskip\noindent
(2) Alle Hauptfilter $F(f)$ von $M$ sind fix.
\par\smallskip\noindent
Ist $M$ zus"atzlich ein kommutativer Ring, so ist dies au\3erdem "aquivalent zu
\par\smallskip\noindent
(3) Alle Primideale von $M$ sind fix.
\end{satz}
Beweis. Von (1) nach (2) ist eine Einschr"ankung, seien also die Hauptfilter
fix und ein beliebiger Filter $F$ vorgegeben. Es ist
$ F = \varphi^\ast \varphi_\ast (F) $ zu zeigen, sei also $ g \in M $ mit
$ \varphi(g) \in \varphi_\ast (F) $. Das hei\3t, es gibt ein $ f \in F $ und
ein $ h \in N $ mit $ \varphi(g) \cdot h = \varphi(f) $, und das bedeutet
$ \varphi(g) \in \varphi_\ast (F(f)) $, woraus nach Voraussetzung
$ g \in F(f) $, und damit $ g \in F $ folgt.
\par\noindent
Sei nun $M$ ein Ring, alle Primideale fix und
$\varphi(g) \in \varphi_\ast(F(f)) $. Sei $g \not\in F(f) $ angenommen. Dann
gibt es ein Primideal mit $ g \in {\bf p} , f \not\in {\bf p} $. Damit ist
$ \varphi_\ast(F(f)) \subseteq \varphi_\ast (M-{\bf p}) $ und damit auch
$ \varphi(g) \in \varphi_\ast (M-{\bf p}) $, woraus aus der Fixheit der
Widerspruch $ g \not\in {\bf p} $ folgt.
\par\bigskip\noindent
{\bf Beispiel 1} Sei $\varphi : A \longrightarrow B $ ein injektiver ganzer
Ringmorphismus. Dann sind alle Filter von $A$ fix. Zun"achst ist n"amlich
der Morphismus lokal, d.h. Nichteinheiten
werden auf Nichteinheiten abgebildet,
oder, in unserer Terminologie, die Einheitengruppe von $A$ ist fix. Dies sieht
man so: sei $a \in A $ eine Einheit in $B$, also $a \cdot b = 1 $ f"ur ein
$b \in B $.
Multipliziert man eine Ganzheitsgleichung
$b^n + a_{n-1}\cdot b^{n-1} + ... +a_1 \cdot b + a_0 = 0$ mit $a^n$, so erh"alt
man $ a \cdot a' = 1 $ mit einem $a' \in A $, und wegen der Injektivit"at gilt
diese Gleichheit bereits in $A$. F"ur einen beliebigen Filter $F \subseteq A $
betrachtet man den Ringmorphismus $A_F \longrightarrow B_{\varphi_\ast(F)} $,
der ebenfalls ganz und injektiv ist, und damit wieder lokal. Das bedeutet aber
$ F=\varphi^{-1}((B_{\varphi_\ast(F)}) ^\times) = \varphi^{-1} \varphi_\ast(F) $
\par\smallskip\noindent
{\bf Beispiel 2} Ist $\varphi : A \longrightarrow B $ eine reine Algebra, so
gilt f"ur ein Ideal ${\bf a} $ aus $A$ die Beziehung
$ {\bf a} = \varphi^{-1} ({\bf a}B)$,
(\cite{sch}, \S  88, Aufgabe 6).
Wendet man dies auf die Hauptideale an, so ergibt sich sofort, da\3 die
Hauptfilter fix sind und damit "uberhaupt alle Filter aus $A$.
\par\bigskip
Sei $M$ ein kommutatives Monoid. Wir betrachten die "Aquivalenzrelation $ \sim $
auf $M$, die durch
$ f \sim g : \Longleftrightarrow F(f) = F(g) $ gegeben ist. Auf
$ \tilde{M} := M _{ / \sim } $ l"a\3t sich in kanonischer Weise eine
Monoidstruktur erkl"aren. Sie ist beispielsweise gegeben durch die Menge aller
Hauptfilter $ \{ F(f) : f \in M \} $ mit der Verkn"upfung
$ F(f) \sqcup F(g) := F(F(f) \cup F(g)) = F(f \cdot g ) $ und dem neutralen
Element $ F(1) = M^\times $, oder durch die Basismengen im Filtrum mit dem
Durchschnitt als Verkn"upfung und dem gesamten Filtrum als neutralem Element.
Die kanonische Abbildung $ \varphi : M \longrightarrow \tilde{M} $ ist dabei
ein surjektiver Monoidmorphismus, daher sind alle Filter von $ \tilde{M} $
Fixfilter. Das gleiche gilt aber auch f"ur die Filter aus $M$, da ja
$ g \in \varphi^\ast \varphi_\ast (F) $ f"ur einen Filter $F$ aus $M$
gerade $ D(g) \supseteq D(f) $ f"ur ein $f \in F $ bedeutet, woraus sich
$ g \in F(f) \subseteq F $ ergibt. Wir halten fest:
\begin{satz}
Ist $M$ ein kommutatives Monoid und $ \tilde{M} := M _{/ \sim} $ mit der oben
beschriebenen "Aquivalenzrelation $ \sim $, so ist {\rm Filt}$M \cong $
{\rm Filt}$ M_{/ \sim} $.
\end{satz}
\section{Charakterisierung der Filtrumsr"aume}
In diesem Abschnitt werden die topologischen R"aume, die als Filtren
kommu\-tativer Monoide auftreten k"onnen, topologisch charakterisiert. Dabei
weist der letzte Satz des vorigen Abschnitts den Weg: das Filtrum eines Monoids
ist gleich dem Filtrum des Monoids der Basismengen. Es geht also darum,
Eigenschaften f"ur topologische R"aume zu finden, mit denen sich gewisse offene
Teilmengen als "Basismengen" auszeichnen lassen, die den Basis\-meng\-en in
einem Filtrum entsprechen. Hierzu sind einige Definitionen n"otig.
\par\bigskip\noindent
{\bf Definition} Sei $X$ ein topologischer Raum. Zu $ x \in X $ bezeichne
$ U(x) $ die Menge aller offenen Umgebungen von $ x $.
Auf der Punktmenge $X$ wird in
folgender Weise eine Ordnungsrelation $ \preceq $ eingef"uhrt:
\begin{displaymath}
x \preceq y : \Longleftrightarrow U(x) \subseteq U(y)  .
\end{displaymath}
Es gilt also $ x \preceq y $ genau dann, wenn jede offene Umgebung von $x$ auch
eine offene Umgebung von $y$ ist. In dieser Situation hei\3t $x$
umgebungs\-"armer als $y$ und entsprechend $y$ umgebungsreicher als $x$. Ist
$ T \subseteq X $ eine beliebige Teilmenge, so hei\3t ein Element $ x \in T $
umgebungs"armstes Element von $T$, wenn $x$ umgebungs"armer ist alle anderen
Elemente aus $T$. Besitzt $X$ selbst einen umgebungs"armsten Punkt $x$, so 
hei\3t $X$ ein lokaler Raum mit globalem Punkt $x$.
\par\bigskip\noindent
{\bf Bemerkung} Ein topologischer Raum $X$ ist genau dann ein $ T_0- $Raum, wenn
zwei verschiedene Punkte nicht umgebungsgleich sind. In einem $ T_1- $Raum ist
ein Punkt bez"uglich der Umgebungsrelation nur mit sich selbst vergleichbar.
Ein Punkt $ x \in X $ ist genau dann ein umgebungs"armster Punkt von $ X$, wenn
der gesamte Raum $X$ die einzige offene Umgebung von $x$ ist. In diesem Fall
ist f"ur eine Garbe der Halm in $x$ gleich der globalen Schnittmenge, was das
Bezeichnungspaar lokal/global erkl"art.
\par\bigskip
Im folgenden Satz werden die wesentlichen topologischen Eigen\-schaf\-ten der
Basismengen eines Filtrums zusammengefa\3t, wobei die ersten beiden Aus\-sa\-gen
lediglich Aussagen aus Proposition 2.1.1 in die neue Ter\-mi\-no\-lo\-gie bringen.
\begin{satz}
Sei $M$ ein kommutatives Monoid, $ X =$ {\rm Filt}$M$,
$F,G$ seien Filter in $M$ und $f, f_i \in M $.
Dann gilt
\par\smallskip\noindent
(1) $ F \preceq G $ genau dann, wenn $ F \subseteq G $.
\par\smallskip\noindent
(2) Eine offene Menge $ U \subseteq X $ ist genau dann eine Basismenge, also
vom Typ $ D(f) $ mit einem $ f \in M $, wenn $ U$ ein lokaler Raum ist.
\par\smallskip\noindent
(3) Sei $ T := \bigcap_ {i \in I}D(f_i) $. Dann besitzt $T$ ein
umgebungs"armstes Element.
\par\smallskip\noindent
(4) Ist $ \bigcap _{i \in I} D(f_i) \subset D(f) $, so gibt es eine endliche
Teilmenge $ J \subseteq I $ mit $ \bigcap_{i \in J} D(f_i) \subseteq D(f) $.
\end{satz}
Beweis. (3) Sei $ F_0 := F(f_i : i \in I )$. Dann geh"ort $ F_0 $ zu $T$ 
und f"ur jeden anderen Filter $F \in T $ ist $ F_0 \subseteq F $.
(4) Aus
der Voraussetzung folgt sofort $ F(f_i : i \in I) \in D(f) $, also ist
$ f \in F(f_i : i \in I ) $. $f$ ist also ein Teiler eines Produktes mit
Faktoren $ f_i $, wobei nat"urlich nur endlich viele $ f_i $ als Faktoren
vorkommen. Daher ist $ f \in F(f_i : i \in J) $, mit endlichem
$ J \subseteq I $, und der Weg r"uckw"arts liefert die Behauptung.
\par\bigskip\noindent
\begin{satz}
Ein topologischer Raum $X$ ist genau dann Filtrum eines kommu\-tativen Monoids,
wenn die Menge $ \cal D $ der lokalen offenen Teilmengen von $X$ folgende
Bedingungen erf"ullt:
\par\smallskip\noindent
(1) $X$ ist ein $ T_o $-Raum.
\par\smallskip\noindent
(2) $ X \in \cal D $.
\par\smallskip\noindent
(3) Sind $ D_1 $ und $ D_2 \in \cal D $, so ist auch $ D_1 \cap D_2 \in \cal D $.
\par\smallskip\noindent
(4) $ \cal D $ bildet eine Basis der Topologie von $X$.
\par\smallskip\noindent
(5) Ist $T = \bigcap_{i \in I} D_i $ mit $ D_i \in \cal D $, so besitzt $T$ ein
umgebungs"armstes Element.
\par\smallskip\noindent
(6) Ist $ \bigcap_{i \in I} D_i \subseteq D $ mit $ D, D_i \in \cal D $, so
gibt es eine endliche Teilmenge $ J \subseteq I $ mit
$ \bigcap_{i \in J } D_i \subseteq D $.
\end{satz}
Beweis. Ist $X = \mbox{ Filt}M $ f"ur ein kommutatives Monoid $M$, so besteht
$ \cal D $ gerade aus allen Basismengen $ D(f) , f \in M $, und f"ur diese
Mengen wurden die Aussagen bereits gezeigt. Sei also $X$ ein topologischer
Raum, der die angegebenen Eigenschaften hat. Wir werden zeigen, da\3 $X$
hom"oomorph zu Filt$ \cal D$ ist. Die Eigenschaften (2) und (3) besagen
gerade, da\3 $ \cal D $ mit dem Durchschnitt als Verkn"upfung ein Monoid ist,
so da\3 man von diesem Filtrum sprechen kann. Wir betrachten die Abbildung
\begin{displaymath}
\Psi : X \longrightarrow \mbox{ Filt}{\cal D} \mbox{ mit }
\Psi (x) := \{ D \in {\cal D} : x \in D \}.
\end{displaymath}
$ \Psi (x) $ ist offenbar ein Filter, da die Teilerbeziehung in $ \cal D $
durch die Obermeng\-enbeziehung gegeben ist. Wegen
\begin{displaymath}
\Psi ^{-1} D(D) = \{ x \in X : \Psi (x) \in D(D) \} =
\{ x \in X : D \in \Psi(x) \} = D 
\end{displaymath}
ist $ \Psi $ stetig, und da $ \cal D $ eine Basis von $X$ bildet, folgt hieraus
auch, da\3 $X$ die Initialtopologie tr"agt; es bleibt also noch die
Bijektivit"at zu zeigen.
\par\noindent
Zur Injektivit"at: Seien $ x,y $ verschiedene Punkte in $X$. Da $X$ ein
$T_0-$Raum ist, gibt es eine offene Menge, die den einen enth"alt, den anderen nicht.
Sei also $ x \in U $, $ y \not\in U $. Da $ \cal D $ eine Basis bildet, gibt es
dann auch eine offene lokale Menge $ D \in \cal D $ mit $ x \in D $,
$y \not\in D $, daher ist $ \Psi (x) \neq \Psi (y) $.
\par\noindent
Zur Surjektivit"at: Sei $F$ ein Filter in $ \cal D $. Nach Bedingung (5) gibt
es in $ T := \bigcap_{D \in F} D $ ein umgebungs"armstes Element, das mit
$x$ bezeichnet sei. Es ist $ \Psi(x) = F $ zu zeigen, wobei `` $ \supseteq $ ''
klar ist. Sei also $ x \in D_0 $, $ D_0 \in \cal D $.  Zun"achst ist
\begin{displaymath}
\bigcap _{D \in F } D \subseteq D_0 ,
\end{displaymath}
da aus $ y \in \bigcap_{D \in F } D $ wegen der Wahl von $x$  $ x \preceq y $
und somit $ y \in D_0 $ folgt. Nach Bedingung (6) gibt es dann endlich viele
Elemente $ D_1, ... , D_n \in F $ mit
\begin{displaymath}
\bigcap_{i=1}^n D_i \subseteq D_0 .
\end{displaymath}
$ D_0 $ ist also ein Teiler eines Produktes von Elementen aus $F$ und geh"ort
daher selbst zu $F$.

%% file: tafi.tex
\chapter{Das Filtrum eines topologischen Raumes}
\section{Topologische Filter}
In metrischen R"aumen lassen sich viele topologische Eigen\-schaften mittels
Folgen charakterisieren. So ist beispielsweise eine Abbildung
$ \phi: X \longrightarrow Y $ zwischen metrischen R"aumen $X$ und $Y$ genau dann
stetig in $x \in X $, wenn f"ur jede gegen $x$ konvergente Folge auch die
Bildfolge gegen $ \phi(x) $ konvergiert, oder eine Teilmenge $M \subseteq X $
ist genau dann abgeschlossen, wenn jede Folge in $M$, die konvergiert, ihren
Limes in $M$ hat.
\par
F"ur allgemeine topologische R"aume lassen sich solche
topologischen Ei\-gen\-schaften gew"ohnlich nicht mehr mit Folgen ausdr"ucken, so da\3
man sich gezwungen sieht, eine geeignete Verallgemeinerung von Folgen
und einen umfassenderen Konvergenzbegriff zu suchen. Eine M"oglichkeit liefern
die Net\-ze, die eine direkte Verallgemeinerung der Folgen sind, indem nicht nur
die nat"urlichen Zahlen, sondern jede gerichtete Indexmenge als
Definitions\-be\-reich zugelassen wird.
\par
Eine andere M"oglichkeit bieten die
mengentheoretischen Filter, die 1937 von H. Cartan eingef"uhrt wurden.
Dabei wird ein Filter in einem topolo\-gi\-schen Raum $X$ definiert als
eine nichtleere,
durchschnittstabile Teilmenge der Potenzmenge von $X$, die mit einer Teilmenge
auch jede gr"o\3ere enth"alt. Schon in Beispiel 1 des ersten Kapitels wurde
gezeigt, da\3 diese Definition mit der hier vertretenen zusammenf"allt, wenn
man die Potenzmenge mit dem Durchschnitt als Verkn"upfung als kommutatives
Monoid auffa\3t. Man sagt, da\3 ein Filter $F$ gegen einen Punkt $x$
konvergiert, wenn jede Umgebung von $x$ zu $F$ geh"ort, und kann zeigen, da\3
eine Abbildung zwischen topologischen R"aumen genau dann in einem Punkt $x$ 
stetig ist, wenn f"ur jeden Filter, der gegen $x$ konvergiert, der Bildfilter
gegen den Bildpunkt konvergiert. So k"onnen viele Eigenschaften von
topologischen R"aumen und Abbildungen mit Filtern beschrieben werden. F"ur
genauere Informationen zu (mengentheoreti\-schen) Filtern in der Topologie
sei auf die Literatur verwiesen, (\cite{que}, Kap.5 und Kap.16.C).
\par
In diesem Kapitel geht es nicht um diese mengentheoretischen Filter, sondern
um topologische Filter, die sogleich definiert werden. Sie haben gegen"uber
den mengentheoretischen Filtern sowohl Nach- als auch Vorteile. Ein Nachteil
liegt beispielsweise darin, da\3 eine sinnvolle Verwendung topolo\-gi\-scher
Filter im Zusammenhang mit Abbildungen zwischen topologischen R"aumen die
Ste\-tig\-keit vor\-aus\-setzt und daher mit ihnen ein Begriff wie Ste\-tig\-keit nicht
definiert werden kann. Der grundlegende Vorteil der topolo\-gi\-schen Filter liegt
darin, da\3 sie sich st"arker an der Topologie des Raumes orien\-tieren, w"ahrend
die mengen\-theo\-re\-tischen nur von der zugrunde\-liegen\-den Punkt\-menge des Raumes
abh"angen. Damit steht in Zusammenhang, da\3 die Menge der topologischen Filter
ein Filtrum bilden, das mit dem Ausgangs\-raum eng verkn"upft ist, und mit dem
viele topologische Konstruktionen in einer einheit\-lichen Sprache beschrieben
werden k"onnen. Schlie\3lich lassen sich auf topologische Filter die
allgemeinen "Uberlegungen der vorausgegang\-en\-en Kapitel effektiver
anwenden.
\par\bigskip\noindent
{\bf Definition}
Sei $X$ ein topologischer Raum und Top($X$) das System der offenen Mengen von
$X$. Top($X$) ist mit dem Durchschnitt als
Verkn"upfung und dem gesamten Raum als neutralem Element ein kommutatives
Monoid; die Filter dieses Monoids hei\3en topologische Filter, und wir schreiben
Filt$X$ statt Filt Top($X$).
\par\bigskip
Es sei noch
einmal ausdr"ucklich erw"ahnt, da\3 immer, wenn es sich, wie im vorliegenden
Fall, um ein Untermonoid
einer Potenzmenge handelt, die Beziehung der Teilbarkeit zur
Obermengenbeziehung "aquivalent ist. Eine stetige Abbildung
$ \phi : X \longrightarrow Y $ definiert einen Monoidmorphismus
\begin{displaymath}
\varphi: {\rm Top}(Y) \longrightarrow {\rm Top}(X)  \mbox{ mit }
\varphi (V) := \phi ^{-1} (V) .
\end{displaymath}
Eine stetige Abbildung ist ja gerade so definiert, da\3 das Urbild einer
offenen Menge $ V \subseteq Y $ eine offene Menge in $X$ ist. Insgesamt haben
wir einen Funktor von der Kategorie der topologischen R"aume in die Kategorie
der kommutativen Monoide, und die Filtren liefern wiederum Funktoren in die
Kategorie der topologischen R"aume. Die Untersuchung dieser Funktoren ist,
kurz gesagt, der Gegenstand dieses Kapitels. Zun"achst geben wir einige 
Beispiele f"ur topologische Filter; teilweise handelt es sich um
Konkretisier\-ung\-en 
bisheriger Begriffsbildungen an topologischen R"aumen, teilweise um Analoga zu
bestimmten mengentheoretischen Filtern. Sei hierzu
$X$ ein to\-po\-lo\-gi\-scher Raum und $ M :=$ Top$(X)$.
\par\bigskip\noindent
{\bf Beispiel 1} Die Einheitengruppe von $M$ besteht einfach aus dem Gesamtraum
$X$.
\par\smallskip\noindent
{\bf Beispiel 2} Das Monoid $M$ besitzt die leere Menge als Nullelement, und
der zu einer stetigen Abbildung geh"orige  Monoidmorphismus bildet die
Null auf die Null ab. Die Menge der Nichtnullteiler $ M^\ast $ ist die Menge
aller dichten, offenen Teilmengen von $X$. Die Dichtheit einer Teilmenge
$ U \subseteq X $ besagt ja, da\3 ihr Durchschnitt mit jeder nichtleeren
offenen Menge nicht leer ist. Das Monoid $M$ ist genau dann integer, wenn der
Raum $X$ irreduzibel ist.
\par\smallskip\noindent
{\bf Beispiel 3} F"ur einen Punkt $ x \in X $ ist
$ U(x) := \{ U \in M : x \in U \} $ offenbar ein Filter. Er hei\3t der
Umgebungsfilter von $x$ und spielt eine wichtige Rolle. Generell hei\3t zu
einer beliebigen Teilmenge $ T \subseteq X $ das System der offenen Umgebungen
von $T$ der Umgebungsfilter von $T$. Es gilt $ U(T) = \bigcap_{x \in T} U(x) $.
\par\smallskip\noindent
{\bf Beispiel 4} Ein topologischer Filter $F$ ist genau dann ein Ultrafilter,
wenn es f"ur jede offene Menge $ V \not\in F $ ein Menge $ U \in F $ gibt mit
$ V \cap U = \emptyset $. Der Durchschnitt aller Ultrafilter eines nichtleeren
Raumes ist die Menge der
Nichtnullteiler.
\par\smallskip\noindent
{\bf Beispiel 5} Ist $ x_n , n \in {\bf N} $ eine Folge in $X$, so bilden die
offenen Mengen, die fast alle Folgenglieder enthalten, einen Filter. Die
Folge $x_n$ konvergiert genau dann gegen $x \in X $, wenn
$U(x) \subseteq U(x_n) $ gilt.
\par\bigskip
Da es auf $ M$ = Top$(X) $ auch noch die Vereinigung als Verkn"upfung gibt, lassen
sich einige besondere Filter definieren.
\par\bigskip\noindent
{\bf Definition} Sei $F$ ein topologischer Filter des topologischen Raumes $X$.
\par\smallskip\noindent
(1) $F$ hei\3t ein quasikompakter Filter, wenn aus $ U \cup V \in F $ f"ur
beliebige offene Mengen $ U $ und $V$ stets $U \in F $ oder $V \in F $ gilt.
\par\smallskip\noindent
(2) $F$ hei\3t ein irreduzibler Filter, wenn aus
$ \bigcup _ {i \in I} U_i \in F $ , $ U_i \in M $, die Existenz eines
$ i \in I $ folgt mit $ U_i \in F $.
\par\bigskip\noindent
{\bf Bemerkung} Ein irreduzibler Filter ist konsistent. Dies folgt formal
aus der Definition, wenn man die Konvention beachtet, da\3 die leere Vereinigung
gleich der leeren Menge ist. Jeder irreduzible Filter $F$ konvergiert; gebe es
n"amlich f"ur jedes $ x \in X $ ein offene Menge $ U_x \in U(x) $ mit
$ U_x \not\in F $, so w"are $ X = \bigcup_{x \in X} U_x \not\in F $.
\par\bigskip\noindent
{\bf Beispiel 6} Der Umgebungsfilter eines Punktes ist ein irreduzibler
und insbe\-son\-dere ein quasikompakter Filter, da, wenn
$ x \in \bigcup_{i \in I} U_i $ gilt, $ x \in U_i $ f"ur ein $i$ sein mu\3.
\par\smallskip\noindent
{\bf Beispiel 7} Ein Ultrafilter $F$ ist quasikompakt. Ist n"amlich
$ V_1 \not\in F $ und $ V_2 \not\in F $, so gibt es nach der Charakterisierung
von Ultrafiltern $ U_1 , U_2 \in F $ mit $ V_1 \cap U_1 = \emptyset $ und
$ V_2 \cap U_2 = \emptyset $. Dann gilt
$ (V_1 \cup V_2 ) \cap U_1 \cap U_2 = \emptyset $ und daher
$ V_1 \cup V_2 \not\in F $.
\par\smallskip\noindent
{\bf Beispiel 8} Der Filter einer Folge $ x_n , n \in {\bf N} $, ist
in der Regel nicht quasi\-kom\-pakt.
\par\bigskip
Die Bezeichnungen quasikompakt und irreduzibel erkl"aren sich durch die
folgenden Aussagen. Sie geben
zugleich Beispiele, wie topologische Eigen\-schaf\-ten durch topologische
Filter beschrieben werden k"onnen.
\begin{satz}
Sei $X$ ein topologischer Raum. Dann sind "aquivalent:
\par\smallskip\noindent
(1) $X$ ist quasikompakt.
\par\smallskip\noindent
(2) Jeder quasikompakte Filter $F$ konvergiert, d.h. es gibt einen Punkt
$ x \in X $ mit $ F \supseteq U(x) $.
\end{satz}
Beweis. Aus (1) folgt (2). Sei $F$ ein quasikompakter Filter, und sei
angenom\-men, da\3 $F$ nicht konvergiert. Dann gibt es zu jedem Punkt $ x \in X$
eine offene Umgebung $ x \in U_x $ mit $ U_x \not\in F $. Diese $ U_x $
"uberdecken den quasikompakten Raum $X$, daher gibt es eine endliche
Teil"uberdeckung $ \bigcup _{i \in I } U_i = X \in F $. Da die Indexmenge
endlich ist und $F$ quasikompakt, geh"ort ein $U_i $ zu $F$, im Widerspruch zur
Voraussetzung. Die ``Quasikompaktheitseigenschaft'' von $F$ gilt ja nicht nur
f"ur Vereinigungen von zwei offenen Mengen, sondern f"ur beliebige endliche 
Vereinigungen, wie man durch Induktion sofort sieht.
\par\noindent
Aus (2) folgt (1). Sei $ X $ nicht quasikompakt. Dann gibt es eine offene
"Uberdeckung $ X = \bigcup _{i \in I} U_i $ ohne endliche Teil"uberdeckung. Zur
Konstruktion eines nicht konvergenten, quasikompakten Filters ziehen wir das
Lemma 1.1.3 "uber die Existenz von Filtern heran. Sei hierzu
\begin{displaymath}
S := \{ V \subseteq X \mbox{ offen } : V \supseteq \bigcup_{i \in J} U_i ,
\mbox{ mit } J \subseteq I \mbox{ kofinit } \}
\end{displaymath}
und
\begin{displaymath}
{\bf a} := \{ V \subseteq X \mbox{ offen } : V \subseteq \bigcup_{i \in H} U_i , 
\mbox{ mit } H \subseteq I \mbox{ endlich } \}.
\end{displaymath}
$S$ ist ein Filter, da aus $ V_1 \supseteq \bigcup_{i \in J_1 } U_i $ und
$ V_2 \supseteq \bigcup_{i \in J_2} U_i $, wobei $J_1 $ und $J_2 $ jeweils
fast alle Elemente aus $I$ enthalten, folgt
\begin{displaymath}
V_1 \cap V_2 \supseteq \bigcup_{i \in J_1}U_i \cap \bigcup_{i \in J_2} U_i
\supseteq \bigcup _{i \in J_1 \cap J_2} U_i ,
\end{displaymath}
wobei $ J_1 \cap J_2 $ wieder fast alle Indices enth"alt, und daher 
$ V_1 \cap V_2 $ zu $F$ geh"ort. ${\bf a}$ enth"alt mit einer Menge auch jede
kleinere und ist daher ein `Pseudoideal'. Ferner sind $S$ und $ {\bf a} $
disjunkt, denn angenommen es w"are $ V \supseteq \bigcup _{i \in J} U_i $ und
$ V \subseteq \bigcup_{i\in H} U_i $ mit $J$ kofinit und $H$ endlich. Dann
w"are $ \bigcup_ {i\in J} U_i  \subseteq \bigcup_{i \in H} U_i $ und $H$ lie\3e
sich mittels dem Komplement von $J$ zu einer endlichen Menge erg"anzen, die
eine endliche Teil"uberdeckung indiziert, im Widerspruch zur Wahl der $U_i$.
Nach dem zitierten Lemma gibt es also einen Filter $F$, der $S$ umfa\3t und zu
${\bf a} $ disjunkt ist mit der Eigenschaft, da\3 eine offene Menge $ V $ genau
dann nicht zu $F$ geh"ort, wenn es
ein $ U \in F $ gibt mit $ V \cap U \in {\bf a} $. Daraus folgt nun, da\3 $F$
quasikompakt ist: Sind n"amlich $V_1 , V_2 \not\in F $ und $ U_1 , U_2 \in F $
Mengen mit $ V_i \cap U_i \in {\bf a} $, $ i = 1,2 $, so ist
\begin{displaymath}
(V_1 \cup V_2) \cap (U_1 \cap U_2) \subseteq (V_1 \cap U_1) \cup (V_2 \cap U_2)
\in {\bf a} ,
\end{displaymath}
da $ {\bf a} $ stabil unter endlichen Vereinigungen ist. Damit geh"ort
$ V_1 \cup V_2 $ nicht zu $F$. $F$ ist aber nicht konvergent, da jeder Punkt
$ x \in X $ in einer der Mengen $ U_i $ der "Uberdeckung liegt, all diese
Mengen aber zu $ {\bf a} $ geh"oren und damit nicht zu $F$.
\par\bigskip
Vor dem n"achsten Satz erinnern wir daran, da\3 eine Teilmenge $T \subseteq X$
eines
topologischen Raumes irreduzibel hei\3t, wenn eine "Uberdeckung
$ T \subseteq A_1 \cup A_2 $ mit abgeschlossenen Mengen
$A_1 $ und $A_2 $
nur m"oglich ist,
wenn bereits $ T \subset A_i $ gilt f"ur ein $i$. Irreduzible Teilr"aume und
irreduzible Filter stehen folgenderma\3en in Beziehung.
\begin{satz}
Sei $X$ ein topologischer Raum. Ist $ A\subseteq X$ eine nichtleere,
irredu\-zi\-ble Teilmenge, so ist
$ F := \{ U \mbox{ offen } : U \cap A \neq \emptyset \} $ ein irreduzibler
Filter. Dabei definiert eine nichtleere, irreduzible Teilmenge denselben
irreduziblen Filter wie ihr ebenfalls irreduzibler Abschlu\3. 
Ist umgekehrt $F$ ein irreduzibler Filter, so ist
$ A := X - \bigcup _{ U \not\in F} U $ eine abgeschlossene, irreduzible
Teilmenge.
\end{satz}
Beweis. Sei $A$ eine nichtleere, irreduzible Teilmenge und $F$ die oben
definierte Menge. Dann ist $ X \in F $ , $F$ ist teilerstabil, und ist
$ U_1 \cap A \neq \emptyset $ und $ U_2 \cap A \neq \emptyset $, so ist mit
$ A_i := X-U_i , i = 1,2 $ , $ A \not\subseteq A_i , i = 1,2 $, und wegen der
Irreduzibilit"at von $A$ ist damit $ A \not\subseteq A_1 \cup A_2 $, und durch
Komplementbildung erh"alt man $ U_1 \cap U_2 \in F $. $F$ ist also ein Filter.
Ist nun $ \bigcup_{i\in I} U_i \cap A \neq \emptyset $, so ist auch
$ U_i \cap A \neq \emptyset $ f"ur ein $i$, was die Irreduzibilit"at von $F$
beweist.
\par\noindent
Sei umgekehrt $F$ ein irreduzibler Filter, $V := \bigcup _{U \not\in F} U $
und $A$ das Komplement von $V$. Dann gilt f"ur eine offene Menge $U$ :
\begin{displaymath}
U \not\in F \mbox{ genau dann, wenn }  U \subseteq V .
\end{displaymath}
Wegen der Irreduzibilit"at von $F$ geh"ort n"amlich $V$ nicht zu $F$, und daher
auch keine kleinere Menge, und eine offene Menge, die nicht zu $F$ geh"ort,
kommt als Bestandteil der Vereinigung vor und ist
deshalb in ihr enthalten. Sei nun $ A \subseteq A_1 \cup A_2 $ eine
"Uberdeckung mit abgeschlossenen Teilmengen. Dann ist
$ (X-A_1) \cap (X-A_2) \subseteq V $, also ist
$ (X-A_1) \cap (X-A_2) \not\in F $ und daher $ (X-A_i) \not\in F $ f"ur $ i=1 $
oder $ =2$. Dann ist wiederum $ X-A_i \subseteq V $, also $ A \subseteq A_i $,
was die Irreduzibilit"at von $A$ beweist.
\par\bigskip\noindent
{\bf Bemerkung} Wie oben kann man jedem Filter $F$ die abgeschlos\-se\-ne
Men\-ge $A := X - \bigcup_{U \not\in F} U $ zu\-ord\-nen. Die\-se ist gleich der
Men\-ge der Kon\-ver\-genzpunkte von $F$. Ist n"amlich $x \not\in A$, also
$ x \in U$ f"ur eine offene Menge $U$, die nicht zu $F$ geh"ort, so kann $x$
kein Konvergenzpunkt sein, und ist $x$ kein Konvergenzpunkt, so gibt es eine
offene Menge $U$, $x \in U$, $U \not\in F$, und $x$ geh"ort nicht zu $A$.
Eine irreduzible, abgeschlossene
Menge besitzt genau dann einen generischen Punkt, wenn der zugeh"orige
irreduzible Filter ein Umgebungsfilter ist. Eine offene Menge
n"amlich, die mit dem Abschlu\3 eines Punktes einen nichtleeren Schnitt hat,
umfa\3t bereits diesen Punkt, und ist $ F= U(x) $, so ist
$ \bar{\{ x \} } = X - \cup_{x \not\in U} U $.
\par\bigskip
Auch noethersche R"aume, also R"aume, in denen jede offene Teilmenge
quasikompakt ist, lassen sich mit topologischen Filtern charakterisieren.
\begin{satz}
F"ur einen topologischen Raum $X$ sind "aquivalent:
\par\smallskip\noindent
(1) $X$ ist noethersch.
\par\smallskip\noindent
(2) Jeder konsistente, quasikompakte Filter ist irreduzibel.
\end{satz}
Beweis. Aus (1) folgt (2). Sei $F$ ein konsistenter, quasikompakter
Filter im noetherschen
Raum $X$ und $ \bigcup_{i \in I} U_i = U \in F $. $I$ ist nicht leer, da $F$
konsistent ist. $U$ ist quasikompakt, daher
gibt es eine endliche Teil"uberdeckung von $ U_i , i \in I $, und aus der
Quasikompaktheit von $F$ ergibt sich $ U_i \in F $ f"ur ein $ i \in I $.
\par\noindent
Aus (2) folgt (1). Der Beweis erfordert wieder das Lemma 1.1.3 , und ist
"uberhaupt dem Beweis des Satzes 3.1.1 sehr "ahnlich. Sei $X$ als
nicht noethersch angenommen, es gebe also eine nicht quasikompakte,
offene Teil\-men\-ge $ U$ und somit eine offene "Uberdeckung
$ U = \bigcup_{i \in I} U_i $ ohne endliche Teil\-"uber\-deck\-ung. Wir setzen
\begin{displaymath}
S := \{ V \subseteq X \mbox{ offen } :
V \cap U \supseteq \bigcup_{i \in J} U_i \mbox{ f"ur ein kofinites }
J \subseteq I \} 
\end{displaymath}
und wieder
\begin{displaymath}
{\bf a} := \{ V \mbox{ offen } : V \supseteq \bigcup _{i \in H} U_i
\mbox{ mit endlichem } H \}.
\end{displaymath}
$ S$ ist ein Filter, der $U$ enth"alt, ${ \bf a} $ ist abgeschlossen
unter Teilmengen
und unter endlichen Vereinigungen, und $S$ und $ {\bf a} $ sind wieder disjunkt,
da sich sonst eine endliche Teil"uberdeckung von $U$ ergebe. Nach dem
angef"uhrten Lemma gibt es dann einen Filter $F$ mit $ S \subseteq F $ und
$ F \cap {\bf a} = \emptyset $ und der Eigenschaft, da\3 eine offene Menge $W$
genau dann nicht zu $F$ geh"ort, wenn es ein Element $ V \in F $ gibt mit
$ W \cap V \in {\bf a} $. Die gleiche Argumentation wie oben zeigt, da\3 $F$
quasikompakt ist. $F$ ist aber nicht irreduzibel. F"ur alle $ i \in I $ ist
$ U_i \in {\bf a} $, also $ U_i \not\in F $, aber
$\bigcup_{i \in I} U_i = U \in F $.
\par\bigskip
Eine weitere Anwendung des Existenzlemmas f"ur Filter ist folgende Aus\-sa\-ge,
die ein topologisches Analogon zu der Aussage ist, da\3 in einem kom\-mu\-ta\-ti\-ven
Ring jeder Filter Durchschnitt von Komplementen von Prim\-ide\-a\-len ist.
\begin{lem}
Jeder topologische Filter ist Durchschnitt von quasikompak\-ten Filtern.
\end{lem}
Beweis. Sei $F$ ein Filter des topologischen Raumes $X$. Die Behauptung ist
$ F = \bigcap_{F \subseteq G,\hspace{0.5em} G \mbox{ \rm quasikompakt} } G $.
Dabei ist die eine Inklusion trivial, zum Beweis der anderen sei
$ U \not\in F $, und ${\bf a} $ sei die Menge aller offenen Teilmengen von $U$.
Dann gibt es zu $F$ einen Oberfilter $G$, disjunkt zu ${\bf a}$, also
$ U \not\in G $, mit der Eigenschaft $ W \not\in G $ genau dann, wenn es ein
$V \in G $ gibt mit $ W \cap V \subseteq U $. Daraus folgt dann wieder sofort,
da\3 $G$ quasikompakt ist, und $U$ geh"ort nicht zum Durchschnitt aller $F$
umfassenden quasikompakten Filter.
\par\bigskip
Mit den bisherigen "Uberlegungen ergibt sich nun eine "uberraschend ein\-fa\-che
Beschreibung der Filter in einem Zariskiraum. Ein Zariskiraum ist
ein noe\-ther\-scher Raum $X$, in dem jede irreduzible abgeschlossene Teilmenge
genau einen generischen Punkt hat, (\cite{har}, p.93). In unserer
Terminologie han\-delt es sich um einen noetherschen Raum, in dem jeder
irreduzible Filter Umge\-bungsfilter genau eines Punktes ist, wobei hier die
Eindeutig\-keits\-aus\-sage zum $T_0$ Axiom "aqui\-va\-lent ist. Beispielsweise
ist der zugrunde liegende Raum eines noether\-schen Schemas ein Zariskiraum.
Die Aussage lautet:
\begin{kor}
Sei $X$ ein noetherscher Raum, in dem jeder irreduzible Filter ein
Umgebungsfilter eines Punktes ist. Dann ist jeder Filter $F$ Umgebungsfil\-ter
einer Teilmenge $ T \subseteq X $, und zwar ist
$ T := \{ x \in X : F \subseteq U(x) \} $
\end{kor}
Beweis. Sei $F$ ein Filter. Nach dem vorangegangenen Lemma ist $F$ der Durch\-schnitt
aller ihn umfassenden quasikompakten Filter. Diese sind aber aufgrund der
Voraussetzungen gleich den
irreduziblen und damit gleich den Umge\-bungs\-filtern oberhalb von $F$, also gleich
den Umgebungsfiltern zu Punk\-ten aus T. Also ist
$F = \bigcap_{x \in T} U(x) = U(T) $.
\section{Filter unter stetigen Abbildungen}
Zu einer stetigen Abbildung
$ \phi : X \longrightarrow Y $ zwischen topologischen R"aumen geh"ort in
funktorieller Weise der Monoidmorphismus
$ \varphi : {\rm Top}(Y) \longrightarrow {\rm Top}(X) $
mit $ \varphi(V) = \phi^{-1}(V) $. Dazu geh"oren wiederum die in Kapitel 2
entwickelten ko- bzw. kontravarianten stetigen Abbildungen $ \varphi_\ast $
und $ \varphi^\ast $, denen wir uns nun zuwenden. Dabei ergibt sich vieles aus
den allgemeinen "Uberlegungen, es treten aber auch
Eigenschaften und Beziehungen auf, die im abstrakten Fall eines
Monoidmorphismus keine Entsprechung haben beziehungsweise gar\-nicht definiert
werden k"onnen. 
\par
Sei zun"achst 
$ \varphi^\ast : \mbox{ Filt}X \longrightarrow \mbox{ Filt}Y $ die
kontravariante Abbildung zum Monoidmorphismus
$ \varphi : \mbox{ Top}(Y) \longrightarrow \mbox{ Top}(X) $.
Von der stetigen Abbildung $ \phi : X \longrightarrow Y $ aus gesehen handelt es
sich nat"urlich um eine kovariante, und wir bezeichnen sie wieder mit
$ \phi $, was aufgrund der folgenden Proposition (1) nicht zu
Mi\3verst"andnissen f"uhren d"urfte. Einem Filter $F$ in $X$ wird also unter
$\phi $ der Filter
\begin{displaymath}
\phi(F) = \{ V \in {\rm Top}(Y) : \varphi(V) = \phi^{-1} (V) \in F \}
\end{displaymath}
in $Y$ zugeordnet. Folgende Aussagen fassen einige Ergebnisse zusammen, wie
sich Filter unter dieser Abbildung verhalten.
\begin{pro}
Sei $ \phi : X \longrightarrow Y $ eine stetige Abbildung, und $F$ ein Filter
in $X$. Dann gilt
\par\smallskip\noindent
(1) Ist $ T \subseteq X $ eine beliebige Teilmenge, so ist
$ \phi(U(T)) = U(\phi(T)) $. Insbeson\-dere ist $\phi(U(x)) = U (\phi(x)) $
f"ur jeden Punkt $x \in X $.
\par\smallskip\noindent
(2) Ist $F$ quasikompakt, so ist auch $ \phi(F) $ quasikompakt.
\par\smallskip\noindent
(3) Ist $F$ irreduzibel, so ist auch $ \phi(F) $ irreduzibel.
\end{pro}
Beweis. (1) ergibt sich direkt aus
\begin{eqnarray*}
\phi(U(T)) = \{ V \in {\rm Top}(Y) : \phi^{-1} (V) \supseteq T \} \\
= \{ V \in {\rm Top}(Y) : V \supseteq \phi(T) \}
= U(\phi(T)) .
\end{eqnarray*}
(3) Sei $F$ irreduzibel und $ \bigcup_{i \in I} U_i \in \phi(F) $. Dann gibt
es wegen
$ \phi^{-1}( \bigcup_{i \in I} U_i ) = \bigcup_{i \in I} \phi^{-1}(U_i) $
und wegen der Irreduzibilit"at von $F$ ein $i \in I $ mit
$ \phi^{-1}(U_i) \in F $, also
$ U_i \in \phi(F) $. (2) wird wie (3) bewiesen, wobei man sich auf endliche
Indexmengen beschr"ankt.
\par\bigskip\noindent
{\bf Bemerkung} Im Hinblick auf eine sp"atere "Uberlegung sei bemerkt, da\3
die Aussage (3) der Proposition auch dann gilt, wenn nur abstrakt eine
Abbil\-dung $ \varphi : {\rm Top}(Y) \longrightarrow {\rm Top}(X) $ gegeben
ist, die endliche
Durchschnitte und be\-lie\-bige Vereinigungen respektiert, ohne da\3 sie durch eine
stetige Abbildung in\-du\-ziert ist.
\par\bigskip
Das Filtrum eines topologischen Raumes ist in nat"urlicher Weise eine
Erweiterung des Ausgangsraumes, und stetige Abbildungen lassen sich auf die
Filtren ausdehnen.
\begin{satz}
Sei $X$ ein topologischer Raum. Die Abbildung
\begin{displaymath}
i: X \longrightarrow \mbox{\rm Filt}X  \mbox{ mit } i(x) := U(x)
\end{displaymath}
ist stetig, $X$ tr"agt die Initialtopologie bez"uglich dieser Abbildung und
das Bild $ i(X) $ ist dicht in {\rm Filt}$^oX $. Ist $X$ ein $ T_0- $ Raum, so ist
die Abbildung injektiv und eine Einbettung. Ist $ \phi : X \longrightarrow Y $
eine stetige Abbildung, so ist das Diagramm
\par\noindent
\unitlength0.5cm
\begin{picture}(10,7)
\put(4,1){\vector(1,0){3}}
\put(4,5){\vector(1,0){3}}
\put(2.8,4){\vector(0,-1){2}}
\put(8.2,4){\vector(0,-1){2}}
\put(2.6,4.7){$X$}
\put(8,4.7){$Y$}
\put(1.8,0.8){{\rm Filt}$X$}
\put(7.5,0.8){{\rm Filt}$Y$}
\put(2.3,3){i}
\put(7.5,3){i}
\put(5.5,5.4){$\phi$}
\put(5.5,1.5){$\phi$}
\end{picture}
\par\noindent
kommutativ.
\end{satz}
Beweis. Es ist $ i^{-1}(D(U)) = U $ f"ur eine offene Menge $ U \subseteq X $,
das zeigt die Stetigkeit, und da\3 $X$ die Initialtopologie tr"agt. Das
Trennungsaxiom $ T_0 $ besagt gerade, da\3 umgebungsgleiche Punkte gleich sind,
und ist daher "aquivalent zur Injektivit"at der Abbildung. Zur Dichtheit
beachte man, da\3 $D(U) = \emptyset $ in Filt$^oX$ genau dann gilt, wenn
$ U = \emptyset $; ist also $ D(U) \neq \emptyset $, so ist $U$ nicht leer und
es gibt einen Punkt in $U$, dessen Umgebungsfilter in $D(U)$ liegt. Die
Kommutativit"at des Diagrammes besagt einfach $U(\phi(x))= \phi(U(x))$, was
in der vorangegangenen Proposition gezeigt wurde.
\par\bigskip
Betrachtet man im obigen Diagramm statt der Abbildung
$ \phi : \mbox{Filt}X \longrightarrow \mbox{Filt}Y $ die sich aus dem
Monoidmorphismus $ \varphi : {\rm Top}(Y) \longrightarrow {\rm Top}(X) $
ergeben\-de kovariante Abbildung
\begin{displaymath}
\phi^{-1} := \varphi_\ast : \mbox{Filt}Y \longrightarrow \mbox{Filt}X
\end{displaymath}
mit $ \phi^{-1}(G) = F(\{ \phi^{-1}(V) : V \in G \} ) $, so kommutiert das
entsprechende Diagramm in der Regel nicht. Es gelten auch nicht die
entsprechenden Aussagen der Proposition 3.2.1 : Ist etwa $ S \subseteq Y $ eine Teilmenge,
so gilt zwar stets $ \phi^{-1} (U(S)) \subseteq U(\phi^{-1}(S)) $, da ja
$ U \in \phi^{-1}(U(S)) $, also $ U \supseteq \phi^{-1} (V) $ f"ur ein
$ V \supseteq S $ sofort $ U \supseteq \phi^{-1}(S) $ impliziert, die umgekehrte
Inklusion gilt aber nicht. F"ur einen Punkt $ x \in X $ mit $ \phi(x) = y $ ergibt
sich speziell $ U(x) \supseteq U(\phi^{-1}(y)) \supseteq \phi^{-1}(U(y)) $.
Die G"ultigkeit der anderen Inklusionen charakterisiert Eigenschaften stetiger
Abbildungen.
\begin{satz}
Sei $ \phi :X \longrightarrow Y $ eine stetige Abbildung.
Dann sind "aquivalent:
\par\smallskip\noindent
(1) $X$ tr"agt
die Initialtopologie bez"uglich $ \phi $, es sind also gerade die Mengen
$ \phi^{-1}(V) , V \in {\rm Top}(Y) $ offen.
\par\smallskip\noindent
(2) Alle topologischen Filter auf $X$ sind fix.
\par\smallskip\noindent
(3) Die Umgebungsfilter $U(x)$ sind Fixfilter f"ur alle $ x \in X $.
\par\smallskip\noindent
(4) Das folgende Diagramm kommutiert.
\par\noindent
\unitlength0.5cm
\begin{picture}(10,7)
\put(7,1){\vector(-1,0){3}}
\put(4,5){\vector(1,0){3}}
\put(2.8,4){\vector(0,-1){2}}
\put(8.2,4){\vector(0,-1){2}}
\put(2.6,4.7){$X$}
\put(8,4.7){$Y$}
\put(1.8,0.8){{\rm Filt}$X$}
\put(7.5,0.8){{\rm Filt}$Y$}
\put(2.3,3){i}
\put(7.5,3){i}
\put(5.5,5.5){$\phi$}
\put(5.5,1.5){$\phi^{-1}$}
\end{picture}
\par\noindent
Ist $X$ zus"atzlich ein $T_0 $ Raum, so ist $ \phi $ injektiv, also eine
Einbettung.
\end{satz}
Beweis. (3) besagt $ U(x) = \phi^{-1} \phi(U(x)) $ und (4) besagt
$ U(x) = \phi^{-1}U(\phi(x)) $, jeweils f"ur alle $x \in X $, was wegen
$ \phi(U(x)) = U(\phi(x)) $ dasselbe ist. Aus (1) folgt (2).
Die Voraussetzung bedeutet, da\3 der Monoidmorphismus
$ V \longmapsto \phi^{-1}(V) $ surjektiv ist, also sind die Filter des
Monoids Top$(X) $ fix. Die Implikation von (2) nach (3) ist eine Einschr"ankung.
Aus (3) folgt (1). Sei $U \subseteq  X $ offen
und ohne Einschr"ankung nicht leer. F"ur $ x \in U $ gibt es dann nach
Voraussetzung ein $ V_x \in U(\phi(x)) $ offen mit
$x \in \phi^{-1}(V_x) \subseteq U $.
Somit ist $U$ die Vereinigung der $ \phi^{-1}(V_x) , x \in U $.
Ist $X$ ein $T_0$-Raum,
so ist $ i: X \longrightarrow \mbox{Filt}X $ injektiv, also auch
$ \phi : X \longrightarrow Y $.
\begin{satz}
Sei $ \phi : X \longrightarrow Y $ eine stetige Abbildung. Dann sind "aquivalent:
\par\smallskip\noindent
(1) $ \phi $ ist abgeschlossen.
\par\smallskip\noindent
(2) F"ur jede Teilmenge $ S \subseteq Y $ ist $ \phi^{-1}U(S) = U(\phi^{-1}(S))$.
\par\smallskip\noindent
(3) F"ur jeden Punkt $ y \in Y $ ist $ \phi^{-1}U(y)  = U(\phi^{-1}(y)) $.
\end{satz}
Beweis. Aus (1) folgt (2). Schon oben wurde bemerkt, da\3
$ \phi^{-1} U(S) \subseteq U(\phi^{-1}(S)) $ stets gilt, sei also
$ \phi^{-1}(S) \subseteq U $ offen. Sei $A := X-U$, dann ist nach Voraussetzung
$ V := Y- \phi(A) $ offen.
Aus $ x \not\in U $ ergibt sich $\phi(x) \not\in S $, es sind also $\phi(A)$
und $S$ disjunkt, also $ S \subseteq V $. Ferner ergibt sich aus $ x \not\in U $
sofort $ \phi(x) \in \phi(A) $, also $\phi(x) \not\in V $, d.h.
$ U \supseteq \phi^{-1}(V) $, und damit ist $ U \in \phi^{-1}U(S) $.
(3) ist ein Abschw"achung von (2), bleibt also noch (3) impliziert (1) zu
zeigen. Sei $ A \subseteq X $ abgeschlossen und somit ist $ W := Y - \phi(A) $
als offen nachzuweisen; sei hierzu $ y \in W $. Dann ist
$ \phi^{-1}(y) \subseteq X-A $, da es andernfalls ein $ x \in A $ gebe mit
$ \phi(x) = y $. Es ist also $ X-A \in U(\phi^{-1}(y)) $ und nach
Voraussetzung gibt es ein $ V \in U(y) $ mit $ \phi^{-1}(V) \subseteq X-A $.
Dann ist aber $ y \in V \subseteq Y - \phi(A) = W $, also ist $W$ offen.
\par\bigskip
H"aufig taucht folgendes Problem auf: man hat eine stetige Abbildung auf
einer Teilmenge eines topologischen Raumes definiert und fragt sich, inwiefern
diese Abbildung auf den gesamten Raum ausgedehnt werden kann. Eine Variante
davon ergibt sich, wenn man nicht mehr von einem festen Oberraum des
Definitionsbereiches ausgeht, sondern nach einem m"oglichst gro\3en
Erweiterungsraum sucht, auf dem sich die Abbildung fortsetzen l"a\3t.
Der Bildraum wird im folgenden als regul"ar vorausgesetzt, d.h. er ist ein $T_0$
Raum und jede Umgebung eines Punktes umfa\3t eine abgeschlossene Umgebung dieses
Punktes. Metrische R"aume und kompakte R"aume sind beispielsweise regul"ar.
Zu einem topologischen Raum $Y$ sei mit Kon$Y$ die Menge der konvergenten,
konsistenten Filter mit der induzierten Topologie bezeichnet. Zum
folgenden vergleiche \cite{que}, Kap 6.C.
\begin{lem}
Sei $Y$ ein regul"arer topologischer Raum. Dann ist die Kon\-ver\-genz\-abbildung
$ k: \mbox{\rm Kon}Y \longrightarrow Y $, die jedem kon\-ver\-gen\-ten Filter seinen
Kon\-ver\-genz\-punkt
zuordnet, stetig.
\end{lem}
Beweis. Da $Y$ insbesondere Hausdorffsch ist, gibt es zu jedem konvergenten,
konsisten Filter genau einen Limespunkt, die Abbildung ist daher wohldefi\-niert.
Sei $V$ eine offene Menge in $Y$, und $G \in k^{-1}(V)$. Sei
$y := k(G) \in V $ der Konvergenzpunkt. Aufgrund der Regularit"at ist
$ y \in U \subseteq A \subseteq V $ mit $U$ offen und $A$ abgeschlossen. Wir
zeigen $ G \in D(U) \subseteq k^{-1}(V) $. Aus $ y \in U $ und
$ U(y) \subseteq G $ folgt $ U \in G $. Sei $ H \in D(U) $ ein anderer
konsistenter, konvergenter Filter mit Konvergenzpunkt $z$ und sei
$ z \not\in V $ angenommen. Dann ist erst recht $ z \not\in A $, also
$ z \in Y-A $ und damit $ Y-A \in H $. Damit enth"alt $H$ die disjunkten offenen
Mengen $U$ und $ Y-A $ im Widerspruch zur Konsistenz.
\begin{satz}
Sei $ \phi : X \longrightarrow Y $ eine stetige Abbildung, $Y$ regul"ar.
Es sei {\rm Kon}$_\phi X := \phi^{-1}(\mbox{\rm Kon}Y) $
das Urbild der konvergenten Filter
auf $Y$ unter der zugeh"origen Filtrumsabbildung
$\phi : \mbox{\rm Filt}X \longrightarrow \mbox{\rm Filt}Y $. Dann ist
$ \tilde{\phi} : \mbox{\rm Kon}_{\phi}X \longrightarrow Y $ mit
$ \tilde{\phi} = k \circ \phi $ eine stetige Fortsetzung
von $ \phi $. Ist
\par\noindent
\unitlength0.5cm
\begin{picture}(10,7)
\put(4,5){\vector(1,0){3}}
\put(2.8,4){\vector(0,-1){2}}
\put(4,1){\vector(1,1){3}}
\put(2.6,4.7){$X$}
\put(8,4.7){$Y$}
\put(2.6,0.8){$Z$}
\put(2.3,3){j}
\put(5.5,5.5){$\phi$}
\put(6,2){$\bar{\phi}$}
\end{picture}
\par\noindent
ein kommutatives Diagramm mit stetigen Abbildungen und mit $j(X) $ dicht in $Z$,
und ist $ \psi $ die stetige Abbildung
$ Z \stackrel{ i}{\longrightarrow} \mbox{Filt}Z \stackrel { j^{-1} }
 {\longrightarrow} \mbox{Filt}X $,
so ist $\psi(Z) \subseteq Kon_\phi X $, und das Diagramm
\par\noindent
\unitlength0.5cm
\begin{picture}(10,7)
\put(4,5){\vector(1,0){3}}
\put(2.8,4){\vector(0,-1){2}}
\put(4,1){\vector(1,1){3}}
\put(2.6,4.7){$Z$}
\put(1.4,0.8){{\rm Kon}$_\phi X $ }
\put(8,4.7){$Y$}
\put(2.1,3){$\psi$}
\put(5.5,5.5){$\bar{\phi}$}
\put(6,2){$ \tilde{\phi}$}
\end{picture}
\par\noindent
kommutiert. In diesem Sinne ist also
$ \tilde{\phi} : Kon_\phi X \longrightarrow Y $
die einzige wesent\-li\-che Fortsetzung von $ \phi : X \longrightarrow Y $ in
regul"ares $Y$ auf einen Oberraum, in dem $X$ dicht ist.
\end{satz}
Beweis. $ \tilde{\phi} = k \circ \phi $ ist stetig, und es ist f"ur $ x \in X $:
\begin{displaymath}
k \circ \phi \circ i(x) = k (U(\phi(x))) = \phi(x) ,
\end{displaymath}
$\tilde{\phi} $ ist also eine Fortsetzung von $ \phi $. Sei nun
$ Z , \bar{\phi} $ wie im Satz beschrieben. Es ist f"ur $z \in Z $ zu zeigen, da\3
$ F := \psi (z) $ zu Kon $_\phi X $ geh"ort, und  da\3 das Bild von $ F $
unter $\phi$ gegen $ \bar{\phi}(z) $ konvergiert. Es ist
$ F = F \{ j^{-1}(W) : z \in W \} $, und da $j$ dichtes Bild in $Z$ hat ist
$F$ konsistent, und damit ist auch $ \phi (F) $ konsistent, und es ist 
lediglich $ U( \bar{\phi}(z) ) \subseteq \phi (F) $ zu zeigen. Sei hierzu
$ \bar{\phi}(z) \in V $, $V$ offen in $Y$. Dann ist
$ \bar{\phi}^{-1}(V) \in U(z) $ und damit
$ j^{-1}(\bar{\phi}^{-1}(V)) = \phi^{-1}(V) \in F $. Also ist $V \in \phi(F) $.
\section{Filtrumserweiterungen}
In diesem und im n"achsten Abschnitt steht die zu einer stetigen Abbildung
$ \phi : X \longrightarrow Y $ geh"orende stetige Abbildung
$ \psi : Y \longrightarrow \mbox{ Filt}X $, die sich durch
Hintereinanderschaltung der Abbildungen $ i: Y \longrightarrow \mbox{ Filt}Y $
und $ \phi^{-1} : \mbox{ Filt}Y \longrightarrow \mbox{ Filt}X $ ergibt,
im Mittelpunkt der Untersuchung. Unter bestimm\-ten Voraussetzungen
stiftet $ \psi $
eine Hom"oomorphie zwischen $Y$ und dem Bild
$ \psi(Y) \subseteq \mbox{ Filt}X $,
und unter weiteren Bedingungen k"onnen topologische Erwei\-ter\-ungen $Y$ von $X$
mit geeigneten Filtern und der induzierten
Filtrumstopo\-lo\-gie be\-schrie\-ben werden.
Solche Erweiterungen hei\3en Filtrumserweiterungen; bei diesen ergeben sich
Eigenschaften wie Dichtheit und Einbettung, aber auch, wenn die Konstruktion
funktoriell ist, die erweiterten Morphismen automa\-tisch aus den Eigenschaften
des Filtrums des Ausgangs\-rau\-mes. "Ahn\-liche "Uberlegungen, wie man 
topologische Erweiterungen systematisch mit Filtern beschreiben kann,
finden sich in
\cite{csa}, Kap. 6, allerdings ohne den kanonischen Rahmen, den das Filtrum
uns zur Verf"ugung stellt.
\par
Zun"achst wird eine hinreichende Bedingung angegeben, die es gestattet, den
Bildraum als Unterraum des Filtrums aufzufassen.
\par\bigskip\noindent
{\bf Definition} Eine stetige Abbildung $ \phi : X \longrightarrow Y $ hei\3t
filterhaft, wenn sie folgender Bedingung gen"ugt:
\par\medskip\noindent
F"ur alle Punkte $ y \in Y $ und alle offenen Mengen $W$
mit
$ y \in W $ gibt es eine offene Menge $V$ mit
$ y \in V \subseteq W $, so da\3 f"ur alle offenen Teilmengen $ V' $
gilt:
Ist $ \phi^{-1}(V') \subseteq \phi^{-1}(V) $, so ist
$ V'\subseteq W $.
\par\medskip\noindent
$Y$ hei\3t eine Filtrumserweiterung von $X$, wenn $Y$ ein $ T_0 $-Raum ist,
und die Abbildung $ \phi $ eine dichte, filterhafte Einbettung ist.
\par\bigskip\noindent
{\bf Bemerkung} Es gen"ugt, die obige Bedingung f"ur offene Mengen $W$
und $ V'$ aus einer Basis zu formulieren beziehungsweise zu testen.
Ist $ \phi : X \longrightarrow Y $ filterhaft, und ist $ Z \subseteq Y $ ein
Unterraum, der das Bild von $ \phi $ umfa\3t, so ist auch
$ \phi : X \longrightarrow Z $ filterhaft. Dagegen ist die Verkn"upfung
zweier filterhafter Abbildungen in der Regel nicht wieder filterhaft.
\par\medskip\noindent
{\bf Beispiel 1} Die kanonische Abbildung
$ i : X \longrightarrow \mbox{ Filt}^\circ X $ ist, wenn $X$ ein $T_0 $-Raum ist,
eine Filtrumserweiterung. Es wurde schon weiter oben gezeigt, da\3 sie eine
dichte Einbettung ist. Die angegebene Eigenschaft zeigt man so: Sei
$ F \in D(U)=W=V $, also $ U \in F $, und sei $ V'= D(U') $. Aus der Beziehung
$ U' \subseteq U $ in $X$ ergibt sich sofort wie gefordert
$ D(U') \subseteq D(U) $.
Dasselbe gilt f"ur jeden Unterraum $Y$ des Filtrums mit
$ X \subseteq Y \subseteq \mbox{ Filt}^\circ X $.
\par\smallskip\noindent
{\bf Beispiel 2} Jede surjektive Abbildung $ \phi : X \longrightarrow Y $ ist
filterhaft, da bei einer surjektiven Abbildung aus
$ \phi^{-1}(V') \subseteq \phi^{-1} (V) $ generell $ V'\subseteq V $ folgt.
\par\smallskip\noindent
{\bf Beispiel 3} Sei $Y$ regul"ar und $ \phi $ eine dichte Einbettung. Dann ist
$ Y $ eine Filtrumserweiterung. Zun"achst sei bemerkt, da\3 generell
f"ur dichtes $ D \subseteq Y $ und offenes $ V \subseteq Y $ gilt
$ \bar{V} = \widehat{D \cap V} $. Dabei ist ``$ \supseteq $''  klar, sei also
$ y \in \bar{V} $, und es ist zu zeigen, da\3 jede offene Umgebung $ y \in U $
mit $ V \cap D $ einen nichtleeren Durchschnitt hat. Es ist aber nach
Voraussetzung $ U \cap V \neq \emptyset $ und daher
$ D \cap U \cap V \neq \emptyset $ wegen der Dichtheit von $D$.
\par
Sei nun $ y \in W $, $W$ offen. Dann gibt es aufgrund der Regularit"at von $Y$
eine offene Menge $V$ mit $ y \in V \subseteq \bar{V} \subseteq W $. Ist nun
$ V'$ eine offene Menge mit $ \phi^{-1} (V') \subseteq \phi^{-1} (V) $, also
mit $ X \cap V' \subseteq X \cap V $, so ist
$ \widehat{ V'\cap X} \subseteq \widehat{ V \cap X} = \bar{V} $
nach der Vor"uberlegung.
Ferner ist $ V'\subseteq \bar{V'} = \widehat{ V'\cap X} $, also insgesamt
$ V' \subseteq \bar{V} \subseteq W $, wie verlangt.
\begin{satz}
Ist $ \phi : X \longrightarrow Y $ filterhaft und $Y$ ein $T_0 $-Raum, so ist
die Abbildung $ \psi : Y \longrightarrow \mbox{ \rm Filt}X $ eine Einbettung. Ist
$ \phi $ zus"atzlich eine dichte Einbettung, dann ist auch $ \psi $ eine dichte
Einbettung in {\rm Filt}$ ^\circ X$ und das Dia\-gramm 
\par\noindent
\unitlength0.5cm
\begin{picture}(10,7)
\put(4,5){\vector(1,0){3}}
\put(2.8,4){\vector(0,-1){2}}
\put(7,4){\vector(-1,-1){3}}
\put(2.6,4.7){$X$}
\put(1.2,0.8){{\rm Filt}$ ^\circ X$}
\put(8,4.7){$Y$}
\put(2.3,3){i}
\put(5.5,5.5){$\phi$}
\put(6,2){$\psi$}
\end{picture}
\par\noindent
kommutiert. In diesem Fall ist $ \psi $ durch diese Eigenschaften eindeutig
be\-stimmt. Jede Filtrumserweiterung $ \phi , Y $ entspricht also einem
eindeutigen
Teilraum $ X \subseteq Y' \subseteq \mbox{ \rm Filt} ^\circ X $ samt der kanonischen
Abbildung i, und $ X \longrightarrow \mbox{ \rm Filt}^\circ X $ ist die
universelle Filtrumserweiterung.
\end{satz}
Beweis. Sei $ y \neq y' $. Dann gibt es aufgrund der $T_0$-Eigenschaft eine
offene Menge $W$ mit $ y \in W $ und $ y' \not\in W $. Sei $ y \in V \subseteq W$
so, wie durch die Filterhaftigkeit beschrieben. Dann ist
$ \psi(y) \in D(\phi^{-1}(V)) $. Sei angenommen, es w"are auch
$ \psi(y') \in D(\phi^{-1}(V)) $, also
$ \phi^{-1}(V) \in \psi(y') = F \{ \phi^{-1}(U) : y'\in U\} $. Das hei\3t es
gibt eine offene Umgebung $ V'$ von $y'$ mit $ \phi^{-1}(V') \subseteq \phi^{-1}(V) $. Nach der
Wahl von $V$ folgt aber daraus $ V'\subseteq W $, ein Widerspruch zu
$ y' \not\in W $. Dies zeigt die Injektivit"at.
\par
Zum Beweis der Einbettungseigenschaft ist zu vorgegebenen $ y \in W $ zu zeigen,
da\3 es eine offene Menge im Filtrum gibt, deren Urbild unter $ \psi $ in $W$
liegt und zu dem $y$ geh"ort. Sei wieder $ y \in V \subseteq W $. Dann leistet
$ D( \phi^{-1}(V)) $ das Gew"unschte. Wie beim Beweis der Injektivit"at folgt
n"amlich aus $ y' \not\in W $ wieder $ \psi(y') \not\in D(\phi^{-1}(V)) $, also
$ \psi^{-1}(D(\phi^{-1}(V))) \subseteq W $.
\par
Ist nun $ \phi $ zus"atzlich eine dichte Einbettung, so ist f"ur $ y \in Y $
der Filter $ \psi(y) $ konsistent, da in diesem Fall $ \phi^{-1}(U) $ f"ur
$ y \in U $ offen nicht leer ist. Da\3 das Diagramm kommutiert wurde f"ur eine
Einbettung $ \phi $ schon in 3.2.3 gezeigt. Die Dichtheit ergibt sich, da
bereits $X$ in Filt$^\circ X$ dicht ist und $ X \subset \psi(Y) $ ist.
\par
F"ur eine weitere Einbettung
$ \tilde{\psi} : Y \longrightarrow \mbox{ Filt}^\circ X $,
f"ur die das Diagramm kommutiert, ist zu zeigen, da\3 sie mit der kanonischen
Einbettung $ \psi $ zusam\-men\-f"allt. Hierzu ist zu zeigen, da\3 f"ur einen
Punkt $ y \in Y $ die Bildfilter $ \tilde{\psi}(y) $ und $ \psi(y) $
"ubereinstimmen.
Sei $ U \in \tilde{\psi}(y) $, also $ \tilde{\psi}(y) \in D(U) $. Wegen der 
Stetigkeit von $ \tilde{\psi} $ ist dann $ y \in \tilde{\psi}^{-1}(D(U)) =: V $
offen. Es ist $ \phi^{-1}(V) = i^{-1}(D(U)) = U $, da das Diagramm nach
Voraussetzung kommutiert, und damit ist $ U \in \psi(y) $.
\par
Umgekehrt gen"ugt es zu zeigen, da\3 f"ur $ y \in V $ gilt
$ \phi^{-1}(V) \in \tilde{\psi}(y) $, da die $ \phi^{-1}(V) $, $ V \in U(y) $, den
Filter $ \psi(y) $ erzeugen. Da $ \tilde{\psi} $ eine Einbettung ist, gibt es
eine offene Menge $ W \subseteq \mbox{ Filt}^\circ X $ mit
$ \tilde{\psi}^{-1}(W) = V $ und damit eine Basismenge $ D(U) $ mit
$ \tilde{\psi}(y) \in D(U) $ und $ y \in \tilde{\psi}^{-1}(D(U)) \subseteq V $.
Damit ist $ U = \phi^{-1} \tilde{\psi}^{-1}(D(U)) \subseteq \phi^{-1}(V) $,
also $ \phi^{-1}(V) \in \tilde{\psi}(y) $, da dies bereits f"ur $U$ gilt.
\par\bigskip
Im folgenden wird anhand einiger Beispiele gezeigt, durch welche topologi\-schen
Filter die Bildr"aume unter filterhaften Abbildungen beschrieben wer\-den. Die
direkte Beschreibung mittels Filtern kann sowohl Vor- als auch Nach\-teile haben.
Vorteile liegen in der Systematik und darin, da\3 gewisse, immer wiederkehrende
"Uberlegungen nicht st"andig wiederholt werden m"us\-sen; Nach\-teile treten
insbesondere dann auf, wenn die Konstruk\-ti\-on nicht nur eine topolo\-gische ist,
sondern zus"atzlich gewisse weitere Struk\-turen
(Metrik, Ver\-kn"upf\-ungen etc.) des
Ausgangsraumes "ubertragen werden sollen.
\par\bigskip\noindent
{\bf Beispiel 4} Auf einem topologischen Raum  $X$ sei eine "Aquivalenzrelation
$ \sim $ gegeben. $Y$ sei die Menge der "Aquivalenzklassen und mit der
Quotiententopo\-logie versehen, und $ \phi : X \longrightarrow Y $ sei die
kanonische, stetige, surjektive Ab\-bil\-dung. Eine Teilmenge von $Y$ ist also genau
dann offen, wenn das Urbild offen ist. Die offenen Urbilder sind dabei genau
die offenen Teilmengen von $X$, die mit der "Aquivalenzrelation vertr"aglich
sind, die also mit einem Punkt $x$ auch dessen "Aquivalenzklasse $ [x] $
enth"alt.
Damit ist
$ \psi(y) = F( \{ U \subseteq X \mbox{ offen und vertr"aglich mit } \sim :
[y] \subseteq U \} ) $, und $Y$ ist hom"oomorph zum Raum der 
Filter von diesem Typ.
Ist die Quotientenabbildung abgeschlossen, so sind diese Filter einfach die
Umgebungsfilter der "Aquivalenzklassen. 
\par\smallskip\noindent
{\bf Beispiel 5} Sei $X$ ein lokalkompakter Raum. Dann erh"alt man die
Einpunkt\-kompaktifizierung $X'$ von $X$, indem man die Menge der Umgebungsfilter
durch den Filter
$ F = \{ U \subseteq X \mbox{ offen } : X - U \mbox{ ist kompakt }\} $
erg"anzt. $F$ ist in der Tat ein Filter, da die Vereinigung zweier kompakter
Mengen wieder kompakt ist und ebenso eine abgeschlossene Teilmenge einer
kompakten Menge. $F$ ist konsistent, wenn $X$ nicht schon selbst kompakt ist,
worauf wir uns beschr"anken. $ X' $ ist Hausdorffsch, da f"ur $ x, y \in X $
trennende Umgebung\-en $ U $ und $V$ sofort zu trennenden Umgebung\-en $D(U) $ und
$D(V) $ f"uhren, da ja $ D(U) \cap D(V) = \emptyset $ in Filt$^\circ X$ ist.
Zu $x \in X $ gibt es $ x \in U \subseteq K $ mit $K$ kompakt, und damit sind
$ D(U) $ und $ D(X-K) $ trennende Umgebungen von $U(x) $ und $F$ . Ist
schlie\3lich $ \bigcup_{i \in I} D(U_i) $ eine "Uberdeckung von $ X'$, so ist
$ F \in D(U_j) $ f"ur ein $ j \in I $, und das Komplement von diesem $ U_j $
ist kompakt, wird also "uberdeckt von endlich vielen der $ U_i $. 
\par
Ist $X$ ein unendlicher, diskreter Raum, so ist $X$ nat"urlich lokalkompakt,
und die kompakten Teilr"aume sind einfach die endlichen Teilr"aume. Der in der
Einpunktkompaktifizierung hinzuzunehmende Filter besteht also aus allen
Komplementen endlicher Teilmengen; dieser Filter wird Frechetfilter auf $X$
genannt.
\par\smallskip\noindent
{\bf Beispiel 6} Sei $X$ ein metrischer Raum und $ Y $ seine
Vervollst"andigung. Dann entspricht einem Punkt $ y \in Y $, also einer
"Aquivalenzklasse von Cauchyfolgen, der Filter
\begin{displaymath}
\psi(y) = F(\phi^{-1}(V) : y \in V ) = \bigcap_{x_n \in y} U(x_n) .
\end{displaymath} 
Dabei bezeichnet $ \phi $ die kanonische Einbettung von $X$ in $Y$, $ x_n $
durchl"auft s"amtliche Cauchyfolgen aus der "Aquivalenzklasse $y$, und
$ U(x_n) $ steht f"ur die offenen Mengen aus $X$, die fast alle Folgenglieder
aus $x_n $ enthalten. Die Gleichheit ergibt sich so:
Ist $ B(y, \varepsilon ) $ eine offene Umgebung von $y$,
die durch $ y_n $ repr"asentiert sei, so enth"alt
$ \phi^{-1} ( B(y , \varepsilon )) $ alle Punkte mit
$ \lim_{n\to\infty} d(x,y_n) < \varepsilon /2 $ und damit auch fast alle Glieder
einer beliebigen Folge $ z_n \in y $. Ist umgekehrt $ U \not\in \psi(y) $, also
$ U \not\supseteq \phi^{-1}(V) $ f"ur alle $ V \in U(y) $, so l"a\3t sich zu
$ n \in {\bf N} $ eine Folge $ z_n $ definieren mit
$ z_n \in \phi^{-1}(B(y, 1/n)) $ und $ z_n \not\in U $. Dies ist eine Cauchyfolge
aus $y$ mit $ U \not\in U(z_n) $.
\par\smallskip\noindent
{\bf Beispiel 7} Eine Kompaktifizierung eines topologischen Raumes $X$ ist
eine dichte Einbettung $ i : X \longrightarrow K $ in einen kompakten Raum $K$.
Da ein kompakter Raum regul"ar ist, ist jede Kompaktifizierung eine
Filtrumserwei\-terung in Filt$^\circ X $. Im Filtrum lassen sich nun
verschiedene
Kompaktifizier\-ung\-en einfach vergleichen. F"ur zwei Kompakti\-fi\-zie\-rung\-en
$i : X \longrightarrow K $ und $j : X \longrightarrow K' $ gilt n"amlich 
die "Aquivalenz
\par\smallskip\noindent
(1) Es gibt eine stetige Abbildung $ f : K \longrightarrow K'$, mit
$ j = f \circ i $.
\par\smallskip\noindent
(2) F"ur jeden Filter $ F \in K $ gibt es einen Filter $ F' \in K' $ mit
$ F \supseteq F' $.
\par\smallskip\noindent
Aus (1) folgt (2). Sei $ y \in K $ und $ f(y) = z $, und wir zeigen, da\3 f"ur
die zugeh"origen Filter die Beziehung $ F_y \supseteq F_z $ gilt. Sei hierzu
$ U \in F_z $ vorgegeben, wobei man sofort $ U = j^{-1}(W) $ mit $ z \in W $
annehmen darf. Wegen der Kommutativit"at ist $ U = i^{-1}( f^{-1}(W)) $ und es
ist $ y \in f^{-1}(W) $, also $ U \in F_y $. Aus (2) folgt (1). Da
$ j : X \longrightarrow K' $ eine stetige Abbildung in einen regul"aren Raum
ist, ist $ j $ auf $K$ ausdehnbar, wenn f"ur jeden Punkt $ y \in K $ der Filter
$ j ( \psi(y)) = j (F_y) $ konvergiert. Sei $ y \in K $ und $ z \in K' $ ein
Punkt mit $ F_y \supseteq F_z $. Dann ist f"ur $ z \in W \subseteq K' $ offen
$ j^{-1}(W) \in F_z $, also $ j^{-1}(W) \in F_y $, und damit ist
$ W \in j(F_y) $, also $ j(F_y) \supseteq U(z) $, was die Konvergenz zeigt.
\par\bigskip
Zum Schlu\3 dieses Abschnitts gehen wir auf das Problem ein, inwiefern sich
ein topologischer Raum $X$ aus seiner Topologie rekonstruieren l"a\3t. Hierbei
ist unter der Topologie die Menge der offenen Mengen von $X$ samt den
Verkn"upfungen endlicher Durchschnitt und beliebige Vereinigung zu ver\-ste\-hen,
und die Frage lautet, ob man aus diesen Informationen den Raum zur"uck\-ge\-win\-nen
kann. Ferner definiert eine stetige Abbildungen von $X$ nach $Y$ ein Abbildung
von Top(Y) nach Top(X), die mit den Verkn"upfungen vertr"aglich ist, und die
Frage ist hier, ob jede solche Abbildung von einer stetigen Abbildung auf
den topologischen R"aumen herr"uhrt.
\par
Sei also $X$ ein topologischer Raum und $M$:=(Top$(X), \cap , \bigcup ) $ die
Topologie. $M$ ist ein kommutatives Monoid, man kann daher von Filtern
spre\-chen und, da wir die Vereinigung als abstrakte Verkn"upfung zur Verf"ugung
haben, auch von irreduziblen Filtern. Die Menge der irreduziblen Filter sei
mit $X'$ bezeichnet und trage die induzierte Filtrumstopologie. $X' $ ist nun
ein topologischer Raum, dessen Topologie mit der vorgegebenen "ubereinstimmt.
Hierzu ist zu zeigen, da\3 die Zuordnung $ U \longmapsto D(U) $ eine
Bijektion ist, die mit den Verkn"upfungen vertr"aglich ist. Wie immer ist
diese Zuordnung ein Monoidmorphismus. Eine offene Menge in $ X' $ ist gegeben
durch $ \bigcup_{i \in I} D(U_i) $. Da in $ X' $ aber nur irreduzible Filter
vorkommen, ist diese Menge gleich $ D(\bigcup_{i \in I} U_i) $. Damit ist die
Zuordnung mit beliebigen Vereinigungen vertr"ag\-lich, alle offenen Mengen sind
Basismengen und die Zuordnung ist surjektiv. Da die vorgegebene Topologie von
einem Raum herr"uhrt, folgt aus $ V \neq U $, da\3 es einen Punkt $ x \in X $
gibt mit $ x \in V$, $ x \not\in U $, und
daher gibt es auch einen irreduziblen Filter $F$ mit
$ F \in D(V)$, $F \not\in D(U) $, was die Injektivit"at zeigt. Damit haben
wir zu $M$ einen topologischen Raum $X'$ gefunden mit $M$ = Top$X'$.
\par
Ist jetzt $ \psi : {\rm Top}(Y) \longrightarrow {\rm Top}(X) $ ein
Monoidmorphismus,
der auch noch mit Vereinigungen vertr"aglich ist, so l"a\3t sich die
Abbildung $ \psi^\ast : \mbox{ \rm Filt}X \longrightarrow \mbox{ \rm Filt}Y $ auf
die irreduziblen Filter einschr"anken, und man gelangt so zu einer stetigen
Abbildung $ \psi^\ast : X' \longrightarrow Y' $. Wegen
$ (\psi^\ast)^{-1}(D(U)) = D(\psi(U)) $ gewinnt man aus dieser stetigen
Abbildung den vorgegebenen Monoidmor\-phis\-mus zur"uck. Dieselbe Problematik
wird auch in \cite{art}, expos\'{e} IV, Abschnit\-te 2.1, 4.2, 7.1,
allerdings in der Sprache der Topoi, behandelt.
\par
Der Ausgangsraum $X$ steht mittels der Abbildung
$ i :X \longrightarrow X' \subseteq \mbox{ Filt(Top(X))} $ mit
$ X' $ in kanonischer Beziehung. Diese Abbildung ist, wie nebenbei gezeigt, ein
Quasi-Hom"oomorphismus und $X'$ ist die Sobrierung von $X$ im Sinne
von \cite{gro}, p.64 und p.67-70. $ X' $ ist ein Raum, in dem
jeder irreduzible Filter Umgebungsfilter genau eines Punktes ist, und $X$
stimmt mit $X'$ genau dann "uberein, wenn bereits $X$ diese Eigen\-schaft hat.
Schon weiter oben wurde gezeigt, da\3 diese Eigenschaft gleich\-be\-deu\-tend ist zu
der Eigenschaft, da\3 jede irreduzible abgeschlossene Teilmenge genau einen
generischen Punkt hat. Dies ist beispielsweise dann erf"ullt, wenn $X$ ein
Hausdorffraum oder der zugrunde liegende Raum eines noetherschen
Schemas ist (allerdings nicht, wenn $X$ nur ein $T_1$-Raum ist, wie
f"alschlich in \cite{art}, IV.4.2.1 behauptet).
\par
Die Sobrierung $X'$ eines topologischen Raumes $X$ spielt in der algebra\-ischen
Geometrie eine Rolle, da damit der funktorielle "Ubergang
von quasi\-pro\-jek\-tiven
Variet"aten zu Schemata auf topologischer Ebene beschrieben wird. Ein Schema 
er\-h"alt man dann, wenn man noch die Struktur\-garbe "uber\-tr"agt,
(\cite{har}, II.2.6).
Wie mir scheint ist die Be\-schrei\-bung von $X'$ mit ir\-re\-du\-zi\-blen Fil\-tzern
ein\-facher und "uber\-sicht\-lich\-er als mit
ir\-re\-du\-zi\-blen, ab\-ge\-schlos\-senen Mengen. 
\section{Das Filtrum als garbentheoretisches Ob\-jekt}
In diesem Abschnitt geht es um Garben, und zwar immer um Garben abel\-scher
Gruppen. In der Kategorie der abelschen Gruppen existieren induktive Limiten
"uber nach oben gerichtete Indexmengen. Ist $ \cal G $ eine Garbe auf dem
topologischen Raum $X$, so bilden zu einem Punkt $ x \in X $ die Gruppen
${ \cal G}(U) ,x \in U $ mit den Restriktionsabbildungen
$ R^U_V : {\cal G}(U) \longrightarrow {\cal G}(V) $ f"ur $ V \subseteq U $ ein
induktives System, und der induktive Limes dieses Systems hei\3t der Halm
der Garbe $\cal G$ im Punkt $x$.
\par
Offenbar ist diese Konstruktion nicht nur
f"ur Umgebungsfilter m"oglich, sondern generell f"ur beliebige topologische
Filter. Zun"achst soll gezeigt wer\-den, da\3 diese verallgemeinerten Halme in
nat"urlicher Weise auftreten, selbst dann, wenn das eigentliche Interesse nur den
Halmen in einem Punkt gilt. Dazu erinnern wir an das Bild- und Urbildnehmen von
Garben.
\par\bigskip
Sei $ \phi: X \longrightarrow Y $ eine stetige Abbildung. Ist $\cal G$ eine
Garbe auf $Y$, so definiert man auf $X$ die Urbildgarbe $\phi^{-1}{\cal G} $ als
die Garbe zur Pr"agarbe
\begin{displaymath}
U \longmapsto \lim_{\phi(U) \subseteq V} {\cal G}(V)
= \lim_{U \subseteq \phi^{-1}(V)} {\cal G}(V) .
\end{displaymath}
Der Halm von $\phi^{-1} {\cal G} $ in einem Punkt $ x \in X $ berechnet sich dann
zu
\begin{displaymath}
\lim_{x \in U} \phi^{-1} {\cal G}(U) =
\lim_{x \in U} \lim_{U \subseteq \phi^{-1}(V)} {\cal G}(V)
= \lim_{x \in \phi^{-1}(V)} {\cal G}(V) = \lim_{\phi(x) \in V} {\cal G}(V)
= {\cal G}_{\phi(x)}.
\end{displaymath}
Ist umgekehrt $\cal G $ eine Garbe auf $X$, so wird die Bildgarbe
$ \phi_\ast {\cal G} $ definiert durch die Zuordnung
\begin{displaymath}
V \longmapsto {\cal G}(\phi^{-1}(V))  ,
\end{displaymath}
die bereits eine Garbe ist. Der Halm in einem Punkt $ y \in Y $ von
$ \phi_\ast {\cal G} $ berechnet sich zu
\begin{displaymath}
(\phi_\ast {\cal G})_y = \lim_{y \in V} {\cal G}(\phi^{-1}(V)) =
\lim_{U \in \psi(y)} {\cal G}(U) ,
\end{displaymath}
wobei wieder $ \psi $ die schon oft benutzte Abbildung
$ Y \longrightarrow \mbox{ Filt}X $ ist. $ \psi(y) $ ist aber im allgemeinen
kein Umgebungsfilter eines Punktes. F"ur den verallge\-mein\-erten Halm einer
Garbe $\cal G$ in einem Filter $F$ schreiben wir $ {\cal G}_F $ oder
$ {\cal G}(F) $.
\par\bigskip\noindent
{\bf Beispiel 1} Sei $ i: X \longrightarrow \mbox{ Filt}X $ die kanonische
Abbildung und $\cal G$ eine Garbe auf $X$. Dann ist der Halm von
$ i_\ast {\cal G} $ in einem Punkt $ F \in \mbox{ Filt}X $ gleich dem
verallgemeinerten Halm von $\cal G$ im Filter $F$. Es ist ja $ \psi(F) = F $, und
so k"onnen "uberhaupt alle verallgemeinerten Halme als Punkthalme auftreten.
\par\smallskip\noindent
{\bf Beispiel 2} Ist $U$ eine offene Teilmenge von $X$, so kann man
$ {\cal G}(U) $ formal als Halm im Filter $F(U)$ auffassen.
\par\smallskip\noindent
{\bf Beispiel 3} Sei $ X =$ Spek$A, \tilde{A} $ das Spektrum eines kommutativen
Ringes $A$, und $ F \subseteq A $ ein ringtheoretischer Filter. Ihm entspricht
der topologische Filter $ F':= F ( \{ D(f) : f \in F \} ) $, und es gilt
$ \tilde{A}_{F'} = A_F $. Im Spektrum finden sich also s"amtliche Bruchringe
von $A$ als gewisse verallgemeinerte Halme wieder. Dasselbe gilt nat"urlich
auch f"ur einen beliebigen quasikoh"arenten $\tilde{A}$-Modul auf $X$. Die
Aussage $\tilde{M}_x = M_{A-{\bf p}_x} $ ist ein Spezialfall davon.
\par\smallskip\noindent
{\bf Beispiel 4} Oben haben wir gezeigt, da\3 f"ur eine stetige Abbildung
$ \phi : X \longrightarrow Y $ und f"ur den Halm einer Garbe $\cal G$ in
einem
Punkt $y \in Y $ $ (\phi_\ast {\cal G})_y = {\cal G}_{\psi(y)} $ gilt. Die obige Argumentation
l"a\3t sich f"ur beliebige Filter $ F$ aus $Y$ aufrechterhalten und man gelangt zu
$ (\phi_\ast {\cal G})_F = {\cal G}_{\phi^{-1}(F)} $.
F"ur $ F = F(V) $ erh"alt
man nat"urlich wegen $ \phi^{-1} F(V) = F(\phi^{-1}(V)) $ die
Ausgangs\-defi\-ni\-tion
$ (\phi_\ast {\cal G})(V) = {\cal G}(\phi^{-1}(V)) $ zur"uck.
\par\smallskip\noindent
{\bf Beispiel 5} Sei $A$ ein Integrit"atsbereich, und
$ X= \mbox{ Spek}A, \tilde{A} $. F"ur die Strukturgarbe in einer offenen Menge
$ U \subseteq X $ gilt hier $ \Gamma(U) = \bigcap_{ {\bf p} \in U} A_{\bf p} $,
wobei die Lokalisierungen $ A_{\bf p} $ als Unterringe des Quotientenk"orpers
$ Q = Q(A) $ aufzufassen sind. F"ur den Halm von $\tilde{A} $ in einem
Umgebungsfilter einer Teilmenge $ T \subseteq X $ gilt nun entsprechend
\begin{displaymath}
 \Gamma (U(T), \tilde{A} ) = \lim_{T \subseteq U} \Gamma(U) =
 \bigcap _{ {\bf p} \in T} A_{\bf p} .
\end{displaymath}
F"ur $ T \subseteq U $ gilt offenbar
$ \bigcap_{{\bf p} \in U } A_{\bf p} \subseteq
\bigcap_{{\bf p} \in T} A_{\bf p} $ und damit
$ \Gamma(U(T) , \tilde{A} ) \subseteq \bigcap_{{\bf p} \in T} A_{\bf p} $. Ist
umgekehrt $ q \in \bigcap_{{\bf p} \in T} A_{\bf p} $, so gibt es f"ur jedes
${\bf p} \in T $ eine Darstellung von $q$ mit einem Nenner, der nicht zu
${\bf p}$ geh"ort, und $q$ ist auf der Vereinigung der Basismengen definiert,
die zu diesen Nennern geh"oren, und damit ist
$ q \in \Gamma( U(T) , \tilde{A} ) $. Ist $A$ zus"atzlich noethersch, so ist
$X$ ein Zariskiraum, und jeder verallgemeinerte Halm ist von diesem Typ, da
dann jeder topologische Filter Umgebungsfilter einer Teilmenge ist.
\begin{satz} Sei $ \phi: X \longrightarrow Y $ eine stetige Abbildung, und
$ \cal G $ eine Garbe auf $X$. Ferner seien wieder
$ i: X \longrightarrow \mbox{ \rm Filt}X $ und
$ \psi: Y \longrightarrow \mbox{ \rm Filt}X $
die kanonischen Abbildungen. Dann ist
$ \phi_\ast {\cal G} = \psi^{-1} (i_\ast {\cal G}) $, d.h. die einzige
we\-sentliche Bildgarbe zu einer Garbe ist die durch i auf das Filtrum induzierte.
\end{satz}
Beweis. Schon in der Vorbemerkung zu diesem Abschnitt wurde die behaup\-te\-te
Gleichheit lokal, also in den Halmen der Punkte $ y \in Y $ gezeigt, und daraus
ergibt sich die globale Aussage.
\par\bigskip
Ist in obiger Situation $ \phi $ eine Einbettung, so gilt f"ur
$ x \in X $ und $ y= \phi(x) $ nach Satz 3.2.3 $ \psi(y) = U(x) $, und damit
$ \phi_\ast {\cal G}_y = {\cal G}_x $.
\par\bigskip
Zu einer Garbe $\cal G $ auf einem topologischen Raum $X$ und einer stetigen
Abbildung $ \phi : X \longrightarrow Y $ lassen sich neben der Bildgarbe
$ \phi_\ast {\cal G} $ auch noch die sogenannten h"oheren Bildgarben
definieren (\cite{har}, III.8). Hierzu geht man wie bei der Definition
der Kohomologiegruppen von
einer injektiven Aufl"osung
\begin{displaymath}
0 \longrightarrow {\cal G} \longrightarrow I_0 \longrightarrow I_1
\longrightarrow \mbox{ etc. }
\end{displaymath}
aus, und betrachtet davon die Bildgarben auf $Y$ unter Auslassung von
$ \phi_\ast {\cal G} $
\begin{displaymath}
0 \longrightarrow \phi_\ast I_o \longrightarrow \phi_\ast I_1 \longrightarrow 
\mbox{ etc. }
\end{displaymath}
Diese Sequenz ist zwar nicht mehr exakt, aber ein Komplex, und so definiert man
die h"oheren Bildgarben $ R^i \phi_\ast ({\cal G}) $ als die i-te Homologiegarbe
des Kom\-plex\-es. F"ur $ i =0 $ ergibt sich einfach $ \phi_\ast {\cal G} $. Man
kann zeigen, da\3 die h"oheren Bildgarben mit den Garben zur Pr"agarbe auf $Y$
\begin{displaymath}
V \longmapsto H^i( \phi^{-1}(V) , {\cal G} \mid _{ \phi^{-1}(V)})
\end{displaymath}
"ubereinstimmen. Daraus ergibt sich nun sofort f"ur den
Halm von $ R^i \phi_\ast ({\cal G})$ in einem Punkt $ y \in Y $ die
Beziehung 
\begin{displaymath}
( R^i \phi_\ast ({\cal G} ) )_y =
\lim_{ y \in V} H^i (\phi^{-1}(V) , {\cal G} \mid _{\phi^{-1}(V)} )
= \lim_{ U \in \psi(y)} H^i (U, {\cal G} \mid_U ),
\end{displaymath}
wobei wieder $ \psi(y) = \phi^{-1}(U(y)) $ gesetzt wurde. Da diese Beziehung
speziell f"ur die h"oheren Bildgarben auf dem Filrum von $X$ bez"uglich der
kanonischen Abbildung $ i : X \longrightarrow \mbox{ \rm Filt} X $ gilt, ergibt sich
der Satz:
\begin{satz}
Ist $ \phi : X \longrightarrow Y $ eine stetige Abbildung und $ \cal G $ eine
Garbe abel\-scher Gruppen auf $X$, so sind die h"oheren Bildgarben 
$ R^i \phi _\ast ({\cal G}) $ auf $Y$ gleich den Urbildgarben unter $ \psi $
der h"oheren Bildgarben $ R^i i_\ast ({\cal G}) $ auf dem Filtrum von $X$. Ist
$ R^i i_\ast ({\cal G}) = 0 $, so gilt dies auch f"ur die i-te h"ohere Bildgarbe
auf $Y$, und in diesem Sinn besitzt die kanonische Abbildung in das Filtrum die
maximale relative Kohomologie.
\end{satz}
\par
Gew"ohnlich besteht zwischen den oben beschriebenen Pr"agarben und den zugeh"origen
h"oheren Bildgarben eine gro\3er Unterschied. Der Extremfall zeigt sich, wenn
man auf $X$ die Identit"at betrachtet. Dann verschwinden die h"oheren
Bildgarben f"ur $ i \ge 1 $, wie sich unmittelbar aus der Definition ergibt,
w"ahrend nat"urlich die Pr"agarben $ U \longmapsto H^i (U, {\cal G} \mid_U ) $
in der Regel nicht verschwinden. Ein derartiger Kollaps einer Pr"agarbe kann nun
auf einem Filtrum nicht passieren, da dort jede Basismenge $D(U)$ den globalen
Punkt $F(U)$ besitzt und somit auf diesen Basismengen die Garbe zu einer
Pr"agarbe mit dieser "ubereinstimmt. Damit ergibt sich insbesondere 
$ R^i i_\ast ({\cal G}) (D(U)) = H^i (U , {\cal G} \mid _U ) $, und die
kanonische Abbildung $ i : X \longrightarrow \mbox{ Filt}X $ hat nur dann keine
relative Kohomologie, wenn f"ur s"amtliche offene Teilmengen $ U \subseteq X $
die Kohomologie verschwindet. Die i-te Bildgarbe zu $\cal G $ auf dem Filtrum
nennen wir im
folgenden einfach die i-te Kohomologiegarbe von $\cal G$.
\par
Im allgemeinen ist es nat"urlich schwierig, Aussagen "uber die Kohomolo\-gie\-garben
zu machen. Neben direkten "Uberlegungen ist es h"aufig sinnvoll, Ergebnisse
"uber die relative Kohomologie auf $Y$ bez"uglich einer stetigen Abbildung
$ \phi : X \longrightarrow Y $ mittels $\psi $ auf das Filtrum zu "ubertragen
um so wenigstens "uber die Halme der Kohomologiegarben in den Urbildfiltern
unter $\psi$ etwas aussagen zu k"onnen. So ergibt sich beispielsweise sofort,
da\3 die Kohomologiegarben in den Umgebungsfiltern von Punkten verschwinden.
Allgemeiner gilt, da\3 die Kohomologie in jedem irreduziblen Filter
ver\-schwin\-det. Dies ergibt sich daraus, da\3 die Abbildung
$ i: X \longrightarrow X' $ keine relative Kohomologie hat, wozu man zu zeigen
hat, da\3 unter dieser Abbil\-dung Bild- und Urbildnehmen von Garben invers
zueinander sind, da sich die offenen Teilmengen bijektiv entsprechen
(\cite{gro}, p.85f). Man kann auch mit folgendem Lemma zum Ziel kommen:
\begin{lem}
Ist $ {\cal F} \longrightarrow {\cal G} \longrightarrow {\cal H} $ eine exakte
Garbensequenz auf dem topologischen Raum $X$ und $F$ ein irreduzibler Filter,
so ist auch die Sequenz
$ {\cal F}_F \longrightarrow {\cal G}_F \longrightarrow {\cal H}_F $
exakt.
\end{lem}
Beweis. Die Abbildungen seien mit $ \alpha $ und $ \beta $ bezeichnet. 
Wir haben zu zeigen, da\3 ein Element $ s \in {\cal G}_F $
mit $ \beta(s) = 0 $ in $ {\cal H}_F $ von einem Element aus ${\cal F}_F $
herr"uhrt. Sei $ s \in {\cal G}(U) , U \in F $. $ \beta(s) $ verschwindet in
$F$ und damit in einer offenen Menge, die zu $F$ geh"ort, und wir k"onnen 
$ \beta(s) = 0 $ in $ {\cal H} (U) $ annehmen. Wegen der Exaktheit in den
Punkthalmen gibt es f"ur $ x \in U $ Elemente $ t _x \in {\cal F}_x $ mit
$ \alpha(t_x) = s_x $, und dies gilt dann wieder in offenen Umgebungen $ U_x $
von $x$, also $ \alpha(t_x) = s $ in $ {\cal G}( U_x) $ und
$ t_x \in {\cal F}(U_x) $. Da $U \in F $ von den $U_x $ "uberdeckt werden
und $F$ irreduzibel ist, ist $ U_x \in F $ f"ur ein $x$, und damit ist
$ t_x \in {\cal F}_F $ und $ \alpha(t_x) = s $ in $F$.
\par\bigskip
\begin{satz}
In einem irreduziblen Filter verschwinden die Kohomologie\-gar\-ben.
\end{satz}
Beweis. Der zu einer injektiven Aufl"osung einer Garbe geh"orende Garben\-komplex auf
dem Filtrum ist nach dem Lemma in den Halmen zu irreduziblen Filtern exakt,
und somit verschwinden an diesen Stellen die Quotientengar\-ben, die ja nach
Definition die Kohomologiegarben sind.

%% file: rang.tex
\chapter{Filter in beringten R"aumen}
\section{Die Invertierbarkeitsmenge einer Funk\-tion}
In diesem Kapitel geht es um die Frage, inwiefern bei einem beringten Raum
$X, {\cal O}_X $ die topologischen Filter auf $X$ mit den Filtern im globalen
Schnittring $ \Gamma (X, {\cal O}_X) = A $ in Beziehung stehen. Generell ist
es eine wichtige Frage, in dieser Situation topologische und algebraische
Objekte zu vergleichen, wobei auf der einen Seite beispielsweise offene Mengen,
Punkte, etc. und auf der anderen Ideale, Primideale, Funktionen etc. zur
Diskussion stehen. Da es gute Gr"unde daf"ur gibt, gewisse topologische und
ringtheoretische Objekte mit dem einen Begriff des Filters zu erfassen, liegt die
Vermutung nahe, da\3 sich diese Begriffskorrespondenz auch in der Sache
best"atigt. Zun"achst mu\3 aber eine Korrespondenz auf der Ebene der Monoide
$A$ und Top$X$ gestiftet werden.
\begin{pro}
Sei $ X, {\cal O}_X $ ein beringter Raum mit $ \Gamma (X, {\cal O}_X ) = A $.
Dann gilt:
\par\smallskip\noindent
(1) F"ur $ f \in A $ ist $ X_f := \{ x \in X : f_x \in {\cal O}_x^\times \}$ eine
offene Menge in $X$. Sie hei\3t die Invertierbarkeitsmenge von $f$.
\par\smallskip\noindent
(2) Die Zuordnung $ \chi : A \longrightarrow \mbox{ \rm Top}(X)$ mit
$ f \longmapsto X_f $ ist ein Monoidmor\-phis\-mus.
Ist $F(f) \subseteq F(g) $, so ist $ X_g \subseteq X_f $.
\par\smallskip\noindent
(3) $f$ ist eine Einheit in $\Gamma(X_f,{\cal O}_X)$. Ist
$ f \in \Gamma(U, {\cal O}_X)^\times $, so ist $ U \subseteq X_f $.
Insbesondere ist
$X_f = X $ genau dann, wenn $f$ eine Einheit in $A$ ist.
\par\smallskip\noindent
(4) Ist $ X, {\cal O}_X $ zus"atzlich lokal beringt, so gilt f"ur $f,g \in A $
die Beziehung $ X_{f+g} \subseteq X_f \cup X_g $. Ist $ f \in A $ nilpotent, so
ist $X_f = \emptyset $.
\end{pro}
Beweis. (1) Sei $ x \in X_f $, also $ f_x \in {\cal O}_x^\times $. Damit gibt
es ein Inverses $ f' \in {\cal O}_x $ mit $ f_x \cdot f' =1 $. Dieses $f'$
kommt nun von einer offenen Umgebung von $x$ her, und die Gleichheit gilt ebenfalls in
einer offenen Umgebung, die man als gleich
annehmen darf. Sei also $ f \cdot f' = 1 $ in $\Gamma (U, {\cal O}_X \mid _U) $
mit $ x \in U $. F"ur $ y \in U $ gilt dann im Halm $ f_y \cdot f'_y =1 $, und
damit ist $ f_y $ eine Einheit in ${\cal O}_y$ und damit
$ x \in U \subseteq X_f $, also ist $X_f$ offen.
\par\noindent
(2) Einheiten aus $A$ bleiben unter den Ringmorphismen
$ A \longrightarrow {\cal O}_x $ Einheiten und sind somit in jedem Halm
invertierbar, insbesondere ist also $ X_1 = X $. Zur Multiplikativit"at ist
$X_{f \cdot g} = X_f \cap X_g $ f"ur beliebige $f,g \in A $ zu zeigen. Ist
$ x \in X _{f \cdot g} $, so ist $ (f \cdot g)_x = f_x \cdot g_x $ eine
Einheit in ${\cal O}_x $ und damit auch $f_x $ und $g_x $, und $x$ geh"ort
zum Durchschnitt $ x \in X_f \cap X_g $. Zum Beweis der anderen Inklusion geht
man diesen Weg r"uckw"arts. Ist $ f \in F(g) $, so kann man $ f \cdot a = g^n $
mit $ n \ge 1 $ ansetzen und erh"alt, da ein Element genau dann Einheit ist,
wenn dies f"ur
eine echte Potenz gilt, $X_g = X_{g^n} = X_{f \cdot a}
= X_f \cap X_a \subseteq X_f $.
\par\noindent
(3) Wie in der Formulierung im Satz bezeichnen wir die Einschr"ankung von
$ f\in A$ auf $\Gamma(X_f,{\cal O}_X) $ wieder mit $f$. Da $f$ in jedem Halm
von Punkten $ x \in X_f $ eine Einheit ist,
gibt es eine "Uberdeckung
$X_f = \bigcup_{i \in I} U_i $ und Elemente $f_i \in {\cal O}(U_i)$ mit
$f \cdot f_i =1 $, wobei die Beschr"ankung von $f$ auf $U_i$ einfach wieder mit
$f$ bezeichnet wird. Damit liefern auf $U_i \cap U_j $ sowohl $f_i$ als auch
$f_j $ Inverse zu $f$ und m"ussen daher gleich sein. Damit ist die
Vertr"aglichkeitsbedingung erf"ullt und es gibt ein
$f'\in \Gamma(X_f,{\cal O}_X) $, dessen
Beschr"ankungen die $f_i $ sind. Damit ist $ f \cdot f' - 1 $ gleich null auf
den $U_i $, und damit gilt dies auf $X_f$, $f'$ ist also Inverses zu $f$.
Ist $f$ auf $U$ invertierbar, so nat"urlich auch in jedem Halm von Punkten aus
$U$, und damit ist $ U \subseteq X_f $. Die letzte Aussage ist klar.
\par\noindent
(4) Ist $ x \in X_{f+g} $, so ist $ f+g $ eine Einheit in dem lokalen Ring
$ {\cal O}_x $. Damit mu\3 bereits $f$ oder $g$ eine Einheit sein, also
$ x \in X_f \cup X_g $. Ein nilpotentes Element $f \in A $ bleibt unter den
Restriktionsmorphismen nilpotent, und kann somit nur im Nullring eine Einheit
sein, ein lokaler Ring ist aber nach Definition vom Nullring verschieden.
\par\bigskip\noindent
{\bf Bemerkung} In (2) wurde gezeigt, da\3 $X_f$ die maximale Menge von $X$
ist, auf der $f$ invertierbar ist. Damit ist wohl die Bezeichnung
Invertierbarkeitsmenge der Funktion $ f\in A $ angebracht.
\par\smallskip\noindent
{\bf Beispiel 1} Ist $ (X, {\cal O}_X) = ({\rm Spek} A, \tilde{A}) $ , so ist
$X_f = D(f) $. Hier lassen sich einige Aussagen von oben umkehren.
Beispielsweise folgt aus $ D(g) \subseteq D(f) $ stets $ F(f) \subseteq F(g) $,
und aus $ D(f) = \emptyset $ folgt, da\3 $f$ nilpotent ist. Man beachte, da\3
im Fall eines projektiven Schemas die Basismengen $D_+(f)$ f"ur ein homogenes
Element $f$ aus dem homogenen Koordinatenring nichts mit den oben definierten
offenen Mengen $X_f$ zu tun hat; $X_f $ ist ja nur f"ur Elemente aus dem
globalen Schnittring definiert, der vom Koordinatenring verschieden ist.
\par\smallskip\noindent
{\bf Beispiel 2} Ist $ (X, {\cal O}_X) = (X, {\cal C}_{\bf R}(X) ) $ ein
topologischer Raum mit der Garbe der reellwertigen, stetigen Funktionen,
so ist $X_f = f^{-1}({\bf R}^\times) $ das Komplement der
Nullstellenmenge von $f$. Ist $ x \in X_f $, so mu\3 sicherlich $f(x) \neq 0 $
sein; ist umgekehrt $f(x) \neq 0 $, so ist $f $ bereits in einer Umgebung von
$x$ von null verschieden und damit invertierbar. Dasselbe gilt f"ur die
komplexen Zahlen und ebenso f"ur Mannigfaltigkeiten mit den Garben der
differenzierbaren oder analytischen Funktionen, da deren Inverses auf einer
nullstellenfreien offenen Menge wieder differenzierbar beziehungsweise
analytisch ist. 
\par\bigskip
Die offenen Mengen $X_f $ h"angen aufs engste mit den Basismengen $D(f) $ im
Spektrum beziehungsweise Filtrum eines kommutativen Ringes zusammen. Um sich
dies klar zu machen, sei zun"achst an folgende universelle Eigenschaft des
Spektrums Spek$B, \tilde{B} $ eines kommutativen Ringes $B$ in der
Kategorie der lokal beringten R"aume erinnert, (\cite{gro} p.210):
Ist $ X, {\cal O}_X $ ein lokal
beringter Raum mit dem globalen Schnittring $ {\cal O}_X(X) = A $, und ist
$ \varphi : B \longrightarrow A $ ein Ringmorphismus, so gibt es genau einen
Morphismus lokal beringter R"aume
\begin{displaymath}
\phi : X, {\cal O}_X  \longrightarrow \mbox{ Spek}B, \tilde{B} ,
\end{displaymath}
dessen globaler Ringmorphismus gerade der vorgegebene ist. Dabei ist ein
Morphismus lokal beringter R"aume definitionsgem"a\3 ein Morphismus bering\-ter
R"aume, dessen s"amtliche Halmabbilungen lokal sind; ohne diese Voraus\-setz\-ung
gilt die Eindeutigkeit der Abbildung nicht. Auf der topologischen Ebene wird
dabei einem Punkt $ x \in X $ das Primideal
\begin{displaymath}
\{g \in B : (\varphi(g))_x \in m_x \} = (R_x \circ \varphi)^{-1} (m_x) ,
\end{displaymath}
zugeordnet, wobei $ R_x $ die Restriktionsabbildung des globalen Schnittrings in den lokalen
Ring $ {\cal O}_x $ und $m_x $ das maximale Ideal davon bezeichnet.
\par
Im folgenden werden wir zeigen, da\3 das Filtrum eines kommutativen Ringes
dieselbe universelle Eigenschaft in der Kategorie aller beringten R"au\-me hat wie
das Spektrum in der Kategorie der lokal beringten R"aume. Hierzu mu\3 das
Filtrum zu einem beringten Raum gemacht und "uber die erlaubten Morphismen
beringter R"aume eine Konvention getroffen werden. Zun"achst l"a\3t sich
mittels der kanonischen Abbildung
$ i: \mbox{ Spek}B \longrightarrow \mbox{ Filt}B $ die Strukturgarbe $\tilde{B}$
auf das Filtrum "ubertragen, die wir wieder mit $\tilde{B} $ be\-zeich\-nen.
F"ur Basismengen $D(g) $ im Filtrum gilt dann
$ \Gamma (D(g), \tilde{B})=B_g $, und f"ur den Halm im Filter $F$ gilt
nat"urlich $ \tilde{B}_F = B_F $.
\par
Ferner legen wir fest, unter einem Morphismus beringter R"aume
$ \phi, \varphi : X,{\cal O}_X \longrightarrow Y,{\cal O}_Y $ nur
solche zuzulassen, die in den
Halmen lokal sind. Es wird also gefordert, da\3 f"ur $ y = \phi(x) $ der
Ringmorphismus $ \varphi_{y,x} : {\cal O}_y \longrightarrow {\cal O}_x $
Nichteinheiten auf Nichteinheiten abbildet. Sind die R"aume lokal beringt, so
f"allt diese Forderung mit der "ublichen Definition lokaler Morphismen
zusammen, da ja in einem lokalen Ring die Nichteinheiten gerade die Elemen\-te
des einzigen maximalen Ideals sind. Das folgende Lemma zeigt bereits eine
w"unschenswerte Konsequenz dieser Festlegung:
\begin{lem}
Ist $ \phi, \varphi : X,{\cal O}_X \longrightarrow Y, {\cal O}_Y $ ein
Morphismus beringter R"aume, mit dem globalen Ringmorphismus
$ \varphi : B=\Gamma(Y, {\cal O}_Y) \longrightarrow A=\Gamma(X, {\cal O}_X) $,
so kommutiert das folgende Diagramm.
\par\noindent
\unitlength0.5cm
\begin{picture}(10,7)
\put(4,1){\vector(1,0){3}}
\put(4,5){\vector(1,0){3}}
\put(2.8,4){\vector(0,-1){2}}
\put(8.2,4){\vector(0,-1){2}}
\put(2.6,4.7){$B$}
\put(8,4.7){$A$}
\put(0.7,0.8){{\rm Top}(Y)}
\put(7.5,0.8){{\rm Top}(X)}
\put(5.5,5.5){$\varphi$}
\put(5.5,1.5){$\phi^{-1}$}
\end{picture} .
\par\noindent
Es gilt also $ \phi^{-1}(Y_g) = X_{\varphi(g)} $ f"ur $ g \in B $.
\end{lem}
Beweis. Man beachte, da\3 der globale Ringmorphismus und die Halmabbil\-dung\-en
mit den Restriktionen vertr"aglich sind, da\3 also
$\varphi(g)_x = \varphi(g_{\phi(x)}) $. Sei nun $ x \in X_{\varphi(g)} $,
also $ \varphi(g)_x $ eine
Einheit in $ {\cal O}_x $. Da der Halmmorphismus
${\cal O}_{\phi(x)} \longrightarrow {\cal O}_x $ nach unserer Definition
lokal ist und $g_{\phi(x)} $
auf $\varphi(g)_x $ abbildet, ist auch
$ g_{\phi(x)}$ eine Einheit in $ {\cal O}_{\phi(x)} $, also ist
$ \phi(x) \in Y_g $ und damit $ x \in \phi^{-1}(Y_g) $.
Ist umgekehrt $ x \in \phi^{-1}(Y_g) $, so wird unter der Halmabbildung die
Einheit $ g_{\phi(x)} $ auf die Einheit $ \varphi(g)_x $ in ${\cal O}_x $
abgebildet, und somit ist $ x \in X_{\varphi(g)} $.
\begin{satz}
Ist $ X,{\cal O}_X $ ein beringter Raum mit $\Gamma(X, {\cal O}_X) = A $ und
ist $ \varphi : B \longrightarrow A $ ein Ringmorphismus kommutativer Ringe,
so gibt es genau einen Morphismus beringter R"aume
\begin{displaymath}
\phi,\varphi : X,{\cal O}_X \longrightarrow \mbox{ \rm Filt}B,\tilde{B} \hspace{1em},
\end{displaymath}
dessen globaler Ringmorphismus gerade der vorgegebene ist.
\end{satz}
Beweis. Aufgrund des Lemmas mu\3 f"ur einen Morphismus die "Aquivalenz
$ x \in X_{\varphi(g)} $ genau dann ,wenn $ \phi(x) \in Y_g = D(g) $
gelten, und damit kommt als Bildfilter von $ x \in X $ nur
$ \phi(x) = \{ g \in B : x \in X_{\varphi(g)} \}$ in Frage. Dies ist der
Urbildfilter des Umgebungsfilters von $x$ unter der zusammengesetzten Abbildung
$ B \longrightarrow A \longrightarrow \mbox{ Top}(X) $, oder, was dasselbe ist,
der Urbildfilter der Einheitengruppe von $ {\cal O}_x $ unter der Abbildung
$ B \longrightarrow A \longrightarrow {\cal O}_x $. Die Abbil\-dung ist ste\-tig,
da ja $ \phi^{-1}(D(g)) = X_{\varphi(g)} $ offen ist. Als glo\-balen
Ring\-mor\-phis\-mus legen wir $ \varphi : B \longrightarrow A $
fest und er\-hal\-ten
f"ur $ g \in B $ Ring\-mor\-phis\-men
$ B \longrightarrow \Gamma(X_{\varphi(g)},{\cal O}_X) $. Da $ \varphi(g) $ in
$ \Gamma(X_{\varphi(g)},{\cal O}_X) $ eine Ein\-heit ist, gibt es,
auf\-grund der uni\-ver\-sel\-len Ei\-gen\-schaft der Bruch\-ringe, genau einen
Ring\-mor\-phis\-mus
$\Gamma(D(g), \tilde{B})=B_g \longrightarrow \Gamma(X_{\varphi(g)},{\cal O}_X)$.
F"ur eine belie\-bige of\-fene Menge $ U \subseteq \mbox{ Filt}B $ erge\-ben sich
dann die noch feh\-len\-den Ring\-mor\-phis\-men in ein\-deu\-tiger Weise durch
Zusam\-men\-kleben der bisher konstru\-ierten Mor\-phis\-men. Die Halmab\-bil\-dungen
$ B_{\phi(x)} =\tilde{B}_{\phi(x)} \longrightarrow {\cal O}_x $ sind dabei
lokal, da ja $ \phi(x) $ der Urbild\-filter der Einhei\-ten\-gruppe von ${\cal O}_x $
ist.
\par
Insbesondere gibt es zu jedem beringten Raum $X,{\cal O}_X $
mit $\Gamma(X,{\cal O}_X) = A$ einen kanonischen Morphismus
$X,{\cal O}_X \longrightarrow \mbox{ Filt}A,\tilde{A} $, der sich ergibt, wenn
man im obigen Satz f"ur $\varphi $ die Identit"at ansetzt. Es ist dann
$ X_f = \phi^{-1}(D(f)) $, und in diesem Sinn sind die $D(f)$ die typischen
Beispiele f"ur die $X_f $. Eine weitere Deutung der $X_f$ ergibt sich, wenn
man $A$ in der eindeutigen Weise als ${\bf Z} $-Algebra auffa\3t und den durch
$f$ definierten ${\bf Z}-$Algebramorphismus $ {\bf Z}[T] \longrightarrow A $
betrachtet und dazu den Morphismus
$ X,{\cal O}_X \longrightarrow \mbox{ Filt}{\bf Z}[T], \widetilde{{\bf Z}[T]} $.
Unter diesem Morphismus ist $ X_f = \phi^{-1}(D(T)) $.
\par
Ist $X,{\cal O}_X $ ein lokal beringter Raum, so sind die Filter,
die einem Punkt zugeordnet werden, Komplemente von Primidealen, und die kanonischen
Morphismen in die Spektren ergeben sich aus den Abbildungen in die Filtren
durch Einschr"anken
der Bildbereiche. Diese Aussage ist ein Spezialfall der folgenden Behaup\-tung.
\begin{satz}
Sei $X,{\cal O}_X $ ein lokal beringter Raum mit globalem Schnittring $A$
und $F$ ein quasikompakter und konsistenter Filter auf $X$. Dann gilt
\par\smallskip\noindent
(1) Unter der zu $\chi : A \longrightarrow \mbox{ \rm Top}(X) $ geh"orenden
kontravarianten Filtrumsab\-bildung
$ \chi^\ast \mbox{ \rm Filt}(X) \longrightarrow \mbox{ \rm Filt}A $ geht $F$
auf das Komplement eines Primideals.
\par\smallskip\noindent
(2) Der Halm $({\cal O}_X)_F $ der Strukturgarbe im Filter $F$ ist ein
lokaler Ring.
\end{satz}
Beweis. (1) ist eine direkte Folgerung aus 4.1.1 (4). (2) $ {\cal O}_F $
ist wegen der Konsistenz von $F$, und da $X, {\cal O}_X $ lokal beringt ist,
nicht der Nullring. Seien $f,g \in {\cal O}_F $ und
$f+g \in {\cal O}_F^\times $. Ohne Einschr"ankung k"onnen wir
$f,g \in {\cal O}_X(U)$ und $f+g \in {\cal O}_X(U)^\times $ mit $ U \in F $
annehmen. In dem lokal beringten Raum $ U,{\cal O}_X \mid_U $ gilt dann
$U_f \cup U_g \supseteq U_{f+g} = U$, und wegen der Quasikompaktheit von
$F$ ist $U_f \in F $ oder $U_g \in F $, und damit ist
$ f$ oder $g \in {\cal O}_F^\times $.
\par\bigskip
Mit der Zuordnung $ f \longmapsto X_f $ lassen sich eine ganze Reihe von
Korrespon\-den\-zen zwischen algebraischen und topologischen Objekten definieren.
Das Kom\-plement von $X_f$ wird im lokalen Fall auch die Nullstellen\-menge von
$f$ genannt. Diese Bezeichnung r"uhrt daher, da\3 man auf einem lokal beringten
Raum die Elemente $ f \in A $ als Funktionen auffassen kann, indem man einem
Punkt $ x \in X $ den Wert $ f(x) := \bar{f_x} $ im Restklassenring
$ {\cal O}_x /m_x $
zuordnet; dann ist $ x \not\in X_f $ genau dann, wenn $ f(x) = 0 $ ist.
Zu einem Ideal $ {\bf a} \subseteq A $ und allgemeiner zu jeder
Teilmenge des Ringes definiert man die offene Menge
$ U := \bigcup_{f \in{\bf a}} X_f $, deren Komplement man im lokalen
Fall die gemeinsame Nullstellenmenge der $ f \in {\bf a} $ nennt.
\par
Ist umgekehrt eine offene Menge $U$ in einem lokal(!) beringten Raum gegeben, so
definiert man ein zugeh"origes Ideal durch
$ {\bf a} := \{f \in A: X_f \subseteq U \} $. Dies ist in der Tat ein Ideal;
zun"achst ist $ X_0 = \emptyset $, also $ 0 \in {\bf a}$, und aus
$ X_{f+g} \subseteq X_f \cup X_g $ ergibt sich die Additivit"at, und die
Abgeschlossenheit unter Skalarmultiplikation ergibt sich aus
$ X_{a\cdot f} \subseteq X_f $. Geht man von einer abgeschlossenen Teilmenge
$ Y\subseteq X $ aus, so nennt man das eben zum Komple\-ment von $Y$ definierte
Ideal das Verschwindeideal von $Y$.
\par
Unser Interesse gilt nat"urlich den topologischen und algebraischen Filtern in
einem beringten Raum, die wir im n"achsten Abschnitt systematisch
unter\-suchen wollen.
\section{Topologische und algebraische Filter}
Wir betrachten wieder einen beringten Raum $ X,{\cal O}_X $ mit
$ \Gamma (X,{\cal O}_X) = A $ mit dem zugeh"origen Monoidmorphismus
$ \chi : A \longrightarrow \mbox{ Top}X $, $ f \longmapsto X_f $ und fragen nach dem
Verhalten der Filter unter dieser Abbildung, und mit welchen Eigenschaften
beringter R"aume die Eigenschaft gewisser Filter, Fixfilter zu sein, in
Zusammenhang steht. Zuerst betrachten wir die Umgebungsfilter von Punkten
$ x \in X $.
\begin{satz}
F"ur einen beringten Raum sind folgende Aussagen "aqui\-valent:
\par\smallskip\noindent
(1) Die Umgebungsfilter $U(x)$ sind fix f"ur alle $x \in X $.
\par\smallskip\noindent
(2) Die offenen Mengen $X_f $ bilden eine Basis der Topologie von $X$.
\par\smallskip\noindent
(3) $X$ tr"agt die Initialtopologie unter der kanonischen Abbildung
$ X,{\cal O}_X \longrightarrow \mbox{ \rm Filt}A,\tilde{A} $.
\end{satz}
Beweis. Aus (1) folgt (2). Sei $ x \in U $ offen. Die Voraussetzung besagt
$ U(x) = \chi_\ast \chi^\ast U(x) $, es ist also
$ U \in F( X_f : X_f \in U(x) ) $ und damit gibt es ein $ f \in A $ mit
$ x \in X_f \subseteq U $. Die Invertierbarkeitsmengen bilden also eine Basis.
Aus (2) folgt (3). Da die D(f) eine Basis der Filtrumstopologie bilden, und da
das Urbild von $D(f)$ unter der kanonischen Abbildung gerade $X_f $ ist, ist
die Initialtopologie unter dieser Abbildung die gleiche wie die von den $X_f$
erzeugte; da diese nach Voraussetzung bereits die Topologie von $X$ ist,
ist sie gleich der Initialtopologie. Aus (3) folgt (1). Jede offene Menge
$ U \subseteq X $ ist Urbild einer offenen Menge im Filtrum, ist also vom Typ
$ U = \bigcup_{i \in I} X_{f_i} $. F"ur $ x \in U $ gibt es also ein $ f \in A $
mit $ x \in X_f \subseteq U $. Somit ist $f \in \chi ^\ast U(x) $ und
$ U \in \chi_\ast \chi ^\ast U(x) $, also ist $U(x)$ fix.
\par\bigskip
Einen beringten Raum, der die "aquivalenten Bedingungen des Satzes
erf"ullt,
nennen wir im folgenden einen reich beringten Raum. Ist $X,{\cal O}_X $ ein
lokal beringter Raum, so kann man in der Formulierung (3) das Filtrum durch
das Spektrum ersetzen. Ist $X$ zus"atzlich ein $T_0$-Raum, so ist die
kanonische Abbildung in das Filtrum (oder Spektrum) injektiv und damit eine
Einbettung. Aus $x \neq y $ ergibt sich n"amlich die Existenz einer offenen
Menge $U$ mit $ x \in U , y \not\in U$ (oder umgekehrt) und damit ein
$f \in \Gamma(X,{\cal O}_X) $ mit $ x \in X_f \subseteq U $, also werden $x$ und
$y$ auf verschiedene Filter abgebildet.
\par
Ist $ \phi : X,{\cal O}_X \longrightarrow Y,{\cal O}_Y $
ein Morphismus beringter R"aume, wobei $Y,{\cal O}_Y $ reich beringt
sei und $X$
die Initial\-topologie trage, so ist auch $X,{\cal O}_X $ ein reich
beringter Raum.
Die Voraus\-setzungen ergeben n"amlich in Zusammen\-hang mit Lemma 4.1.2 , da\3
$ X_{\varphi(f)} $, $ f \in \Gamma(Y,{\cal O}_Y) $, eine Basis bilden, also erst
recht alle $X_g $, $ g \in \Gamma(X,{\cal O}_X) $. Speziell ist f"ur einen offenen
Unterraum $ V \subseteq Y $ mit $ Y,{\cal O}_Y $ auch $ V, {\cal O}_Y\mid_V $
reich beringt.
\par\bigskip\noindent
{\bf Beispiel 1} Im Spektrum (und ebenso im Filtrum) eines kommutativen Ring\-es
bilden die Mengen $X_f = D(f)$ eine Basis der Topologie. Affine Schemata sind
demnach reich beringt, und dasselbe gilt f"ur quasiaffine Schemata, also offene
Unterschemata affiner Schemata. Eine projektive Variet"at "uber einem
alge\-bra\-isch abgeschlossenen K"orper $K$ hat als globalen Schnittring $K$ selbst,
und ist daher nur dann reich beringt, wenn sie einpunktig ist.
\par
Ein reich beringtes Schema $X,{\cal O}_X $ ist bereits dann quasiaffin, wenn
$X$ quasikompakt ist. Zu jedem Punkt $x \in X$ gibt es eine offene affine
Umgebung $ x \in U $ und damit ein $ f_x \in A= \Gamma(X,{\cal O}_X) $ mit
$ x \in X_{f_x} \subseteq U $, wobei
 $X_{f_x}$ dann
auch affin ist, da sie mit der Basismenge $D(f_x) $ in $U$ zusammenf"allt.
Deshalb sind auch die offenen Mengen
$ X_{f_x \cdot f_y} = X_{f_x} \cap X_{f_y} $ affin. Aufgrund der
Quasikompaktheit gibt es somit eine endliche "Uberdeckung von $X$ mit offenen
affinen Mengen vom Typ $X_{f_i} $, $ i=1, ... ,n $, mit affinen Durchschnitten.
Unter diesen Voraussetzungen gilt dann aber f"ur beliebiges $ g \in A $ f"ur
die Schnittringe $ \Gamma ( X_g, {\cal O}_X ) = A_g $
(\cite{har}, II.2., exercise 16), und insbesondere haben wir Isomorphien
\begin{eqnarray*}
( X_{f_i} , {\cal O}_X \mid_{X_{f_i}}) \longrightarrow
\mbox{(Spek}A_{f_i}, \tilde{A_{f_i}}) \\ = (D(f_i),\tilde{A} \mid_{D(f_i)}) .
\end{eqnarray*}
Da unter dem kanonischen Morphismus
$ X,{\cal O}_X \longrightarrow \mbox{ Spek}A, \tilde{A} $
die $X_{f_i} $ die Urbilder der Basismengen $D(f_i) $ sind, liegt das Bild von
$X$ unter dieser Abbildung in $D(f_1, ... , f_n) $, und sie stiftet eine
Isomorphie, da die Einschr"an\-kung\-en auf $X_{f_i} $ Isomorphien sind.
\par\smallskip\noindent
{\bf Beispiel 2} Zu einem topologischen Raum $ X$ betrachten wir als
Strukturgarbe die Garbe der stetigen, reellwertigen Funktionen.
Ein topologischer Raum $X$
wird ein $ T_{3a}$-Raum genannt, wenn es zu jeder abgeschlos\-sen\-en Menge
$ A \subseteq X $
und jedem Punkt $ x \not\in A $ eine stetige Funktion
$f: X \longrightarrow [0,1] $ gibt mit
$f(x) =1 $ und $ f \mid_A = 0 $, (\cite{que}, 6.A). Diese Eigenschaft ist "aquivalent dazu, da\3
der beringte Raum $X,{\cal C}_{\bf R}(X) $ reich beringt ist, da man aus einer
stetigen Funktion $f$ mit $ x \in X_f \subseteq U $ leicht eine Funktion $f'$
konstruiert, deren Bild im Einheitsintervall liegt und in $x$ den Wert 1
annimmt.
\par
Ist $X$ zus"atzlich ein $T_0$-Raum, so nennt man einen solchen Raum
vollst"an\-dig regul"ar. Die vollst"andige Regularit"at ist eine starke
Trennungseigenschaft und impliziert beispielsweise die Regularit"at.
Normale R"aume, also insbeson\-dere metrische und kompakte R"aume sind
vollst"andig regul"ar; dies ist der
Inhalt von Urysohn's Lemma, (\cite{que}, 7.A.).
In dem Buch von Gillman/Jerison, Rings of Continuous Functions (\cite{gil}),
in dem es um die beringten R"aume $(X,{\cal C}_{\bf R})$
geht, stehen die vollst"andig regul"aren
R"aume im Mittelpunkt der Untersuch\-ung, 
da nur dann der Ring der stetigen Funktio\-nen interessante Infor\-ma\-ti\-o\-nen "uber
den topologischen Raum bein\-hal\-tet.
\par
Als Unterraum seiner Einpunktkompaktifizierung ist auch ein lokal kom\-pak\-ter
Raum vollst"andig regul"ar, und ein vollst"andig regul"arer Raum $X$ ist genau
dann lokal kompakt, wenn das Bild von $X$ unter der
Stone/$\check{\rm C}$ech\-Kompakti\-fi\-zie\-rung
$ X \longrightarrow \mbox{Spm}C(X) $ offen ist. Somit besteht eine ge\-wis\-se
Ana\-logie zwischen quasi\-affinen Schemata und den beringten R"aumen
$X,{\cal C}(X)$ f"ur lokal kompakte R"aume $X$.
\par\smallskip\noindent
{\bf Beispiel 3} Sei $X$ eine reelle Mannigfaltigkeit und ${\cal O}_X$ die
Garbe der unendlich oft differenzierbaren Funktionen. Dann ist $X,{\cal O}_X $
reich beringt. Ist n"amlich $ x \in U $, so gibt es zun"achst eine zu
$ {\bf R}^n $ hom"oomorphe Umgebung $ x \in V \subseteq U $. Dann gibt es
eine sogenannte Hutfunktion $f$ auf $ {\bf R}^n $, also eine unendlich oft
differenzierbare Funktion mit kompaktem Tr"ager
und mit $f(x) \neq 0 $, (\cite{die}, 16.4.1.4). Diese
Funktion l"a\3t sich au\3erhalb von $V$ mittels der Nullfunktion auf die
gesamte Mannigfaltigkeit differenzierbar fortsetzen und man erh"alt eine
Funktion $ f \in {\cal O}_X(X) $ mit $ x \in X_f \subseteq V $.
\par
Die komplexen Zahlen sind, versehen mit der Garbe der
holomorphen Funktio\-nen,
nicht reich beringt. Beispielsweise gibt es f"ur einen Punkt und eine offene
Ballumgebung keine holomorphe Funktion mit $ x \in X_f \subseteq B $,
da eine solche Funktion auf der offenen Menge $ {\bf C} - \bar{B} $ 
verschwinden m"u\3te, und dann w"are sie nach dem Identit"atssatz bereits
die Nullfunktion. Dies gilt generell f"ur analytische Mannigfaltigkeiten.
\par\bigskip
Aus der kommutativen Algebra kennt man die Aussage, da\3 ein Ideal eines
kommutativen Ringes, das in einer endlichen Vereinigung von Primide\-a\-len
enthalten ist, bereits in einem der Primideale enthalten ist. Sie ist ein
Spezialfall der folgenden Aussage.
\begin{satz}
Ist $X,{\cal O}_X $ ein lokal und reich beringter Raum, dann sind die
Umgebungsfilter endlicher Teilmengen fix.
\end{satz}
Beweis. Die Punkte der endlichen Teilmenge seien mit $x_1, ... ,x_n $
bezeichnet, und es ist zu zeigen, da\3 es f"ur eine offene Menge $U$, in der
diese Punkte liegen, eine Funktion $ f \in A=\Gamma(X,{\cal O}_X) $ gibt mit
$ x_i  \in X_f \subseteq U $ f"ur $i = 1, ... ,n $. Ohne Einschr"ankung
gilt $ U(x_i) \not\subseteq U(x_j) $ f"ur $ i \neq j $, da man sich auf die minimal
umgebenen Punkte beschr"anken kann. Zun"achst gibt es zu jedem Punkt $x_i$
eine Funktion $h_i$ mit $ x_i \in X_{h_i} \subseteq U $. Ferner gibt es zu
jedem $i$ und jedem $ j \neq i $ wegen $ U(x_i) \not\subseteq U(x_j) $ eine
offene Menge
$U_{i,j}$ mit $ x_i \in U_{i,j} , x_j \not\in U_{i,j} $ und damit eine Funktion
$g_{i,j} \in A $ mit $ x_i \in X_{g_{i,j}} , x_j \not\in X_{g_{i,j}} $. F"ur
die Funktion
\begin{displaymath}
 g_i := h_i \cdot \prod_{j \neq i} g_{i,j} 
\end{displaymath}
 ist
$ x_i \in X_{g_i} \subseteq U $ und $ x_j \not\in X_{g_i} $ f"ur $ j \neq i $,
also $ g_i \in {\cal O}_{x_i}^\times $ und $ g_j \in m_{x_i} $. Daraus ergibt
sich mit $ f := g_1 + ... + g_n $ f"ur alle Punkte $x_i$ zun"achst
$ f \in {\cal O}_{x_i}^\times $ und daraus $ x_i \in X_f $ f"ur alle $i$.
Es ergibt sich aber auch, und zwar wieder wegen `lokal beringt' die Inklusion
$X_f \subseteq X_{g_1} \cup ... \cup X_{g_n} \subseteq U $, womit alles
gezeigt ist.
\par\bigskip
Eine extreme Situation liegt vor, wenn s"amtliche topologischen Filter eines
beringten Raumes fix sind. Dies ist "aquivalent dazu, da\3 der Monoid\-mor\-phismus
$ f \longmapsto X_f $ surjektiv ist, da\3 also jede offene Menge
Invertierbar\-keits\-menge einer globalen Funktion ist. Diese Bedingung ist
offenbar hin\-rei\-chend, da f"ur einen surjektiven Monoid\-mor\-phis\-mus stets alle
Filter des Zielmonoids fix sind. Sie ist aber in der Situation eines beringten
Raumes auch notwendig, da f"ur eine offene Menge $ U \subseteq X $ der
Umgebungs\-filter $F(U)$ genau dann fix ist, wenn $U=X_f $ gilt f"ur eine
globale Funktion.
\par\bigskip\noindent
{\bf Beispiel 4} In einem affinen Schema Spek$A, \tilde{A} $ sind genau
dann alle topologi\-schen
Filter fix, wenn $A$ ein Hauptradikalring ist, wenn also jedes Radikal\-ideal
als Radikal von einem einzigen Element erzeugt wird.
\par\smallskip\noindent
{\bf Beispiel 5} Sei $X$ ein metrischer Raum. Dann sind s"amtliche
topologischen Filter fix im beringten Raum $ X,{\cal C}_{\bf R} $. Zu einer
offenen Teilmenge $ U \subseteq X $ mit $A := X-U $ ist die Funktion
$ f(x) := d(x,A) := \mbox{ inf}_{y \in A} d(x,y) $ stetig und verschwindet genau auf
den Punkten aus $A$.
\par\smallskip\noindent
{\bf Beispiel 6} Sei $X= {\bf R}^n $ versehen mit der Garbe der rationalen
Funktionen und der Zariskitopologie. Eine offene Menge $U \subseteq X$ ist
dann Vereinigung endlich vieler Invertierbarkeitsmengen
$ U = \bigcup_{i=1,...,r} X_{f_i} $ mit Funk\-ti\-onen
$ f_i \in \Gamma(X) = {\bf R}[X_1, ... ,X_n]_F $, wobei $F$ aus den
in $X$ nullstellenfreien Polynomfunk\-tionen besteht.
Mit $ f := f_1^2 + ... + f_r^2 $ ist dann
aber $ U =X_f $, da $f$ in einem Punkt $ x \in {\bf R}^n $ genau dann
verschwindet, wenn dort alle $f_i $ verschwin\-den. Also sind auch in dieser
Situation alle topologischen Filter fix.
\par\bigskip
F"ur affine Schemata betrachten wir jetzt die einem topologischen Filter
auf Spek$A$ zugeordneten algebraischen Filter genauer, wobei sich Beziehungen
zu den "Uberlegungen im ersten Kapitel ergeben. Ist
$ Y \subseteq X = \mbox{Spek}A $ eine beliebige Teilmenge, bestehend aus den
Primidealen ${\bf p}_i $, $ i \in I $, so ist der funktionale Umgebungsfilter
einfach gleich $ \bigcap_{i \in I} A-{\bf p}_i $. $ Y \subseteq D(f) $ besagt
ja gerade, da\3 $f$ in keinem ${\bf p_i} $ enthalten ist. Umgekehrt ist der
zu einem Ringfilter $F \subseteq A $ definierte topologische Filter
$F_\ast = F\{ D(f) : f \in F \} $ der Umgebungsfilter der Menge
$Y := $
Spek$A_F = \{ {\bf p} \in \mbox{ Spek}A : {\bf p} \cap F = \emptyset \}$,
da ja $ Y \subseteq D(f) $ genau dann gilt, wenn $f$ zu $F$ geh"ort, was sich
mit dem Lemma von Krull sofort ergibt, vergl. Satz 1.1.1.
\par
Der Umgebungsfilter einer abgeschlossenen Menge $V({\bf a})$ im Spektrum ist fix
und entspricht den Einheiten modulo ${\bf a} $, also dem Filter $F(1+{\bf a}) $.
Ist n"amlich $ V({\bf a}) \subseteq D({\bf b}) $, so "uberdecken
$ D({\bf a})$ und $D({\bf b}) $ das Spektrum, und somit gibt es Elemente
$ g \in {\bf a} $ und $ f \in {\bf b} $ mit $ g+f = 1 $. Damit ist einerseits
$ D(f) \subseteq D({\bf b}) $ und andererseits $ D({\bf a}) \cup D(f) = X $,
also insgesamt $ V({\bf a}) \subseteq D(f) \subseteq D({\bf b})$. Der
Umgebungsfilter der abgeschlossenen Menge hat also eine funktionale Basis und
ist fix. Nimmt man in obiger "Uberlegung speziell $D(f)$ f"ur $D({\bf b}) $, so
gelangt man zur Existenz eines $ x \in A$ mit $ g+ x \cdot f = 1 $, also
$ f \in F(1+{\bf a}) $. Umgekehrt folgt aus $ f \in F(1+{\bf a}) $ zun"achst
$ D(f) \cup D(g) = X $ f"ur ein $ g \in {\bf a} $ und daraus
$V({\bf a}) \subseteq D(f) $, wie behauptet.
\par
Allgemeiner ist der Umgebungsfilter einer Menge
$ Y = \mbox{ Spek}A_F \cap V({\bf a}) $ fix und entspricht dem algebraischen
Filter $ F(F+{\bf a}) $. Das kann man direkt so wie oben zeigen, oder man
f"uhrt es darauf zur"uck, indem man die Situation in Spek$A_F $ betrachtet.
F"ur einen Ringmorphismus $ \varphi : A \longrightarrow B $ ergibt sich die
Faser "uber $ {\bf p} \in \mbox{ Spek}A $ als
\begin{displaymath}
\varphi^ {\ast -1}({\bf p}) =
\mbox{ Spek}B_{\varphi_\ast (A-{\bf p})} \cap V(\langle \varphi({\bf p}) \rangle ).
\end{displaymath}
Der Umgebungsfilter der Faser entspricht also
$ F(\varphi_\ast(A-{\bf p})) + \langle \varphi({\bf p}) \rangle ) $, und der
Urbildfilter des Umgebungsfilters von $ {\bf p} $ ist ebenfalls fix und
entspricht dem Filter $ \varphi_\ast (A-{\bf p}) $. Da eine stetige Abbildung
genau dann abgeschlossen ist, wenn f"ur jeden Punkt des Bildraumes der
Urbildfilter des Umgebungsfilters mit dem Umgebungsfilter der Faser
"ubereinstimmt, und da in unserer Situati\-on wegen der Fixheit beider Filter
dies algebraisch getestet werden kann, ist f"ur einen Ringmorphismus wie oben die
Spektrumsabbildung genau dann abgeschlossen, wenn f"ur jedes Primideal
${\bf p} \in \mbox{ SpekA} $ der Filter $ \varphi_\ast (A-{\bf p}) $ fix
modulo $ \langle \varphi ({\bf p}) \rangle $ ist.
\par
Der Umgebungsfilter einer offenen Menge im Spektrum ist nur dann fix, wenn es
sich um eine Basismenge handelt. Der algebraische Filter $F$, der dem
Umgebungsfilter einer offenen Menge $ U = D({\bf a}) $ zugeordnet wird,
l"a\3t sich schreiben als
\begin{eqnarray*}
F = F(U)= \bigcap _{{\bf p} \in U} A-{\bf p} =
\{ g \in A : {\bf a} \subseteq {\rm Rad}(g)\}
\\ = \bigcap_{D(f) \subseteq U} F(f) = \bigcap _{\bigcup D(f_i) = U} F(f_i).
\end{eqnarray*} 
Dabei ist die erste Beschreibung f"ur beliebiges $Y$ richtig, insbesondere
also auch f"ur offenes $U$, und die zweite ist nat"urlich auch richtig, also
sind die ersten beiden gleich. Aufgrund von Lemma 1.2.8 stimmen die drei
letzten Filter "uberein.
\par
Die Zuordnung $ U \longmapsto A_{F(U)} $ ist eine Pr"agarbe auf Spek$A$, da f"ur
$V \subseteq U $ die Inklusion $ F(U) \subseteq F(V) $ den
Ringmorphismus $ A_{F(U)} \longrightarrow A_{F(V)} $ liefert. Ist
$ U=D(f) $, so ist $ F(U) =F(f) $, und somit stimmt diese Pr"agarbe auf den
Basismengen mit der Strukturgarbe "uberein, die damit die Garbe zu dieser
Pr"agarbe ist.
\par
Die Abbildungen $ A_{F(U)} \longrightarrow \Gamma(U,\tilde{A}) $,
die sich aus der Vergarbung der Pr"agarbe ergeben, sind lokal, d.h. $F(U)$ ist
das Urbild der Einheitengruppe unter dem Restriktionsmorphismus
$R: A \longrightarrow \Gamma(U,\tilde{A}) $, da ja f"ur ein $f \in A $ gilt
$ U \subseteq D(f) $ genau dann, wenn $ f \in \Gamma(U,\tilde{A})^\times $.
Insbesondere gilt, da\3, wenn $\Gamma(U,\tilde{A}) $ ein Bruchring von $A$ ist,
es sich bereits um $A_{F(U)} $ handeln mu\3. Die Aussage 1.2.11 ist ein
Spezialfall dieser "Uberlegungen, da f"ur einen Integrit"atsbereich $A$ die
Teilerfremdheit von $f,g \in A^\ast$ gerade
$\Gamma(D(f) \cup D(g), \tilde{A}) = A $
besagt, woraus sich dann $F(f) \cap F(g) = A^\times $ ergibt.
\par
Es erhebt sich die Frage, unter welchen Bedingungen diese Pr"agarbe eine Garbe
ist und welche Konsequenzen das hat. F"ur einen Integrit"atsbereich haben wir
folgendes Resultat.
\begin{lem}
Ist $A$ ein Integrit"atsbereich mit Quotientenk"orper $Q=Q(A)$, dann ist die
Pr"agarbe genau dann eine Garbe, wenn der maximale Definiti\-ons\-bereich jeder
Funktion aus dem Quotientenk"orper eine Basismenge ist.
\end{lem}
Beweis. F"ur ein $g \in Q(A) $
versteht man unter dem maximalen oder exakten Definitionsbereich Def($q$)
die Vereinigung der
Basismengen $D(f)$ "uber alle $f \in A$, so da\3 es eine Darstellung $a/f = q $
gibt. Diese offene Menge ist identisch mit der Menge der Primideale
${\bf p} \in \mbox{ Spek}A $ mit $ q \in A_{\bf p} $. Sei die Pr"agarbe eine
Garbe und $q \in Q$ vorgegeben mit
$ U =\mbox{ Def}(q) = \bigcup_{i \in I} D(f_i) $. Die rationale Funktion $q$ ist
auf $U$ definiert, also $ q \in \Gamma(U,{\cal O}_X) = A_{F(U)} $, also
$q =a/f $ mit $f \in \bigcap_{i \in I} F(f_i) $. Damit ist $U \subseteq D(f) $,
wobei $D(f)$ aber auch nicht gr"o\3er sein kann, da ja $U$ bereits der maximale
Definitionsbereich von $a/f $ ist.
\par
Sei umgekehrt der Definitionsbereich jeder rationalen Funktion eine Basis\-menge,
und $U$ eine offene Menge. $B:=\Gamma(U,{\cal O}_X)$ besteht aus allen rationalen
Funktionen, die auf $U$ definiert sind, also $q \in B$ genau dann, wenn 
$U \subseteq \mbox{ Def}(q) $; ist dies der Fall, so gibt es nach Voraussetzung
eine Darstellung $q=a/f$ mit $U \subseteq D(f)$, und $q$ liegt bereits in dem
Pr"agarbenring, und somit sind die Vergarbungsmorphismen surjektiv und auch
bijektiv.
\par\bigskip
Die folgende Aussage
gibt f"ur einen Integrit"atsbereich eine hinreichende Bedingung daf"ur, da\3 die
Pr"agarbe eine Garbe ist. Sie ist zugleich
eine Verallgemeinerung von Satz 1.2.12 .
\begin{satz}
Sei $A$ ein Integrit"atsbereich, in dem jedes Element des Quotien\-tenk"orpers
$Q$ eine Darstellung mit teilerfremden Z"ahler und Nenner besitzt. Dann ist die
oben definierte Pr"agarbe bereits eine Garbe und somit gleich der Strukturgarbe.
\end{satz}
Beweis. Wir haben lediglich zu zeigen, da\3 f"ur
$ U = D({\bf a}) = D( f_i : i \in I) $ der Ringmorphismus
$A_F(U) \longrightarrow \Gamma(U,\tilde{A}) = \bigcap_{i \in I} A_{f_i} $
bijektiv ist, wobei der Durchschnitt im Quotientenk"orper zu bilden ist.
Die Injektivit"at ist dabei klar, sei also $a/g \in Q $ vorgegeben, mit
$ a/g \in A_{f_i} $ f"ur alle $i \in I$. Dabei seien $a$ und $g$ teilerfremd.
Aus den Darstellungen $ a/g = a_i/f_i^n $, also $ a \cdot f_i^n = a_i \cdot g $
ergibt sich aufgrund der Teilerfremdheit sofort, da\3 $f_i^n$ von $g$ geteilt
wird, also $g \in F(f_i) $ f"ur alle $ i  \in I $, also
$ a/g \in A_{ \bigcap_{i \in I} F(f_i) } = A_{F(U)} $
\par\bigskip
Die Bedingung des Satzes ist beispielsweise dann erf"ullt, wenn $A$ faktoriell
ist, da in diesem Fall ein Element aus dem Quotientenk"orper in eine Form
gebracht werden kann, wo Z"ahler und Nenner keine gemeinsamen Primfakto\-ren
haben, oder in einem Bezoutbereich, da es dort f"ur ein $q \in Q $ sogar eine
Darstellung $a/g$ gibt mit $ Aa + Ag = A $. Die Aussage 1.2.12 ergibt sich
hieraus als Korollar, da unter der angegebenen Bedingung speziell
$ A_{F(f) \cap F(g)} = A_f \cap A_g $ gilt, und der linke Ring ist gleich $A$,
wenn $F(f) \cap F(g) =A^\times $ ist, und der rechte, wenn $f$ und $g$
teilerfremd sind.
\par\bigskip
Als n"achstes wollen wir f"ur den Fall, da\3 $A$ ein noetherscher, normaler
Integrit"atsbereich ist, sehen, wie sich die Eigenschaft, da\3 die Pr"agarbe
eine Garbe ist, auf die Divisorenklassengruppe auswirkt. Hierzu sind zwei
Lemmata n"otig.
\begin{lem}
Sei $A$ noethersch, normal und integer mit Quotientenk"orper $Q$. Dann ist der
Definitionsbereich einer rationalen Funktion $f \in Q$ gleich einem endlichen
Durchschnitt von Komplementen irreduzibler Hyperfl"achen, und zwar ist
\begin{displaymath}
\mbox{\rm Def}(f) = \bigcap_{ht({\bf p}) =1 , \hspace{0.5em} \nu_{\bf p}(f) <0 } D({\bf p}),
\end{displaymath}
wobei $\nu_{\bf p}$ die Bewertung f"ur ein Primideal der H"ohe eins bezeichnet.
\end{lem}
Beweis. `$\subseteq$' Sei $ f \in A_{\bf q} $ und
${\bf p}$ ein Primideal der H"ohe 1 mit $\nu_{\bf p}(f) <0 $. Dann ist
$ f \not\in A_{\bf p} $, also $A_{\bf q} \not\subseteq A_{\bf p} $, also
${\bf p} \not\subseteq {\bf q} $, also ${\bf q} \in D({\bf p}) $.
\par
`$\supseteq$' Sei ${\bf q}$ im rechten Schnitt, also
${\bf q} \not\supseteq {\bf p}_i$, wo die ${\bf p}_i$ die (endlich vielen)
Polstellen von $f$ sind.
$A_{\bf q}$ ist als Lokalisierung eines normalen Rings wieder normal, und daher
Durchschnitt seiner Lokalisierungen der H"ohe 1. Es ist also
\begin{eqnarray*}
A_{\bf q} =
\bigcap_{ ht({\bf r})=1 , \hspace{0.5em}{\bf r} \in \mbox{ Spek}A_{\bf q} } (A_{\bf q})_{\bf r}
\\ = \bigcap_{ ht({\bf r})=1, \hspace{0.5em}{\bf r} \subseteq {\bf q} } A_{\bf r}  .
\end{eqnarray*}
Aufgrund der Voraussetzung "uber ${\bf q}$ sind die $ {\bf r} \subseteq {\bf q}$
von allen ${\bf p}_i $ verschieden, und damit ist $\nu_{\bf r}(f) \ge 0 $, also
$ f \in A_{\bf r} $ f"ur alle ${\bf r} \subseteq {\bf q} $, und damit
$ f \in A_{\bf q} $. 
\begin{lem}
Sei $A$ noethersch, normal und integer mit dem Quo\-tien\-ten\-k"or\-per $Q$.
${\bf p}_1, ... , {\bf p}_n $ seien Primideale der H"ohe 1. Dann gibt es eine
rationale Funktion $f \in Q$, deren Defini\-tions\-bereich gerade gleich
$ D({\bf p}_1) \cap ... \cap D({\bf p}_n) $ ist.
\end{lem}
Beweis. Wir zeigen, da\3 es eine rationale Funktion gibt, die genau an den
vorgegebenen Primidealen eine Polstelle hat, woraus aufgrund des
vorausge\-gang\-enen Lemmas die Behauptung folgt. Sei $g \in A $ mit
$0 \neq g \in \bigcap_{i=1,...,n} {\bf p}_i $, also
$\nu_{{\bf p}_i}(g^{-1}) <0$.
${\bf q}_1, ... ,{\bf q}_k $ seien die endlich vielen anderen Primideale der
H"ohe 1 mit $ \nu_{{\bf q}_j }(g^{-1}) < 0 $. F"ur alle $i,j$ ist
${\bf q}_j \not\subseteq {\bf p}_i $, also ${\bf p}_i \in D({\bf q}_j )$, und
damit ist $U:=D({\bf q}_1) \cap ... \cap D({\bf q}_k) $ eine offene Umgebung der
${\bf p}_i , i=1,...,n $, und es ist ${\bf q}_j \not\in U $ f"ur $j=1,...,k$.
 Damit gibt es in $U$ auch eine funktionale Umgebung
dieser Punkte, also ein $h \in A$ mit $ {\bf p}_i \in D(h) $ und
${\bf q}_j \not\in D(h)$. Es ist also $ \nu_{{\bf p}_i}(h)=0 $ und
$ \nu_{{\bf q}_j} (h) >0 $. Durch Wahl einer geeigneten Potenz von $h$ kann man
erreichen, da\3 $ \nu_{{\bf q}_j}(h^n \cdot g^{-1}) $ nichtnegativ ist f"ur alle
$j=1,...,k $. Wir setzen also $f := h^n \cdot g^{-1}$ und haben damit 
$\nu_{{\bf q}_j}(f) \ge 0 $, f"ur ein ${\bf p}_i $ ist
$ \nu_{{\bf p}_i} (f) = \nu_{{\bf p}_i}(g^{-1}) <0 $ , und f"ur ein sonstiges
Primideal ${\bf p}$ der H"ohe 1 haben wir
$ \nu_{\bf p}(f) = n\cdot\nu_{\bf p}(h) + \nu_{\bf p}(g^{-1}) \ge \nu_{\bf p}(g^{-1}) \ge 0 $.
$f$ hat also genau an den vorgegebenen Primidealen negative Divisorwerte.
\begin{satz}
F"ur einen noetherschen, normalen Integrit"atsbereich $A$ ist die Pr"agarbe genau
dann gleich der Strukturgarbe auf dem Spektrum, wenn $A$ fastfaktoriell ist, d.h.
die Divisorenklassengruppe eine Torsionsgruppe ist.
\end{satz}
Beweis. Die beiden obigen Lemmata haben gezeigt, da\3 unter den
Voraus\-set\-zung\-en
die Menge der Definitionsbereiche rationaler Funktionen und die Men\-ge der
Divisorbereiche, also Mengen vom Typ $D({\bf p}_1)\cap...\cap D({\bf p}_n)$ mit
Prim\-ide\-a\-len der H"ohe 1, "ubereinstimmen. Da die Pr"agarbe genau dann eine
Garbe ist, wenn s"amtliche Definitionsbereiche rationaler Funktionen
Basismengen sind, ist zu zeigen, da\3 die Fastfaktorialit"at zu der Eigenschaft
"aquivalent ist, da\3 alle Divisorenbereiche Basismengen sind.
\par
Sei $A$ fastfaktoriell und ${\bf p}$ ein Primideal der H"ohe 1, und $D$ der
Divisor, der genau an dieser Stelle den Wert 1 hat, und sonst "uberall 0. Da
die Divisoren\-klassen\-gruppe torsorisch ist, ist ein Vielfaches von $D$ ein
Haupt\-divisor einer ratio\-nalen Funktion $f$, die "uberdies
(da der Divisor nicht negativ ist) aus $A$ sein mu\3. Damit ist
Rad$(f) ={\bf p}$ und $D(f) =D({\bf p}) $, und somit sind auch die Durch\-schnitte
von Kom\-ple\-menten irre\-duzi\-bler Hyper\-fl"achen Basis\-mengen.
\par
Ist umgekehrt $D({\bf p}) =D(f) $, so hat der Divisor zu $f$ genau an der Stelle
${\bf p}$ einen positiven Wert, und damit verschwindet ein Vielfaches des
Basisdivisors an der
Stelle ${\bf p}$ in der Divisorenklassengruppe, und mit diesem Erzeugendensystem
ist diese insgesamt eine Torsionsgruppe.
\par\bigskip
Die Frage, welche offenen Teilmengen eines affinen Schemas selbst affin
sind, hat im Fall, da\3 Strukturpr"agarbe und
Strukturgarbe "ubereinstimmen, die gleiche einfache Antwort wie im
faktoriellen Fall.
\begin{satz}
Sei $A$ ein kommutativer Ring, in dem die Pr"agarbe auf
dem Spek\-trum eine Garbe ist. Dann sind die
Basismengen die einzigen affinen Teil\-men\-gen von {\rm Spek}$A$. 
\end{satz}
Beweis. Sei $ U = D({\bf a}) = \bigcup_{i \in I} D(f_i) $ eine
offene, affine Menge im
Spektrum. D.h. die kanonische Abbildung
\begin{eqnarray*}
U,\tilde{A} \mid _U \longrightarrow \mbox{ \rm Spek} \Gamma(U,\tilde{A}) =
\mbox{ Spek}A_F
\end{eqnarray*}
mit $ F :=  \bigcap_{i \in I} F(f_i) $
ist ein Isomorphismus, also insbesondere surjektiv. W"are $ {\bf a} \cap F $
leer, gebe es nach dem Lemma von Krull ein Primideal ${\bf p}$ mit
$ {\bf a} \subseteq {\bf p} $ und ${\bf p} \cap F = \emptyset $, also 
$ {\bf p } \in \mbox{ Spek}A_F $, aber ${\bf p} \not\in U $, was nach
Voraussetzung nicht sein kann. Sei also $ f \in {\bf a} \cap F $. Aus
$ f \in {\bf a} $ folgt $ D(f) \subseteq D({\bf a}) $, und $ f \in F $
bedeutet $ f \in F(f_i) $ f"ur alle $i \in I $, also
$ D({\bf a}) \subseteq D(f) $.
\par\bigskip
H"aufig spielt in beringten R"aumen die Eigenschaft, da\3 sich je zwei Punkte
durch eine globale Funktion trennen lassen, eine wichtige Rolle. In der Theorie
der komplexen R"aume spricht man von holomorph separabel , und im
Appro\-xi\-mationssatz von Stone-Weierstra\3 von punktetrennenden Unteralge\-bren.
Dabei kann
diese Eigenschaft in der allgemeinen Situation eines bering\-ten Raumes in
verschiedener Weise pr"azisiert werden; man kann die
Injektivi\-t"at der kanonischen Abbildung in das Spektrum (bzw. Filtrum) zur
Defini\-ti\-on nehmen, oder aber f"ur $ U(x) \not\subseteq U(y) $ die Existenz
einer Funk\-tion $f$ verlangen mit $ x \in X_f , y \not\in X_f $, oder man kann
sich
auf die abgeschlossenen Punkte beschr"anken wollen. All diese Formulierungen
sind schw"acher als reich beringt (das $T_0$-Axiom vorausgesetzt), und
Um\-kehr\-ung\-en kann es allen\-falls in Spe\-zial\-f"al\-len ge\-ben.
Wir geben zun"achst ein Bei\-spiel eines Schemas $X,{\cal O}_X $, wo
der globale Schnitt\-ring $A$ die Punkte von $X$ trennt, aber $X$ trotz\-dem nicht
reich be\-ringt ist.
\par\bigskip\noindent
{\bf Beispiel 8}
Sei  $ X:= \mbox{ Spek}{\bf Z} $ als Menge, aber
versehen mit folgender Topolo\-gie: neben der leeren Menge sei jede Teilmenge
von $X$ offen, die das Nullideal enth"alt. Als Garbe definieren wir
$ \Gamma(U) := {\bf Z}_{F(U)} $ mit
$F(U) = \{ f \in {\bf Z} : f \not\in {\bf Z}p $ f"ur alle ${\bf Z}p \in U \}$.
Dies ist aufgrund der
Faktorialit"at von ${\bf Z}$ eine Garbe, und die Halme sind die des richtigen
Spektrums von ${\bf Z}$. F"ur eine Primzahl $p$ ist
$ U= \{ \langle 0 \rangle , \langle p \rangle \} $ eine offene Umgebung mit
dem Schnittring $\Gamma(U,{\cal O}_X) = {\bf Z}_{ {\bf Z}-{\bf Z}p } $, und damit
ist $U$ isomorph zum Spektrum dieses diskreten Bewertungsringes, und damit ist
$X$ ein Schema. Die kanonische Abbildung
$X,{\cal O}_X \longrightarrow \mbox{ Spek}{\bf Z},\tilde{\bf Z} $
ist bijektiv, aber keine Einbettung, da etwa das Null\-ide\-al im Spektrum nicht
offen ist.
\par\bigskip
Das letzte Beispiel war nicht noethersch, und ich wei\3 nicht, ob es auch ein
noethersches Beispiel f"ur diese Situation gibt. Insbesondere w"are es
interessant zu wissen, ob Variet"aten, deren abgeschlossene Punkte durch globale
Funktionen getrennt werden, reich beringt und da\-mit be\-reits
qua\-si\-affin sind. Im Eindimensionalen haben wir folgendes Resultat.
\begin{satz}
Ist $X,{\cal O}_X$ ein eindimensionales, integres, noethersches Schema, wo
die globalen Funktionen die abgeschlossenen Punkte trennen, so ist $X$ affin.
\end{satz}
Beweis. Zun"achst sei daran erinnert, da\3 in einem affinen, noetherschen,
eindimensionalen, integren Schema neben dem gesamten Raum nur die end\-li\-chen
Teilmengen von maximalen Idealen abgeschlossen sind. Dasselbe gilt dann auch
in unserer Situation, da $X$ von endlich vielen affinen offenen Mengen von
diesem Typ "uberdeckt wird, und eine abgeschlossene Menge $ Y \neq X$
in jeder von diesen nur endlich viele Punkte hat. Ist nun $ x \in U $ ein
abgeschlossener Punkt und $U$ offen, so gibt es zu jedem der endlich vielen
abgeschlossenen Punkte au\3erhalb von $U$ jeweils eine globale Funktion, die
in $x$ invertierbar ist, in jenem Punkt aber nicht
\footnote{Wir verstehen hier also unter ``die abgeschlossenen Punkte trennen'',
da\3 es zu abgeschlossenen Punkten $x \neq y $ eine globale Funktion $f$ gibt
mit $x \in X_f , y \not\in X_f $ und eine globale Funktion $g$ mit
$ x \not\in X_g , y \in X_g $. In vielen Situationen folgt
das eine aus dem andern.},
und f"ur das Produkt $f$ dieser
Funktion gilt $ x \in X_f \subseteq U $, und in $X_f$ liegt auch der generische
Punkt, und somit ist $X$ reich beringt und nach Beispiel 1 quasiaffin.
\par
Sei also
$X \subseteq \mbox{ Spek}A $ ein offenes Unterschema. Wir geben einen
kohomo\-lo\-gischen Beweis f"ur die Affinit"at von $X$, indem wir
$H^1(X,{\cal G}) = 0$
f"ur jede quasi\-koh"arente Garbe ${\cal G}$ auf $X$ zeigen.
$ \widetilde{H^1(X,{\cal G})} $ ist auf Spek$A$ gleich der ersten h"oheren
Bildgarbe $R^1i_\ast({\cal G})$ (\cite{har} III.8). F"ur alle
${\bf m} \in$ Spm$A$
ist also nach\-zu\-weisen, da\3
$ \lim_{U \in \psi({\bf m})} H^1(U,{\cal G}) =0 $ gilt.
Geh"ort ${\bf m}$ zu $X$, so ist
$\psi({\bf m}) = U({\bf m}) $ in $X$, und damit
ist der Limes gleich null. Geh"ort ${\bf m}$ nicht zu $X$, so besteht
$\psi({\bf m})$ aus
allen nichtleeren offenen Teilmengen von $X$, da eine solche Komplement einer
endlichen Teilmenge in $X$ ist und damit auch Komplement einer endlichen
Teilmenge in Spek$A$; f"ugt man zu dieser offenen Menge den Punkt ${\bf m}$ hinzu,
erh"alt man eine offene Umgebung von ${\bf m}$, deren Urbild in $X$ die vorgegebene
ist. Die Kohomologie in einem irreduziblen Filter ist aber null.
\par
F"ur den letzten Teil des Beweises wird noch eine Variante folgen, die auf
Kohomologie verzichtet.
\section{Nullstellensatz und algebraische Filter}
Der Hilbertsche Nullstellensatz besagt f"ur einen beringten Raum
$X,{\cal O}_X $, da\3 f"ur beliebige Elemente
$g, f_i \in A = \Gamma(X,{\cal O}_X) $ aus
$X_g \subseteq \bigcup_{i \in I} X_{f_i} $ stets $ g \in {\rm Rad}(f_i : i \in I) $
folgt. Er gilt beispielsweise f"ur affine algebraische Variet"aten "uber einem
algebraisch abgeschlossenen K"orper, und f"ur jedes affine Schema. Der
Nullstellensatz ist eine recht starke Anforderung an einen beringten Raum,
zumindest, wenn die durch die globalen Funktionen definierte Topologie
einigerma\3en fein ist.
\par
In diesem
Abschnitt geht es um einige Spezialf"alle des Nullstellensatzes, und
insbesondere um deren Beziehung zu der Frage, welche algebraischen Filter
fix sind. Betrachtet man in der obigen Formulierung statt einer belie\-bi\-gen
Familie $f_i$ jeweils nur ein einziges $f$, so lautet die Aussage:
wenn $X_g \subseteq X_f $, dann ist $g \in {\rm Rad}(f) $, also $ f \in F(g) $
f"ur beliebige $f,g \in A$.
\begin{satz} Sei $ X,{\cal O}_X $ ein beringter Raum mit
globalem Schnittring $A$.
Dann sind "aquivalent:
\par\smallskip\noindent
(1) F"ur $f ,g \in A $ folgt aus $X_g \subseteq X_f $ die Beziehung
{\rm Rad}$(g) \subseteq $ {\rm Rad}$(f) $.
\par\smallskip\noindent
(2) Bez"uglich $ \chi: A \longrightarrow \mbox{ \rm Top}X $,
$f \longmapsto X_f $
sind
alle Filter aus $A$ fix.
\par\smallskip\noindent
(3) Die kanonische Abbildung $ \phi: X \longrightarrow {\rm Filt}A $
ist filterhaft.
\end{satz}
Beweis. (1) ist "aquivalent zu  (2). F"ur festes $g \in A$ besagen beide Aussagen,
da\3 f"ur beliebiges $f \in A $ aus $X_g \subseteq X_f $ im Ring $ f \in F(g) $
folgt, da man sich bei (2) auf die Hauptfilter $F(g)$ beschr"anken kann.
Aus (1) folgt (3). Wir behaupten, da\3 jede Basismenge im Filtrum
die f"ur die Filterhaftigkeit charakteristische Eigenschaft hat.
Sei $ D(f) $ vorgegeben, und sei $ D(g) $ eine andere Basismenge mit
$\phi^{-1}(D(g)) \subseteq \phi^{-1}(D(f))$. Das hei\3t $X_g \subseteq X_f$ und
daraus folgt nach Voraussetzung $ f \in F(g) $, also $ D(g) \subseteq D(f) $.
Aus (3) folgt (1). Da die kanonische Abbildung filterhaft ist, gibt es
insbesondere f"ur die Situation $ F(f) \in D(f) $ eine offene, filterhafte
Menge $V$ mit
$F(f) \in V \subseteq D(f) $, wof"ur aber nur $D(f)$ selbst in Frage
kommt. Also gilt f"ur beliebiges $g \in A$ die Implikation aus
$ X_g \subseteq X_f $ folgt $D(g) \subseteq D(f) $, also $f \in F(g) $.
\par\bigskip
{\bf Bemerkung} In der Formulierung (2) kann man sich auf die Primideale von $A$
beschr"anken, ebenso wie auf die Hauptfilter. Dagegen gen"ugt es, wie mir
scheint, im Fall eines lokal beringten Raumes nicht, in (3) nur die
Filterhaftigkeit der kanonischen Abbildung in das Spektrum zu fordern.
\begin{pro}
Sei $X,{\cal O}_X $ ein beringter Raum. In folgenden F"allen sind f"ur $X$ alle
algebraischen Filter fix.
\par\smallskip\noindent
(1) $X$ ist ein affines Schema.
\par\smallskip\noindent
(2) $X$ ist ein quasikompaktes und integres Schema.
\par\smallskip\noindent
(3) F"ur jede Funktion $f\in A$ aus dem globalen Schnittring ist
$\Gamma(X_f,{\cal O}_X) = A_f $. In diesem Fall ist
$ \psi: \mbox{ \rm Filt}A,\tilde{A} \longrightarrow \mbox{ \rm Filt}X,i_\ast{\cal O}_X $
sogar eine halmtreue Einbettung.
\end{pro}
Beweis. (2) Sei $X_g \subseteq X_f $. Ohne Einschr"ankung k"onnen wir $f$ als
von null verschieden annehmen, und damit ist $f$ auch in jedem Schnittring
"uber einer nichtleeren Menge ungleich null. Es sei
$X= \bigcup_{i \in I}U_i $ eine
endliche affine "Uberdeckung. In $U_i$ gilt $D(g_i) \subseteq D(f_i)$, wobei der
Index die Restriktion auf $U_i$ anzeigt, und man hat in
$A_i = \Gamma(U_i,{\cal O}_X)$ wegen der
Affinit"at von
$U_i$ die Beziehung $(g_i)^{n_i} = f_i \cdot a_i $. Wegen der Endlichkeit von
$I$ k"onnen wir alle $n_i$ als gleich annehmen. Dann gelten in
$ \Gamma(U_i \cap U_j,{\cal O}_X) $ die Beziehungen $ g^n = f \cdot a_i $ und
$g^n = f \cdot a_j $, also $ 0 = f( a_i -a_j) $, woraus aus der Voraussetzung
"uber $f$ und der Integrit"at $a_i =a_j$ folgt, also die Vertr"aglichkeit der
$a_i$, und man gelangt so zu einer globalen Gleichung $ g^n= f \cdot a $.
\par
(3) Unter der Voraussetzung ist nat"urlich auch zu verstehen, da\3 die
kanonischen Abbildungen $A \longrightarrow A_f $ die
Restriktionsmorphismen sind. Aus $ X_g \subseteq X_f $ ergibt sich dann der
Restriktionsmorphismus $ A_f \longrightarrow A_g $, der "uber $A$ vertr"aglich
ist, und daher folgt $f \in F(g) $, und die angegebene Abbildung ist eine
Einbettung.
Ferner ist f"ur einen Filter $ F\subseteq A$
\begin{eqnarray*}
\lim_{U \in \psi(F)} \Gamma(U) = \lim_{f \in F} \Gamma(X_f) =
\lim_{f \in F}A_f = A_F ,
\end{eqnarray*}
und damit ist Filt$A$ nicht nur topologisch, sondern auch als beringter Raum
ein Unterraum von Filt$X$.
\par
F"ur ein Schema ist die Bedingung erf"ullt,
wenn es eine endliche affine "Uberdeckung gibt, deren paarweise Durchschnitte
quasikom\-pakt sind. Das ist beispielsweise dann der Fall, wenn das Schema
noethersch oder quasikom\-pakt und separiert ist.\footnote{Es haben also sehr viele
Schemata die Eigenschaft, da\3 alle algebraischen Filter fix sind; ich wei\3
aber nicht, ob dies f"ur alle Schemata gilt und w"are f"ur einen Beweis oder
ein Gegenbeispiel dankbar.}
\par\bigskip
H"aufig sind in einem beringten Raum nur die maximalen Ideale (als Filter) fix,
aber nicht alle algebraischen Filter und damit nicht alle Primideale. Diese
Eigenschaft ist "aquivalent damit, da\3 aus einer Inklusion $X_f \subseteq X_g$
im globalen Ring die Beziehung $f \in \mbox{ Jak}(g) $ folgt, da\3 also $f$ in
allen maximalen Idealen enthalten ist, die $g$ umfassen. Ist der globale Ring
eines beringten Raumes eine $K$-Algebra von endlichem Typ, $K$ K"orper, so ist
diese Eigenschaft "aquivalent zu der im allgemeinen st"arkeren, da\3 alle
algebra\-isch\-en Filter fix sind, da ja in dieser Situation die Jakobsonradikale
mit den Ra\-di\-ka\-len "ubereinstimmen.
\par\bigskip\noindent
{\bf Beispiel 1} Ist $X$
ein topologischer Raum, so sind alle maximalen Ideale aus
$C(X,{\bf R})$ fix. Sei
n"amlich ${\bf m} $ ein maximales Ideal, $ f \not\in {\bf m} $,
$X_f \subseteq X_g $. Es gibt dann Funktionen $h \in {\bf m}$ und $a$ mit
$a \cdot f + h =1 $ und somit ist $ X_f \cup X_h = X $, also auch
$X_g \cup X_h = X$. Dies gilt dann auch f"ur die Quadrate und wir erhalten
$ X_{g^2 + h^2} = X_{g^2} + X_{h^2} = X $, und damit ist $g^2 + h^2 $ eine
Einheit, und $g$ kann nicht zu ${\bf m}$ geh"oren.
\par\smallskip\noindent
{\bf Beispiel 2}
Die gleiche Argumentation l"a\3t sich auch f"ur die Garbe der rationalen
Funktionen auf ${\bf R}^n $ durchf"uhren; sie zeigt, da\3 auch in dieser
Situation alle maximalen Ideale fix sind. Es sind aber nicht alle
Filter fix. Seien $f = X^2 + Y^2 $ und
$g = X^2 + Y^4 $. Beide Funktionen haben in ${\bf R}^2 $ die gleiche
Invertierbarkeitsmenge ${\bf R}^2 -\{(0,0)\}$ und sind prim "uber ${\bf R}$, wie 
man sieht, wenn man die Zerlegungen $f =(X+iY)\cdot(X-iY) $ und
$g=(X+iY^2)\cdot(X-iY^2)$ "uber ${\bf C}$ betrachtet. $f$ und $g$ bleiben im
globalen Schnittring ${\bf R}[X,Y]_F $,
wobei $F$ die Menge der
nullstellenfreien Polynome auf ${\bf R}^2 $ ist, prim, und es ist daher
$F(f) \neq F(g)$. Man beachte, da\3 in diesem Beispiel der globale Schnittring
nicht von endlichem Typ ist.
\begin{satz}
Sei $X,{\cal O}_X $ ein lokal beringter Raum mit globalem Schnittring $B$ und
$\varphi : X,{\cal O}_X \longrightarrow \mbox{ \rm Spek}A $ ein Morphismus
mit dem globalen Ringmorphis\-mus $ R: A \longrightarrow B $. Dann gilt
\par\smallskip\noindent
(1) Das Diagramm
\par
\unitlength0.5cm
\begin{picture}(12,7)
\put(4,5){\vector(1,0){3}}
\put(2.8,4){\vector(0,-1){2}}
\put(4.5,1.2){\vector(1,1){3}}
\put(2.6,4.7){$X$}
\put(7.7,4.7){{\rm Spek}$A$}
\put(1.3,1){{\rm Spek}$B$}
\put(2.2,2.4){$j$}
\put(6.5,2.3){$R^\ast$}
\put(5.5,5.5){$\varphi$}
\end{picture}
\par\noindent
kommutiert.
\par\smallskip\noindent
(2) Ist $ {\bf p} \in \mbox{ \rm Spek}B $, so ist
$\varphi(U({\bf p})) \supseteq U(R^\ast({\bf p})) $.
Insbesondere konvergie\-ren
also die Bildfilter von $ U({\bf p}) = F( X_f : f \not\in {\bf p} ) $ in einem
affinem Schema.
\par\smallskip\noindent
(3) Sind in $X$ die algebraischen Filter fix,
so stimmen die beiden Filter in (2) zu einem Primideal
${\bf p} \in \mbox{ \rm Spek}B$ funktional "uberein.
\end{satz}
Beweis. (1) $\varphi $ und $ R^\ast \circ j $ sind beide Morphismen lokal
beringter R"aume in ein 
affines Schema mit dem gleichen globalen Ringmorphismus und sind daher gleich.
\par\noindent
(2) Die $D(g), g \not\in R^\ast ({\bf p})$, bilden eine Umgebungsbasis von
$R^\ast({\bf p}) $. F"ur ein solches $g$ ist also $ R(g) \not\in {\bf p} $ und
damit $ X_{R(g)} \in U({\bf p}) $. Andererseits ist
$\varphi^{-1} (D(g)) = X_{R(g)} $ und damit $D(g) \in \varphi(U({\bf p})) $.
\par\noindent
(3) Sei $D(g) \in \varphi(U({\bf p})) $. D.h. $ X_{R(g)} \in U({\bf p}) $, woraus
sich $R(g) \not\in {\bf p} $ aufgrund der Fixheit von ${\bf p}$ ergibt. Also
ist $g \not\in R^\ast({\bf p}) $ und $D(g) \in U(R^\ast({\bf p})) $.
\par\bigskip
Ein anderer Spezialfall des Hilbertschen Nullstellen\-satz\-es er\-gibt sich, wenn
man f"ur $g$ eine Einheit nimmt. Dann erh"alt man, da\3 eine "Uberdeckung
$X=X_1 = \bigcup_{i \in I} X_{f_i} $ impliziert, da\3 die $f_i$ das Einheitsideal
erzeugen. Um zu "aquivalenten Formulierungen zu gelangen, dient das folgende
Lemma.
\begin{lem}
Sei $X,{\cal O}_X $ ein lokal beringter Raum mit globalem Schnitt\-ring $A$,
und $\bf m $ ein maximales Ideal aus $A$.
Dann sind folgende Aus\-sa\-gen "aqui\-va\-lent:
\par\smallskip\noindent
(1) F"ur Funktionen $f_i \in {\bf m} , i \in I $, ist
$ X \neq \bigcup_{i \in I} X_{f_i} $.
\par\smallskip\noindent
(2) Das maximale Ideal $\bf m $ ist fix und konvergiert in dem Sinn, da\3
es einen Punkt $x \in X $ gibt, so da\3 aus $x \in X_f $ folgt
$ X_f \in U({\bf m}) = F(X_g : g \not\in {\bf m})$.
\par\smallskip\noindent
(3) $\bf m $ konvergiert in dem Sinn, da\3 es einen Punkt $x \in X $ gibt,
dessen funktionaler Umgebungsfilter $\{ f \in A : x \in X_f \}$ in $A-{\bf m}$
enthalten ist.
\par\smallskip\noindent
(4) Unter der Strukturabbildung $j: X \longrightarrow \mbox{ \rm Spek}A $
wird ${\bf m}$ erreicht.
\end{lem}
Beweis. Aus (1) folgt (2). Sei $g \not\in {\bf m} $. Dann gibt es
$h \in {\bf m}$ und $a \in A $ mit $ g+ a\cdot h = 1$ und damit
$ X_g \cup X_h = X $. Ist nun $X_g \subseteq X_f $, so ist ebenfalls
$X_f \cup X_h = X$, und $f$ kann nach Voraussetzung nicht zu ${\bf m}$
geh"oren. Somit ist $\bf m $ fix. W"urde $\bf m$ nicht in dem
beschriebenen Sinn konvergieren, so gebe es zu jedem Punkt $x \in X $ eine
funktionale Umgebung $ x \in X_f $ mit
$ X_f \not\in U({\bf m})$, also mit $f \in {\bf m}$,
und damit eine "Uberdeckung von $X$ mit $X_{f_i} , f_i \in {\bf m}  $, was nach
Voraussetzung nicht sein kann.
\par\noindent
Aus (2) folgt (3). Sei $x \in X$ der Konvergenzpunkt im Sinne von (2), und 
sei $ x \in X_f $. Nach Voraussetzung ist dann $X_f \supseteq X_g $ f"ur
ein $g \not\in {\bf m} $, woraus wegen der Fixheit $ f \not\in {\bf m} $
folgt.
\par\noindent
Aus (3) folgt (4). Sei $x$ der Konvergenzpunkt im Sinne von (3). Aus
$ \{ f \in A : x \in X_f \} \subseteq A-{\bf m} $ folgt
$j(x) = \{ f \in A : x \not\in X_f \} \supseteq {\bf m} $, und da $j(x)$
ein Ideal $ \neq A $ ist und $\bf m $  maximal, gilt hier die Gleichheit.
\par\noindent
Aus (4) folgt (1). Wegen ${\bf m} = j(x) = \{ f \in A : x \not\in X_f \} $
wird $X$ nicht von Invertierbarkeitsmengen zu Funktionen aus $\bf m $
"uberdeckt. 
\par\bigskip\noindent
{\bf Bemerkung} Ist $X,{\cal O}_X $ zus"atzlich reich beringt, so l"a\3t sich
die Formulierung (2) des Lemmas in die angenehmere Gestalt bringen,
da\3 $\bf m $ fix ist und echt konvergiert, also
$ U({\bf m}) \supseteq U(x) $ ist f"ur einen Punkt $ x \in X $.
\begin{satz}
F"ur einen lokal beringten Raum $X,{\cal O}_X$ mit globalem Schnittring $A$ sind
folgende Aussagen "aquivalent:
\par\smallskip\noindent
(1) In $X,{\cal O}_X $ gilt der Nullstellensatz global, d.h. aus einer
"Uberdeckung $X = \bigcup_{i \in I} X_{f_i} , f_i \in A $ folgt
$ 1 \in \langle f_i , i \in  I \rangle $.
\par\smallskip\noindent
(2) Alle maximalen Ideale ${\bf m}$ aus $A$ sind fix und konvergieren
in dem Sinn, da\3
es einen Punkt $x \in X$ gibt, so da\3 aus $x \in X_f$ folgt
$X_f \in U({\bf m})$.
\par\noindent\smallskip
(3) Alle maximalen Ideale ${\bf m}$ konvergieren in dem Sinne, da\3 es einen Punkt
$x \in X$ gibt, dessen funktionaler Umgebungsfilter $ \{f \in A: x \in X_f \}$ in
$A-{\bf m}$ enthalten ist.
\par\smallskip\noindent
(4) Unter der Strukturabbildung $ j : X \longrightarrow \mbox{ \rm Spek}A $ werden alle
maximalen Ideale erreicht.
\end{satz}
Beweis. Nach dem vorausgegangenen Lemma ist lediglich zu zeigen, da\3 der
globale Nullstellensatz "aquivalent ist zu den Formulierungen des Lemmas
f"ur alle maximalen Ideale. Ist $\bf m $ ein maximales Ideal, so be\-sagt der
Null\-stel\-len\-satz un\-mit\-tel\-bar, da\3 die $X_f , f \in {\bf m} $ den Raum nicht
"uber\-decken. Ist umgekehrt $X = \bigcup_{i \in I} X_{f_i} $, ohne da\3
diese Funktionen das Einheitsideal erzeugen, so liegen sie in einem maximalen
Ideal.
\par\bigskip\noindent
{\bf Bemerkung}  Aus dem Lemma 4.3.4 kann man nat"urlich auch eine zu Satz 
4.3.5 analoge Aussage gewinnen, in der jeweils nur
auf alle end\-lich er\-zeug\-ten ma\-xi\-ma\-len Ide\-ale
Bezug genommen wird. Diese Eigenschaft eines lo\-kal be\-ring\-ten Rau\-mes ist aber
nicht "aquivalent zu dem globalen Nullstellensatz f"ur endlich
viele Funktionen, also zu `aus $X = \bigcup_{i = 1,...,n} X_{f_i} $ folgt
$ 1 \in \langle f_1 , ... , f_n \rangle $'. In dieser Form spielt der
Nullstellensatz in der Theorie der komplexen R"aume eine wichtige Rolle. Mit
ihm lassen sich beispielsweise die Steinschen Gebiete
$ G \subseteq {\bf C}^n $ charakterisieren.
\begin{satz}
Ein Schema $X,{\cal O}_X $ ist genau dann affin, wenn es reich beringt ist und
der Nullstellensatz in der Form des vorausgegangenen Satzes gilt.
\end{satz}
Beweis. Die eine Richtung ist klar, sei also $X$ reich beringt, und es gelte
der globale Hilbertsche Nullstellensatz. Zu jedem Punkt $x \in X $ gibt es eine
offene, affine Umgebung $U_i$ und wegen reich beringt auch eine globale Funktion
$f_i$ mit $x \in X_{f_i} \subseteq U $. $X_{f_i}$ ist dann ebenfalls affin.
Somit gibt es eine "Uberdeckung $X = \bigcup_{i \in I} X_{f_i}$ mit affinen
Invertierbarkeitsmengen globaler Funktionen. Aufgrund der anderen Voraussetzung
erzeugen die $f_i$ im globalen Schnittring das Einheitsideal und damit bereits
eine endliche Auswahl davon. Damit wird $X$ bereits von endlich vielen der
$X_{f_i}$ "uberdeckt, und es liegt eine Situation vor, die es gestattet, das
bekannte Affinit"atskriterium anzuwenden. $X$ ist also affin.
\par\bigskip
Aufgrund des letzten Satzes definieren wir:
\par\bigskip\noindent
{\bf Definition} Ein lokal beringter $T_0$-Raum hei\3t affin beringt, wenn er reich
beringt ist und der globale
Hilbertsche Nullstellensatz im Sinn von Satz 4.3.5 gilt.
\par\bigskip
Ein affin beringter Raum ist stets quasikompakt. Die "Uberdeckungseigen\-schaft
braucht nur f"ur die funktionale Basis gezeigt zu werden, was man wie im
obigen Satz macht. Alle maximalen Ideale aus dem globalen Schnitt\-ring sind in
einem affin beringten Raum fix und konvergieren im eigent\-lichen Sinn, d.h. es
gibt einen Punkt $ x \in X $ mit $U(x) \subseteq U({\bf m}) $, wobei dann sogar
die Gleichheit gilt. Ferner gilt der allgemeine Hilbertsche Nullstellensatz in
der Jakobsonform: aus $ X_g \subseteq \bigcup_{i \in I} X_{f_i} $ folgt
$g \in \mbox{ Jak}( f_i : i \in I ) $. Gebe es n"amlich ein maximales Ideal
${\bf m}$ mit $f_i \in {\bf m}, g \not\in {\bf m} $, so auch einen Punkt
$x \in X$ mit $U(x) = U({\bf m})$, also $ x \in X_g $ und $ x \not\in X_{f_i} $
f"ur $ i \in I$ im Widerspruch zur vorausgesetzten Inklusion.
Abgeschlossene Mengen eines affin beringten Raumes haben eine funktionale
Umgebungsbasis. Ist n"amlich $ A \subseteq U = \bigcup_{i \in I} X_{f_i} $ und
$ V := X-A = \bigcup_{j \in J} X_{g_j} $ , so "uberdecken die $X_{f_i}$ und die
$X_{g_j} $ zusammen den Raum $X$. Damit gibt es
$ f \in \langle f_i :i \in I \rangle $
und $ g \in \langle g_j : j \in J \rangle $ mit $ f+g = 1 $. Damit ist
$X_f \subseteq U $ und $X_f \cup X_g = X $, und $ X_g \subseteq V = X-A $.
Also ist $A \subseteq X-X_g \subseteq X_f $.
\par\bigskip\noindent
{\bf Beispiel 3} Ein topologischer Raum $X$ bildet mit der Garbe der stetigen
reellwertigen Funktionen genau dann einen affin beringten Raum, wenn $X$
kompakt ist. Ein affin beringter Raum ist, wie oben bemerkt, stets quasikom\-pakt,
und reich beringt ist in dieser Situation "aquivalent zum Trennungs\-axi\-om
$T_{3a}$, was zusammen mit $T_0$ die vollst"andige Regularit"at bedeutet, aus
der bekanntlich die Hausdorffeigenschaft folgt. Ist $X$ kompakt, so ist $X$ nach
dem Lemma von Urysohn vollst"andig regul"ar, also reich beringt und $T_0$. Sei
$ X = \bigcup_{i \in I} X_{f_i} $ eine "Uberdeckung, wobei ohne Einschr"ankung
aufgrund der Kompaktheit $I$ endlich sei. Die $f_i$ haben also keine gemeinsame
Nullstelle, und damit hat die Summe ihrer Quadrate "uberhaupt keine Nullstelle
und ist daher eine Einheit im globalen Schnittring. An diesem Beispiel sieht
man auch, da\3 die Einschr"ankung eines affin beringten Raumes auf eine
Basismenge $X_f$ in der Regel nicht wieder affin beringt ist, da eine offene
Teilmenge eines kompakten Raumes nur dann wieder kompakt ist, wenn sie auch
abgeschlossen ist.
\par\smallskip\noindent
{\bf Beispiel 4} Sei $ X = {\bf R}^n $, versehen mit der Zariskitopologie
und der Garbe der rationalen Funktionen. Dann ist $X$ und jede offene Teilmenge
davon affin beringt. Jede offene Menge ist hier die Invertierbarkeitsmenge
einer Polynomfunktion, sei also $ U = X_f $ mit
$f \in {\bf R}[X_1,...X_n]=A $. Es ist $ \Gamma(U,{\cal O}_X) = A_F $, wobei $F$
aus allen Polynomen besteht, die in $X_f$ keine Nullstelle haben. F"ur
eine (endliche) "Uberdeckung $X_f = \bigcup_{i \in I} X_{g_i} $ ist dann wieder
wie im vorausgegangenen Beispiel $g = g_1^2+ ... +g_k^2 $ eine nullstellenfreie
rationale Funktion auf $X_f$, also eine Einheit.
\par\smallskip\noindent
{\bf Beispiel 5} Sei $K$ ein K"orper und
$X = \mbox{ Spek}K[S,T] - \langle S,T \rangle $. $X$ ist ein nicht affines
Schema mit dem Polynomring als globalem Schnittring. $X$ ist reich beringt,
jeder algebraische Filter ist fix, aber der zum maximalen Ideal
$\langle S,T \rangle $ geh"orende Filter konvergiert nicht.
Die Korres\-pondenz zwischen Idealen und abge\-schlos\-senen Mengen liefert also f"ur
dieses maximale Ideal die leere (Nullstellen\-)menge, w"ahrend es sich als Filter
aufgefa\3t topo\-logisch als der `Umgebungs\-filter'
eines `Loches' repr"a\-sentieren l"a\3t. In diesem Loch liegt auch die
Kohomologie. Ist ${\cal G} $ eine quasikoh"arente Garbe auf $X$, so ist
bez"uglich $ \varphi : X \longrightarrow \mbox{ Spek}K[S,T] $ bekanntlich
$ R^i \varphi_\ast {\cal G} = \widetilde{ H^i(X,{\cal G}) }$ und damit
\begin{displaymath}
H^i(X,{\cal G})_{\bf m} = \lim_{U \in U({\bf m})} H^i(U,{\cal G}) =0 
\end{displaymath}
f"ur
alle von $ \langle S,T \rangle $ verschiedenen maximalen Ideale ${\bf m}$, die
ja in $X$ konvergie\-ren, und
$H^i(X,{\cal G}) $ ist genau dann null, wenn die Kohomologie auch noch im
Lochfilter verschwindet.
\par\smallskip\noindent
{\bf Beispiel 6} Das im vorangegangenen Abschnitt diskutierte Beispiel
Spek ${\bf Z}$ mit der `fast diskreten' Topologie zeigt, da\3 f"ur ein Schema
die kanonische Strukturabbildung
$ \varphi: X,{\cal O}_X \longrightarrow \mbox{ Spek}\Gamma(X,{\cal O}_X) $
bijektiv sein kann, ohne da\3 $X$ affin ist. Es erhebt sich die Frage, ob dies
auch unter gewissen Zusatzbedingungen an $X$ der Fall sein kann, etwa wenn $X$
noethersch oder eine Variet"at ist. In der Theorie der
komplexen R"aume gibt es ein
verwandtes, offenes Problem: Folgt f"ur einen komplexen Raum $ X,{\cal O}_X $
aus der Bijektivit"at der Strukturabbildung
$\varphi : X \longrightarrow \sigma {\cal O}(X) $ , da\3 $X$ steinsch ist?
Dabei besteht $ \sigma {\cal O}(X) $ aus den stetigen
${\bf C}$-Algebramorphismen, versehen mit der Gelfandtopologie,
(\cite{kra}, S.142).
\par\smallskip\noindent
{\bf Beispiel 7} Mit einem Konvergenzargument l"a\3t sich manchmal die
Affinit"at eines Schemas einfacher als sonst beweisen. Als Beispiel beweisen
wir hier noch einmal ohne kohomologische Mittel, da\3 jede offene Teilmenge $U$
eines affinen, noetherschen, integren, eindimensionalen Schemas $X $ affin ist.
Der Beweis ist nahezu rein topologischer Art. Da $ U, {\cal O}_U $ reich
beringt ist und jedes
maximale Ideal ${\bf m} \in \mbox{ Spek}\Gamma(U,{\cal O}_U) $ fix ist, haben
wir nur zu zeigen, da\3 es in $U$
($U \neq \emptyset $)
konvergiert. Sei $F:= U({\bf m}) = F(U_f : f \not\in {\bf m} ) $ und
$ F_\ast =i(F)$ der zugeh"orige topologische Filter in $X$ bez"uglich der Einbettung
$ i : U \longrightarrow X $. $F$ ist hierbei fix.
Nach Satz 4.3.3 (2)
konvergiert $F_\ast$ in $X$. Sei $x$ ein Konvergenzpunkt.
Ist $x \in U $, so gilt die Konvergenzbeziehung auch in $U$, also
$F \supseteq U(x) $ und $F$ konvergiert. Ist $ x \not\in U $, so ist $x$ ein
abgeschlossener Punkt in $X$, und es ist $ U \subseteq X- \{ x \} $. Wegen
$ U \in F_\ast $ ist dann auch $X - \{ x \}  \in F_\ast $.
Ist $W$ eine nichtleere, offene Teilmenge von $X$,
also $W= X- \{y_1,...,y_n \} $, so ist $W \in F_\ast$: Ist $x \in W $, so folgt
dies aus der Konvergenz; ist $ x \not\in W $, so ist
$ W = X- \{x,y_2,...,y_n \} $ und $W':= X- \{y_2,...,y_n\} $ eine offene
Umgebung von $x$, und es ist $ W = W' \cap X- \{ x \} $ der Durchschnitt zweier
offener Mengen aus $F_\ast $, also $ W \in F_\ast$. $F_\ast$ und damit auch $F$
bestehen also aus allen offenen, nichtleeren
Mengen von $X$ bzw $U$. Dann konvergiert $F$ gegen den generischen Punkt von
$U$.

%% file: ausb.tex
\newpage
\section*{Ausblick}
\addcontentsline{toc}{section}{Ausblick}
Zum Schlu\3 m"ochte ich kurz andeuten, in welche Richtung die angestellten
"Uberlegungen weitergef"uhrt werden k"onnten. Neben den an einzelnen Stel\-len
im Text erw"ahnten offenen Fragen, f"ur die es vielleicht ganz ein\-fache
Ant\-wor\-ten gibt, scheint mir vor allem die Fragestellung des letzten Kapitels
nach der Korrespondenz von algebraischen und topologischen Objekten in
beringten R"aumen Ankn"upfungspunkte zu bieten. Zu untersuchen w"are etwa, welche
Eigenschaften einer bestimmten Klasse von beringten R"aumen durch ein gewisses
Korrespondenzverhalten beschrieben werden k"onnen, wobei Filter keineswegs
die Hauptrolle spielen m"ussen.
\par
Beispielsweise l"a\3t sich fragen, ob man
unter den komplexen R"aumen die Steinschen R"aume mit einem "ahnlichen Vokabular
charakterisieren kann wie die affinen Schemata durch den Begriff des affin
beringten Raumes, oder ob die formalen Schemata, die man als Vervollst"andigung
eines affinen Schemas erhalten kann, durch affin beringt erfa\3t werden
k"onnen, (\cite{har}, II.9). Allgemeiner geht es um das Problem, ob und
inwiefern das Verschwinden gewisser Kohomologiegruppen in einem lokal beringten
Raum generell mit einer engen Beziehung gewisser topologischer und
algebraischer Objekte in Zusammenhang steht. 
\par
Es ist nat"urlich ungewi\3, ob die vorliegende Arbeit zu dieser Problematik
einen Beitrag leistet, doch auch wenn dies nicht der Fall sein sollte,
glaube ich, da\3 beim Studium beringter R"aume Filter in der hier
vorgestellten
Weise miteinzubeziehen sind.

%% file: lita.tex
\addcontentsline{toc}{section}{Literaturverzeichnis}

%% file: orga.bbl
\newpage
\begin{thebibliography}{99}
\bibitem{art} Artin,M.; Grothendieck, A; Verdier, J.L.: Th\'{e}orie des topos et
co\-ho\-mo\-lo\-gie \'{e}tale des sch\'{e}\-mas, Band 1. Berlin-Heidelberg-New York 1972.
\bibitem{csa} Cs\'{a}sz\'{a}r, \'{A}.: General topology. Budapest 1978.
\bibitem{die} Dieudonn\'{e}, J.: Grundz"uge der modernen Analysis, Band III.
Braun\-schweig 1976.
\bibitem{gro} Dieudonn\'{e}, J.; Grothendieck, A.: El\'{e}ments de g\'{e}om\'{e}trie
alg\'{e}brique I. Heidelberg 1971.
\bibitem{gil} Gillman, L.; Jerison, M.: Rings of continuous functions.
Berlin-Heidelberg-New York 1976.
\bibitem{har} Hartshorne, R: Algebraic geometry. New York 1977.
\bibitem{kra} Kramm, B.: Komplexe Funktionen-Algebren. Bayreuth 1980.
\bibitem{que} von Querenburg, B.: Mengentheoretische Topologie.
Berlin-Heidelberg-New York 1979.
\bibitem{mat} Matsumura, H.: Commutative Algebra. New York 1970.
\bibitem{rei} Reiffen, H.-J.; Scheja, G.; Vetter, U.: Algebra. Mannheim 1969.
\bibitem{sch} Scheja, G.; Storch, U.: Lehrbuch der Algebra, Teil 1,2.
Stuttgart 1980,1988.
\end{thebibliography}
